\documentclass[leqno,12pt]{article}

\newlength{\textlarg}

\usepackage{amssymb, ulem}
\usepackage{euscript}
\usepackage{mathrsfs} 
\usepackage[latin1]{inputenc}
\usepackage{amsmath}
\usepackage{amsthm}
\usepackage{amsfonts}
\usepackage{array}
\usepackage{enumerate}

 \usepackage{graphicx}
 \usepackage{float} 

 \usepackage{pgf,tikz}
\usepackage{mathrsfs}
\usetikzlibrary{arrows}
\usetikzlibrary{shapes}
\usetikzlibrary{backgrounds}

\evensidemargin\oddsidemargin

\def\thesection{\arabic{section}}
\renewcommand{\theequation}{\thesection.\arabic{equation}}

\newtheorem{theorem}{Theorem}[section]
\newtheorem{lemma}[theorem]{Lemma}
\newtheorem{proposition}[theorem]{Proposition}
\newtheorem{corollary}[theorem]{Corollary}

\newtheorem{definition}[theorem]{Definition}
\theoremstyle{definition}   
\newtheorem{remark}[theorem]{Remark}

\newcommand{\eqnsection}{
\renewcommand{\theequation}{\thesection.\arabic{equation}}
    \makeatletter
    \csname  @addtoreset\endcsname{equation}{section}
    \makeatother}
\eqnsection

\def\r{{\mathbb R}}

\def\d{\mathrm{d}}
\def\G{\mathcal G}

\def\L{\mathcal L}

\def\n{\mathbf n}

\begin{document}


 \vglue30pt

\centerline{\large\bf Cluster explorations of the loop soup on a metric graph}
\centerline{\large \bf  related to the Gaussian free field }

\bigskip
\bigskip

 \centerline{by}

\medskip

 \centerline{Elie A\"{i}d\'ekon\footnote{\scriptsize LPSM, Sorbonne Universit\'e Paris VI, NYU Shanghai, and Institut Universitaire de France, {\tt elie.aidekon@upmc.fr}}}

\bigskip
\bigskip

\bigskip

{\leftskip=2truecm \rightskip=2truecm \baselineskip=15pt \small

\noindent{\slshape\bfseries Summary.} 
We consider the loop soup at intensity ${1\over 2}$ conditioned on having local time $0$ on a set of vertices with positive occupation field in their vicinities. We  give a relation between this loop soup and the usual loop soup conditioned on its local times. We deduce a domain Markov property for the loop soup, in the vein of the discrete Markov property  proved by Werner \cite{wernermarkov}: when exploring a cluster, the bridges outside the cluster form a Poisson point process. We show how it is related to the  property due to Le Jan \cite{lejan} that the local times of the loop soup are distributed as the squares of a Gaussian free field. Finally, our results naturally give the law of the loop soup conditioned on its occupation field via Fleming--Viot processes. The discrete analog of this question was addressed by Werner \cite{wernermarkov} in terms of the random current model, and by Lupu, Sabot and Tarr\`es \cite{lst19} by means of a self-interacting process.

\medskip

\noindent{\slshape\bfseries Keywords.} Loop soup, Gaussian free field, random currents, Fleming--Viot processes. 

\medskip
 
\noindent{\slshape\bfseries 2010 Mathematics Subject
Classification.} 60J80, 60K35, 82B20.

} 

\bigskip
\bigskip

\section{Introduction}
 \label{s:intro}

The Brownian loop soup  introduced by Lawler and Werner  in \cite{LW04} is defined as a Poisson point process of Brownian loops in the plane (see \cite{sheffield-werner} for its link with conformal loop ensembles). The discrete space analogs are given in Lawler and Limic \cite{LL10} in discrete time, and in Le Jan \cite{lejan} in continuous time. Le Jan showed that this last version of the loop soup is connected to the Gaussian free field: the local times of the loop soup are distributed as the squares of a Gaussian free field, which can be interpreted as a version of Dynkin's isomorphism theorem (see \cite{eisenbaum94}, \cite{ekmrs},\cite{MR06} for  other versions of Dynkin's isomorphism theorems related to Ray--Knight theorems). We refer to Lawler \cite{lawler18}, Le Jan \cite{lejan} and Powell and Werner \cite{pw20} for background on loop soups. Later, Lupu \cite{lupu16} by considering the loop soup on a metric graph,  was  able to complete the connection with the Gaussian free field by recovering the signs (and not only the squares).

\begin{figure}[h!]
\begin{tikzpicture}[ thick,main node/.style={circle,draw},x=0.5cm,y=1cm]
\tikzset{every loop/.style={min distance=2cm,in=0,out=60}}
\clip(-13.5,-2.5) rectangle (6.5,0.5);
  \node[main node] (1) at (0,0){1};
  \node[main node] (2) at (-2,-2) {2};
  \node[main node] (3) at (3,-2) {3};
  \node[main node] (4) at (0.,-1.3) {4};
  \node[main node] (5) at (3.2,-0.6) {5};
  \path[every node/.style={font=\sffamily\small}]
    (1) edge  (2)
        edge [bend right]  (2)
        edge (5)
    (2) edge  (3)
     (1)   edge  (4)
     (4) edge (3)
     (4) edge (5)
     (5) edge (3)
     (3) edge[loop] node {} (3);
\end{tikzpicture}

\begin{tikzpicture}[ thick,main node/.style={rectangle,draw},x=1.5cm,y=1.5cm]
\tikzset{every loop/.style={min distance=2cm,in=0,out=60}}
\clip(-5,-4) rectangle (5,0.5);
  \node[draw,star,star points=7] (1) at (0,0) {};
  
      \node[main node] (11) at (-1,-0.1) {1};
      \node[main node] (12) at (-0.8,-0.8) {1};
      \node[main node] (13) at (0.3,-0.6) {1};
	  \node[main node] (14) at (0.6,-0.25) {1};
  \node[draw,star,star points=7] (2) at (-3,-3) {};
      \node[main node] (21) at (-3,-2.2) {2};
      \node[main node] (22) at (-2.4,-2.4) {2};
      \node[main node] (23) at (-2.2,-3) {2};

  \node[draw,star,star points=7] (3) at (3,-3) {};	
       \node[main node] (31) at (2.2,-3) {3};
      \node[main node] (32) at (2.4,-2.6) {3};
      \node[main node] (33) at (3.2,-2.5) {3};
      \node[main node] (34) at (3.5,-3) {3};
      \node[main node] (35) at (2.8,-2.4) {3};
      
  \node[circle,draw] (4) at (0.9,-1.8) {4};
  \node[circle,draw] (5) at (2.4,-1.0) {5};
  
   \path[every node/.style={font=\sffamily\small}]
      (1) edge[dotted] node [left] {} (11)
          edge[dotted] node [left] {} (12)
          edge[dotted] node [left] {} (13)
          edge[dotted] node [left] {} (14)
      (2) edge[dotted] node [left] {} (21)
          edge[dotted] node [left] {} (22)
          edge[dotted] node [left] {} (23)
      (3) edge[dotted] node [left] {} (31)
          edge[dotted] node [left] {} (32)
          edge[dotted] node [left] {} (33)
     	  edge[dotted] node [left] {} (34)
    	  edge[dotted] node [left] {} (35)
    (11) edge[bend right]  (21)
    (12)   edge  (22)
    (13)   edge   (4)
    (4) edge (32)
    (4) edge (5)
    (5) edge (35) 
    (14) edge (5)
    (23) edge  (31)
    ;   
    
    \draw (33)   to [out=60,in=90] (4.,-2.5) to  [out=270,in=0] (34);
\end{tikzpicture}

\caption{Example of a metric graph  (top) and its associated star graph (bottom) with $W=\{1,2,3\}$. Vertices of $W$ are replaced by stars, star edges (in dotted lines) are added to the graph, old edges (in solid lines) are kept, replicas are represented by squares encapsulating the  label of the vertex they are a copy of. Unlike what the representation may suggest, old edges should have the same length in both graphs.}
\label{f:stargraph}
\end{figure}
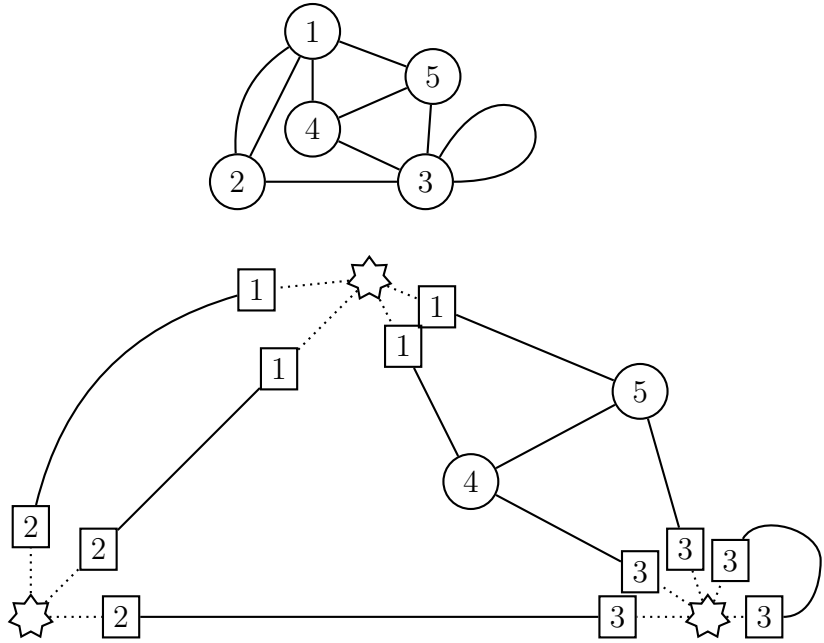

The purpose of this paper is to present  a correspondence between the loop soup on the metric graph and a conditioned loop soup on an extended graph with ``local time zero'' at vertices, and show how it sheds lights on some remarkable properties of the loop soup. The setting is the following. Let $\G$ be a connected finite metric graph. We consider the loop soup on $\G$ at intensity $\frac12$, which is a Poisson point process of Brownian loops traveling along the edges of $\G$ (see Section \ref{s:model} for a description).  Let $W$ be a set of vertices of $\G$. At each vertex  $v$ of $W$, extend each incident edge by a small amount. We then call $v$ a {\it star vertex} (represented by a star in Figure \ref{f:stargraph}) and we put at the other extremity of the added line segment a (non-star) vertex called a {\it replica} of $v$. The added line segments are called {\it star edges}, and are represented by dotted lines in Figure \ref{f:stargraph}. Notice that edges in $\G$ are still present in the extended graph, we  refer to them as {\it old edges}. We call  {\it star graph} this extended graph. See Figure \ref{f:stargraph} for an example.

In Section \ref{s:n}, we  define a probability measure $\n$  as the limiting distribution  of the loop soup on the star graph conditioned on the event that the local time $\widehat \L(v)$ at any star vertex $v\in W$ is smaller than $x_v$, and the local time fields on the star edges do not hit zero (i.e. that all star edges rooted at a star vertex are contained in the same cluster), as $(x_v,\, v\in W) \to 0$. Our result reads as follows. All trajectories are considered in this paper unoriented, unless stated otherwise.

\begin{theorem}\label{t:stargraph}
The trace of the loop soup (on $\G$) outside the vertex set $W$ conditioned on $\{\widehat \L(v)=x_v,\, v\in W\}$ has the same law as the trace of the loop soup outside  the star edges on the star graph under $\n$ conditioned on the event that for any star vertex $v\in W$, all replicas of  $v$ have local time $x_v$.
\end{theorem}

\begin{figure}[h!]
\begin{minipage}{0.5\textwidth}
\begin{tikzpicture}[line cap=round,line join=round,>=triangle 45,x=0.5cm,y=0.5cm]
\draw [xstep=1.0cm,ystep=1.0cm];
\clip(-2.5,-1.8) rectangle (12.12,6.);
\begin{scope}[shift={(7,0)}]
 \node[circle,draw,fill=black,label={[label distance=0.1cm]below:$x_v$}] (1) at (-2.4,0.3) {};
 \node[circle,draw,fill=black,label={[label distance=0.1cm]right:$x_w$}] (2) at (2.7,4.7) {};

\draw [line width=1.pt] (1)--(-2.42,1.16) to[bend left] (-1.6,2.7) to[bend left] (-0.9,2.5) to[bend left] (-1.32,2.1) to[bend left] (-1.2,3.8) to[bend left] (-0.78,4.1) to[bend left] (0.06,3.76) to[bend left]  (-0.14,3.12) to[bend left] (-0.24,4.22)to[bend left] (0.58,4.68) to (2);
\end{scope}

 \begin{scope}[rotate=200, shift={(-1.6,-1.3)}]
\draw [line width=1.pt] (-2.7,2.5)--(-1.72,2.16)--(-1.72,1.8)--(-1.74,0.52)-- (-1.8,0.42)-- (-1.84,0.36)-- (-1.88,0.3)-- (-1.92,0.24)-- (-1.94,0.3)-- (-1.88,0.34)-- (-1.82,0.36)-- (-1.72,0.38)-- (-1.62,0.38)-- (-1.54,0.38)-- (-1.46,0.38)-- (-1.36,0.38)-- (-1.28,0.38)-- (-1.16,0.38)-- (-1.1,0.36)-- (-1.02,0.36)-- (-0.92,0.32)-- (-0.86,0.3)-- (-0.8,0.28)-- (-0.74,0.24)-- (-0.66,0.2)-- (-0.6,0.16)-- (-0.54,0.12)-- (-0.48,0.04)-- (-0.42,-0.02)-- (-0.4,-0.08)-- (-0.38,-0.16)-- (-0.38,-0.26)-- (-0.38,-0.36)-- (-0.38,-0.44)-- (-0.44,-0.5)-- (-0.5,-0.56)--  (-0.66,-0.34)-- (-0.64,-0.26)-- (-0.6,-0.18)-- (-0.52,-0.12)-- (-0.46,-0.06)-- (-0.4,-0.04)--  (1.42,0.28)-- (1.52,0.34)-- (1.56,0.4)-- (1.6,0.46)-- (1.64,0.52)-- (1.66,0.58)-- (1.68,0.66)-- (1.7,0.74)-- (1.7,0.82)-- (1.64,0.88)-- (1.58,0.9)-- (1.5,0.9)-- (1.42,0.9)-- (1.38,0.84)-- (1.38,0.76)-- (1.4,0.68)-- (1.44,0.6)-- (1.48,0.54)-- (1.56,0.5)-- (1.62,0.48)-- (1.68,0.46)-- (2.64,1.46)-- (2.68,1.52)-- (2.7,1.6)-- (2.7,1.84)--  (2.,2.32)-- (1.94,2.36)-- (1.66,2.7)-- (1.64,2.76)-- (1.64,2.84)-- (1.66,2.9)-- (1.72,2.94)-- (1.78,2.98)-- (1.88,2.98)-- (1.96,2.98)-- (2.02,2.96)-- (2.04,2.9)-- (2.06,2.84)-- (2.08,2.78)-- (2.08,2.7)-- (2.08,2.6)-- (2.08,2.52)-- (2.08,2.44)-- (2.06,2.36)-- (2.02,2.3)-- (1.96,2.2)-- (1.88,2.14)-- (1.82,2.08)-- (1.76,2.06)-- (1.7,2.04)-- (1.6,2.04)-- (1.54,2.08)-- (1.48,2.12)-- (1.42,2.18)-- (1.34,2.26)-- (1.3,2.32)-- (1.26,2.4)-- (1.2,2.46)-- (1.14,2.52)-- (1.08,2.58)-- (1.02,2.62)-- (0.94,2.66)-- (0.88,2.68)-- (0.8,2.68)-- (0.72,2.68)-- (0.64,2.68)-- (0.54,2.66)-- (0.46,2.58)-- (0.38,2.5)-- (0.32,2.4)-- (0.28,2.34)-- (0.24,2.26)-- (0.2,2.2)-- (0.16,2.12)-- (0.14,2.06)-- (0.14,1.98)-- (0.14,1.9)-- (0.22,1.84)-- (0.28,1.82)-- (0.34,1.8)-- (0.42,1.8)-- (0.48,1.84)-- (0.52,1.9)-- (0.54,1.96)-- (0.54,2.04)-- (0.52,2.1)-- (0.46,2.14)-- (0.4,2.16)-- (0.32,2.16)-- (0.24,2.16)--(-1.72,2.16);
\end{scope}

\begin{scope}[shift={(8,0)}]
\draw [line width=1.pt] (-3.4,0.3) -- (-2.14,0.3);
\draw [line width=1.pt] (1.68,4.6)--(1.68,3.96)-- (1.68,3.96)-- (1.68,3.88)-- (1.68,3.8)-- (1.7,3.72)-- (1.74,3.66)-- (1.8,3.6)-- (1.88,3.6)-- (1.94,3.58)-- (2.04,3.58)-- (2.1,3.6)-- (2.16,3.62)-- (2.18,3.68)-- (2.18,3.76)-- (2.12,3.8)-- (2.04,3.8)-- (1.96,3.8)-- (1.9,3.74)-- (1.84,3.68)-- (1.82,3.62)-- (1.82,3.54)-- (1.82,3.46)-- (1.82,3.36)-- (1.82,3.26)-- (1.82,3.16)-- (1.82,3.06)-- (1.82,2.94)-- (1.82,2.86)-- (1.82,2.78)-- (1.82,2.68)-- (1.82,2.58)-- (1.82,2.5)-- (1.78,2.44)-- (1.74,2.38)-- (1.68,2.34)-- (1.62,2.3)-- (1.6,2.38)-- (1.64,2.44)-- (1.72,2.5)-- (1.78,2.56)-- (1.84,2.6)-- (1.92,2.64)-- (2.,2.64) -- (2.08,2.64)-- (2.14,2.6)-- (2.18,2.54)-- (2.2,2.48)-- (2.22,2.4)-- (2.26,2.32)-- (2.28,2.24)-- (2.28,2.12)-- (2.28,2.04)-- (2.28,1.92)-- (2.28,1.82)-- (2.28,1.74)-- (1.4,1.72)-- (1.3,1.72)-- (1.2,1.72)-- (1.12,1.72)-- (1.04,1.72)-- (0.98,1.7)-- (0.92,1.68)-- (0.86,1.64)-- (0.8,1.6)-- (0.78,1.54)-- (0.78,1.46)-- (0.8,1.4)-- (0.88,1.42)-- (0.96,1.48)-- (1.02,1.54)-- (1.06,1.6)-- (1.08,1.68)-- (1.1,1.78)-- (1.1,1.9)-- (1.1,1.98)-- (1.1,2.06)-- (1.1,2.14)-- (1.1,2.22)-- (1.08,2.28)-- (1.06,2.36)-- (1.04,2.42)-- (0.98,2.5)-- (0.92,2.58)-- (0.86,2.64)-- (0.78,2.7)-- (0.72,2.74)-- (0.66,2.78)-- (0.6,2.82)-- (0.54,2.84)-- (0.46,2.84)-- (0.38,2.82)-- (0.32,2.8)-- (0.24,2.76)-- (0.14,2.74)-- (0.06,2.74)-- (-0.04,2.74)-- (-0.12,2.74)-- (-0.18,2.76)-- (-0.24,2.8)-- (-0.28,2.86)-- (-0.3,2.92)-- (-0.3,3.)-- (-0.28,3.06)-- (-0.24,3.12)-- (-0.18,3.16)-- (-0.12,3.18)-- (-0.06,3.2)-- (0.02,3.2)-- (0.1,3.2)-- (0.2,3.14)-- (0.26,3.08)-- (0.3,3.02)-- (0.34,2.94)-- (0.36,2.88)-- (0.38,2.82)-- (0.38,2.74)-- (0.4,2.66)-- (0.4,2.56)-- (0.4,2.44)-- (0.38,2.36)-- (0.34,2.28)-- (0.28,2.18)-- (0.24,2.1)-- (0.1,1.84)-- (0.08,1.78)-- (0.02,1.7)-- (-0.04,1.64)-- (-0.12,1.58)-- (-0.18,1.56)-- (-0.22,1.62)-- (-0.16,1.64)-- (-0.08,1.64)-- (0.,1.64)-- (0.08,1.64)-- (0.14,1.6)-- (0.2,1.56)-- (0.26,1.5)-- (0.3,1.42)-- (0.34,1.36)-- (0.36,1.28)-- (0.38,1.2)-- (0.4,1.12)-- (0.42,1.04)-- (0.42,0.96)-- (0.42,0.88)-- (0.42,0.8)-- (0.42,0.72)-- (0.38,0.62)-- (0.34,0.56)-- (0.3,0.5)-- (0.22,0.44)-- (0.14,0.36)-- (0.04,0.32)-- (-0.06,0.3)-- (-0.16,0.3)-- (-0.26,0.3)-- (-0.34,0.34)-- (-0.4,0.38)-- (-0.46,0.44)-- (-0.5,0.5)-- (-0.54,0.56)-- (-0.6,0.6)-- (-0.66,0.64)-- (-0.72,0.68)-- (-0.78,0.7)-- (-0.86,0.72)-- (-0.96,0.72)-- (-1.04,0.72)-- (-1.12,0.72)-- (-1.18,0.66)-- (-1.22,0.6)-- (-1.26,0.52)-- (-1.3,0.4)-- (-1.34,0.32)-- (-1.36,0.26)-- (-1.36,0.16)-- (-1.4,0.08)-- (-1.4,-0.04)-- (-1.4,-0.12)-- (-1.38,-0.18)-- (-1.32,-0.2)-- (-1.24,-0.22)-- (-1.18,-0.24)-- (-1.12,-0.22)-- (-1.08,-0.16)-- (-1.06,-0.1)-- (-1.06,-0.02)-- (-1.06,0.06)-- (-1.12,0.12)-- (-1.18,0.18)-- (-1.24,0.22)-- (-1.3,0.24)-- (-1.38,0.26)-- (-1.46,0.26)-- (-1.54,0.26)-- (-1.64,0.26)-- (-1.72,0.26)-- (-1.8,0.26)-- (-2.14,0.3);
\end{scope}
\end{tikzpicture}
\end{minipage}
\begin{minipage}{0.5\textwidth}
\begin{tikzpicture}[line cap=round,line join=round,>=triangle 45,x=0.5cm,y=0.5cm]
\clip(-2.5,-3.08) rectangle (13.12,7.7);

\node[draw,star, star points=7,fill=black,label={[label distance=0.1cm]below:$0$}]  at (4.6,0.3) {};
\node[draw, star, star points=7, fill=black,label={[label distance=0.1cm]45:$0$}]  at (10.7,5.7) {};
\draw[dotted] (4.6,0.3)--(3.5,0.3);
\draw[dotted] (4.6,0.3)--(4.6,1.3);
\draw[dotted] (4.6,0.3)--(5.6,0.3);
\draw[dotted] (10.7,5.7)--(9.5,5.7);
\draw[dotted] (10.7,5.7)--(10.7,4.7);

\begin{scope}[shift={(7,1)}]
 \node[rectangle,draw,fill=black,label={[label distance=0.1cm]left:$x_v$}] (1) at (-2.4,0.3) {};
 \node[rectangle,draw,fill=black,label={[label distance=0.1cm]above:$x_w$}] (2) at (2.7,4.7) {};

\draw [line width=1.pt] (1)--(-2.42,1.16) to[bend left] (-1.6,2.7) to[bend left] (-0.9,2.5) to[bend left] (-1.32,2.1) to[bend left] (-1.2,3.8) to[bend left] (-0.78,4.1) to[bend left] (0.06,3.76) to[bend left]  (-0.14,3.12) to[bend left] (-0.24,4.22)to[bend left] (0.58,4.68) to (2);
\end{scope}

 \begin{scope}[rotate=200, shift={(-0.7,-1.58)}]
\draw node[rectangle,fill=black,label={[label distance=0.1cm]below:$x_v$}] at (-2.7,2.5){};
\draw [line width=1.pt] (-2.7,2.5)--(-1.72,2.16)--(-1.72,1.8)--(-1.74,0.52)-- (-1.8,0.42)-- (-1.84,0.36)-- (-1.88,0.3)-- (-1.92,0.24)-- (-1.94,0.3)-- (-1.88,0.34)-- (-1.82,0.36)-- (-1.72,0.38)-- (-1.62,0.38)-- (-1.54,0.38)-- (-1.46,0.38)-- (-1.36,0.38)-- (-1.28,0.38)-- (-1.16,0.38)-- (-1.1,0.36)-- (-1.02,0.36)-- (-0.92,0.32)-- (-0.86,0.3)-- (-0.8,0.28)-- (-0.74,0.24)-- (-0.66,0.2)-- (-0.6,0.16)-- (-0.54,0.12)-- (-0.48,0.04)-- (-0.42,-0.02)-- (-0.4,-0.08)-- (-0.38,-0.16)-- (-0.38,-0.26)-- (-0.38,-0.36)-- (-0.38,-0.44)-- (-0.44,-0.5)-- (-0.5,-0.56)--  (-0.66,-0.34)-- (-0.64,-0.26)-- (-0.6,-0.18)-- (-0.52,-0.12)-- (-0.46,-0.06)-- (-0.4,-0.04)--  (1.42,0.28)-- (1.52,0.34)-- (1.56,0.4)-- (1.6,0.46)-- (1.64,0.52)-- (1.66,0.58)-- (1.68,0.66)-- (1.7,0.74)-- (1.7,0.82)-- (1.64,0.88)-- (1.58,0.9)-- (1.5,0.9)-- (1.42,0.9)-- (1.38,0.84)-- (1.38,0.76)-- (1.4,0.68)-- (1.44,0.6)-- (1.48,0.54)-- (1.56,0.5)-- (1.62,0.48)-- (1.68,0.46)-- (2.64,1.46)-- (2.68,1.52)-- (2.7,1.6)-- (2.7,1.84)--  (2.,2.32)-- (1.94,2.36)-- (1.66,2.7)-- (1.64,2.76)-- (1.64,2.84)-- (1.66,2.9)-- (1.72,2.94)-- (1.78,2.98)-- (1.88,2.98)-- (1.96,2.98)-- (2.02,2.96)-- (2.04,2.9)-- (2.06,2.84)-- (2.08,2.78)-- (2.08,2.7)-- (2.08,2.6)-- (2.08,2.52)-- (2.08,2.44)-- (2.06,2.36)-- (2.02,2.3)-- (1.96,2.2)-- (1.88,2.14)-- (1.82,2.08)-- (1.76,2.06)-- (1.7,2.04)-- (1.6,2.04)-- (1.54,2.08)-- (1.48,2.12)-- (1.42,2.18)-- (1.34,2.26)-- (1.3,2.32)-- (1.26,2.4)-- (1.2,2.46)-- (1.14,2.52)-- (1.08,2.58)-- (1.02,2.62)-- (0.94,2.66)-- (0.88,2.68)-- (0.8,2.68)-- (0.72,2.68)-- (0.64,2.68)-- (0.54,2.66)-- (0.46,2.58)-- (0.38,2.5)-- (0.32,2.4)-- (0.28,2.34)-- (0.24,2.26)-- (0.2,2.2)-- (0.16,2.12)-- (0.14,2.06)-- (0.14,1.98)-- (0.14,1.9)-- (0.22,1.84)-- (0.28,1.82)-- (0.34,1.8)-- (0.42,1.8)-- (0.48,1.84)-- (0.52,1.9)-- (0.54,1.96)-- (0.54,2.04)-- (0.52,2.1)-- (0.46,2.14)-- (0.4,2.16)-- (0.32,2.16)-- (0.24,2.16)--(-1.72,2.16);
\end{scope}

\begin{scope}[shift={(9,0)}]
\node[rectangle,draw,fill=black,label={[label distance=0.1cm]below:$x_v$}] at (-3.4,0.3){};
\node[rectangle,draw,fill=black,label={[label distance=0.1cm]right:$x_w$}] at (1.68,4.6){};
\draw [line width=1.pt] (-3.4,0.3) -- (-2.14,0.3);
\draw [line width=1.pt] (1.68,4.6)--(1.68,3.96)-- (1.68,3.96)-- (1.68,3.88)-- (1.68,3.8)-- (1.7,3.72)-- (1.74,3.66)-- (1.8,3.6)-- (1.88,3.6)-- (1.94,3.58)-- (2.04,3.58)-- (2.1,3.6)-- (2.16,3.62)-- (2.18,3.68)-- (2.18,3.76)-- (2.12,3.8)-- (2.04,3.8)-- (1.96,3.8)-- (1.9,3.74)-- (1.84,3.68)-- (1.82,3.62)-- (1.82,3.54)-- (1.82,3.46)-- (1.82,3.36)-- (1.82,3.26)-- (1.82,3.16)-- (1.82,3.06)-- (1.82,2.94)-- (1.82,2.86)-- (1.82,2.78)-- (1.82,2.68)-- (1.82,2.58)-- (1.82,2.5)-- (1.78,2.44)-- (1.74,2.38)-- (1.68,2.34)-- (1.62,2.3)-- (1.6,2.38)-- (1.64,2.44)-- (1.72,2.5)-- (1.78,2.56)-- (1.84,2.6)-- (1.92,2.64)-- (2.,2.64) -- (2.08,2.64)-- (2.14,2.6)-- (2.18,2.54)-- (2.2,2.48)-- (2.22,2.4)-- (2.26,2.32)-- (2.28,2.24)-- (2.28,2.12)-- (2.28,2.04)-- (2.28,1.92)-- (2.28,1.82)-- (2.28,1.74)-- (1.4,1.72)-- (1.3,1.72)-- (1.2,1.72)-- (1.12,1.72)-- (1.04,1.72)-- (0.98,1.7)-- (0.92,1.68)-- (0.86,1.64)-- (0.8,1.6)-- (0.78,1.54)-- (0.78,1.46)-- (0.8,1.4)-- (0.88,1.42)-- (0.96,1.48)-- (1.02,1.54)-- (1.06,1.6)-- (1.08,1.68)-- (1.1,1.78)-- (1.1,1.9)-- (1.1,1.98)-- (1.1,2.06)-- (1.1,2.14)-- (1.1,2.22)-- (1.08,2.28)-- (1.06,2.36)-- (1.04,2.42)-- (0.98,2.5)-- (0.92,2.58)-- (0.86,2.64)-- (0.78,2.7)-- (0.72,2.74)-- (0.66,2.78)-- (0.6,2.82)-- (0.54,2.84)-- (0.46,2.84)-- (0.38,2.82)-- (0.32,2.8)-- (0.24,2.76)-- (0.14,2.74)-- (0.06,2.74)-- (-0.04,2.74)-- (-0.12,2.74)-- (-0.18,2.76)-- (-0.24,2.8)-- (-0.28,2.86)-- (-0.3,2.92)-- (-0.3,3.)-- (-0.28,3.06)-- (-0.24,3.12)-- (-0.18,3.16)-- (-0.12,3.18)-- (-0.06,3.2)-- (0.02,3.2)-- (0.1,3.2)-- (0.2,3.14)-- (0.26,3.08)-- (0.3,3.02)-- (0.34,2.94)-- (0.36,2.88)-- (0.38,2.82)-- (0.38,2.74)-- (0.4,2.66)-- (0.4,2.56)-- (0.4,2.44)-- (0.38,2.36)-- (0.34,2.28)-- (0.28,2.18)-- (0.24,2.1)-- (0.1,1.84)-- (0.08,1.78)-- (0.02,1.7)-- (-0.04,1.64)-- (-0.12,1.58)-- (-0.18,1.56)-- (-0.22,1.62)-- (-0.16,1.64)-- (-0.08,1.64)-- (0.,1.64)-- (0.08,1.64)-- (0.14,1.6)-- (0.2,1.56)-- (0.26,1.5)-- (0.3,1.42)-- (0.34,1.36)-- (0.36,1.28)-- (0.38,1.2)-- (0.4,1.12)-- (0.42,1.04)-- (0.42,0.96)-- (0.42,0.88)-- (0.42,0.8)-- (0.42,0.72)-- (0.38,0.62)-- (0.34,0.56)-- (0.3,0.5)-- (0.22,0.44)-- (0.14,0.36)-- (0.04,0.32)-- (-0.06,0.3)-- (-0.16,0.3)-- (-0.26,0.3)-- (-0.34,0.34)-- (-0.4,0.38)-- (-0.46,0.44)-- (-0.5,0.5)-- (-0.54,0.56)-- (-0.6,0.6)-- (-0.66,0.64)-- (-0.72,0.68)-- (-0.78,0.7)-- (-0.86,0.72)-- (-0.96,0.72)-- (-1.04,0.72)-- (-1.12,0.72)-- (-1.18,0.66)-- (-1.22,0.6)-- (-1.26,0.52)-- (-1.3,0.4)-- (-1.34,0.32)-- (-1.36,0.26)-- (-1.36,0.16)-- (-1.4,0.08)-- (-1.4,-0.04)-- (-1.4,-0.12)-- (-1.38,-0.18)-- (-1.32,-0.2)-- (-1.24,-0.22)-- (-1.18,-0.24)-- (-1.12,-0.22)-- (-1.08,-0.16)-- (-1.06,-0.1)-- (-1.06,-0.02)-- (-1.06,0.06)-- (-1.12,0.12)-- (-1.18,0.18)-- (-1.24,0.22)-- (-1.3,0.24)-- (-1.38,0.26)-- (-1.46,0.26)-- (-1.54,0.26)-- (-1.64,0.26)-- (-1.72,0.26)-- (-1.8,0.26)-- (-2.14,0.3);
\end{scope}
\end{tikzpicture}
\end{minipage}
\caption{Illustration of Theorem \ref{t:stargraph}. We condition the loop soup on having local times $x_v$ and $x_w$ at vertices $v$ and $w$ (left). Excursions made by loops which hit these two vertices have the same law as excursions away from the replica (right) on the star graph under the measure $\n$, when we condition the replicas on having the local time of the original vertex, and the star edges on having positive local time fields (except at their root, the local time at a star vertex is $0$ under $\n$ by construction).}
\label{f:loopsoupcond}
\end{figure}
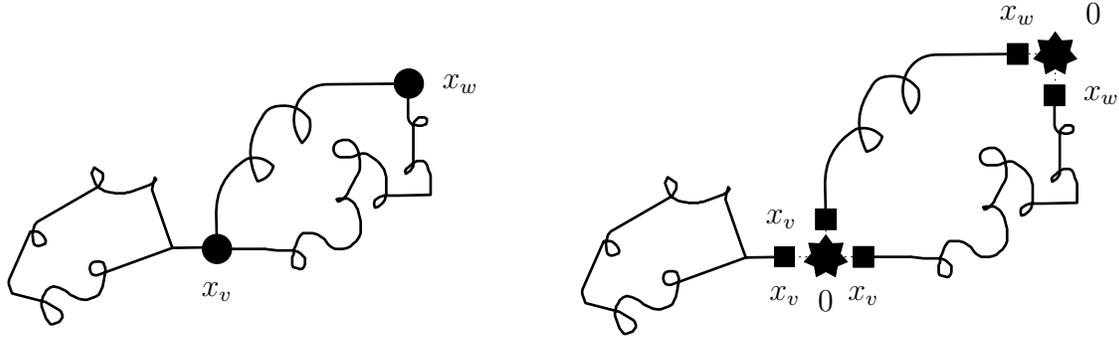

Let us specify the statement of the theorem. The trace of the loop soup outside the vertex set $W$ means that we look at loops which do not hit any vertex of $W$, and excursions between vertices of $W$.  We insist that the trajectories are unoriented. Similarly, the loop soup outside the star edges is made of loops which do not hit the star edges, and excursions between replicas which lie outside the star edges, again considered unoriented. We do not keep the information on whether the excursions belong to the same loop or not. See Figure \ref{f:loopsoupcond} for an illustration. We draw some consequences of Theorem \ref{t:stargraph}.


\paragraph{Spatial Markov property}
 
A subgraph $\widetilde \G$ of $\G$ is the data of a set of vertices  and of edges of $\G$ such that any edge of $\widetilde{\G}$ has its endpoints in $\widetilde{\G}$. When a vertex of $\widetilde{\G}$ is adjacent in $\G$ to an edge which is not in $\widetilde{\G}$, we say that the vertex lies on the boundary of $\widetilde{\G}$.  The trace of the loop soup outside $\widetilde{\G}$ is the collection of loops which do not hit $\widetilde{\G}$, and of excursions away from vertices on the boundary of $\widetilde{\G}$ which do not hit $\widetilde{\G}$ except at their endpoints. Again, all trajectories are considered unoriented, and we lose track on whether two excursions belong to the same loop or not.

\begin{theorem}\label{t:markov}
Let $\widetilde{\G}$ be a subgraph of $\G$.  Conditioned  on the local times at vertices which lie on the boundary of $\widetilde{\G}$ and on the event that the local times inside the edges of $\widetilde{\G}$ do not hit zero,  the trace of the loop soup outside $\widetilde{\G}$ is independent of the occupation field inside $\widetilde{\G}$.
\end{theorem}

 If  $\widetilde{\G}$ is connected, we have the following  simple description of the conditional law. Write $\partial \widetilde{\G}$ for the set of vertices which lie on the boundary of $\widetilde \G$.

\begin{theorem}\label{t:markov2}(Domain Markov property)
Let $\widetilde{\G}$ be a {\it connected} subgraph of $\G$.  Let $H_{\widetilde {\G}^c}$ be the boundary Poisson kernel for the graph outside $\widetilde{\G}$: for vertices $v,w \in \partial \widetilde{\G}$, $H_{\widetilde {\G}^c}(v,w)$ is the mass under the excursion measure at $v$ of the excursions from $v$ which stay away from $\widetilde {\G}$ before hitting $w$. Conditioned  on the local times at vertices $v\in \partial \widetilde{\G}$ being $(x_v,\, v\in \partial \widetilde{\G})$ and on the event that the local times inside the edges of $\widetilde{\G}$ do not hit zero, the numbers of excursions between pairs of vertices $v\neq w \in  \partial \widetilde{\G}$ are independent Poisson distributed random variables with respective parameter $2 H_{\widetilde {\G}^c}(v,w) \sqrt{x_v x_w}$.
\end{theorem}

These theorems  are closely related to a result of Werner \cite{wernermarkov}  which states that when conditioning on the local times on a set $W$ of vertices, the law of the trace of the  loop soup outside $W$ is given by a loop soup (for loops which do not hit $W$) plus a Poisson point process of bridges with  parameters $2 H_{\widetilde {\G}^c}(v,w) \sqrt{x_v x_w}$, conditioned on having an even number of bridges at each vertex of $W$. This last condition may give in our setting a dependence between the outside of $\widetilde{\G}$ and the inside. This dependence is only a matter of parity. From this point of view, Theorem \ref{t:markov}  says  that conditioned on the local times on the boundary of $\widetilde{\G}$, the configuration of even/odd crossings inside $\widetilde{\G}$ is actually independent of the occupation field, when we suppose that the occupation field stays nonzero.  We mention that we cannot remove the nonzero condition for the local times (when the graph has cycles). When looking for example at a loop soup on a circle,  the presence of a zero in the local time field inside some explored arc will indicate the absence of cycles so that the number of bridges in the unexplored arc is necessarily even. This information cannot be retrieved by looking solely at the local times at the boundary of the arc. This issue disappears when looking at the edge-occupation fields of the discrete loop soup. In fact, Werner \cite{wernermarkov} (see \cite{CaLis} for non-backtracking loop soups) shows a general spatial Markov property in the discrete case (as well as in the case of the loop soup at intensity $1$): trajectories inside and outside a domain are independent conditionally on the number of crossings of the edges on the boundary. This statement is not true when we rather condition on the local times at vertices on the boundary.  Theorem \ref{t:markov2} is also reminiscent of the decomposition of a cluster of the Brownian loop soup in the plane proved by Qian and Werner \cite{qw19}, \cite{qw18}: while the outer boundary is a loop of a ${\rm CLE}_4$, the excursions of the cluster which touch the boundary form a Poisson point process. See the work of Aru, Lupu and Sep\'ulveda \cite{als20} for a description of this decomposition in terms of first passage sets.

We prove our theorems  using the description of the loop soup given in Theorem \ref{t:stargraph}, when taking for $W$ the set of vertices of $\widetilde{\G}$. The following observation is at the heart of the proof. Under $\n$, loops are only allowed to cross at most once  edges adjacent to $W$. When an edge is crossed, its occupation field is a ${\rm BESQ}^3$ bridge (squared Bessel bridge of dimension $3$, see chap. XI in \cite{revuz-yor}). Indeed the loop soup naturally gives a ${\rm BESQ}^1$ bridge, and the crossing adds a ${\rm BESQ}^2$ bridge. When an edge is not crossed, we are asking the ${\rm BESQ}^1$ bridge from the loop soup to stay positive. This process is again a ${\rm BESQ}^3$ bridge. We conclude that the occupation field will not give any information on the crossing of an edge under $\n$.


\paragraph{Le Jan's isomorphism theorem}

As mentioned earlier, Le Jan's isomorphism theorem \cite{lejan} states that  the occupation times of a discrete  loop soup at intensity ${1\over 2}$ on a graph are distributed as the squares of a Gaussian free field (GFF) (see \cite{FR14} for a general setting).  Lupu \cite{lupu16} gave a signed version of this isomorphism, by assigning signs on the vertices so as to recover the Gaussian free field itself instead of its square. This signed isomorphism is obtained  by considering the  loop  soup on the associated metric graph. Signs are then chosen independently on each cluster of the loop soup.  The Markovian properties described above reflect those of the Gaussian free field. 

We show that we can go in the reverse direction, and see Le Jan's isomorphism theorem as a consequence of the domain Markov property of the loop soup stated in Theorem \ref{t:markov2}. Note that the original proof by Le Jan \cite{lejan} of the isomorphism theorem via the Feynman--Kac formula is direct and more general. The same claim holds for the proof which can be found in the lecture notes of Powell and Werner \cite{pw20}. Our goal is to see how one-dimensional arguments can be transferred to the metric graph. It also somehow makes clear the link between the Markov property of the GFF and that of the loop soup: when exploring a cluster, the square of the harmonic extension of the GFF at a vertex $v$ of the unexplored part corresponds for the loop soup to the mean number of excursions from the boundary which go through $v$, multiplied by the Green function at that point. This observation is related to the isomorphism theorem proved by Aru, Lupu and Sep\'ulveda in \cite{als20}.

On the real line, the generalized  Ray--Knight theorems say that the local times of a loop soup are {\rm BESQ} processes (\cite{legall-yor},\cite{carmona-petit-yor94b},\cite{werner95},\cite{eyzburglar}). In the setting of the metric graph, we will do the analog and locally explore the edges of a cluster of the loop soup, write the associated SDE for the local time via one-dimensional Ray--Knight theorems and verify that it is indeed that followed by the GFF, described in Lupu and Werner \cite{lw18}. In the GFF setting, various explorations have been studied in the literature, see for example  \cite{schramm-sheffield}, \cite{lw18} and \cite{als20}.

\paragraph{Random current model}

We give a  surprisingly simple description of the loop soup conditioned on the clusters.

For $\ell_1,\ell_2,\rho>0$, call ${\mathcal B}(\rho,\ell_1,\ell_2)$ the Brownian loop soup with   intensity  $\frac12$,  on the interval $[0,\rho]$, conditioned on having local times $\ell_1$ and $\ell_2$ at the endpoints of the interval. Observe that the number of crossings of the interval is necessarily even. Call ${\mathcal C}(\rho,\ell_1,\ell_2)$ the analog when we add an extra crossing so that the number of crossings is now odd. A quick way to define it is to take a Brownian loop soup with intensity $\frac12+1$, and concatenate the extra loops inside $[0,\rho]$  to form the extra crossing. We refer to Section \ref{s:loopsoupcond} to specify the definitions of these objects. We call cluster a maximal set of adjacent edges with positive local time field. 
A configuration $(\alpha_e)_e$  of integers in $\{0,1\}$ on the edges of $\G$ is said admissible if :   when $e$ is not contained in a cluster, necessarily $\alpha_e=0$;  for every vertex $v$, the sum of the $\alpha_e$'s over all incoming edges $e$ of $v$ is even (a self-loop is counted twice). For an edge $e$, we call $\rho(e)$ the length of the edge, $e_1$ and $e_2$ its endpoints, and $\widehat {\mathcal L}(e_1)$, resp. $\widehat {\mathcal L}(e_2)$, denotes the local time at $e_1$, resp. $e_2$.

\begin{theorem}\label{t:cluster}
 Consider the loop soup at intensity $\frac12$ on the metric graph $\G$. Set $\alpha_e=1$ if the number of crossings of the edge $e$ is odd, and $\alpha_e=0$ if it is even. 
\begin{enumerate}[(1)]
\item Conditionally on the clusters of the loop soup, the configuration $(\alpha_e)_e$ is uniform among all admissible configurations and is independent of the occupation field of the loop soup.
\item Conditionally on $(\alpha_e)_e$  and on the local times at vertices, the traces of the loop soup on the edges are independent:
\begin{enumerate}
\item if $\alpha_e=0$,  the trace on the edge $e$ is distributed as ${\mathcal B}(\rho(e),\widehat \L(e_1),\widehat \L(e_2))$; 
 \item if $\alpha_e=1$, the trace on the edge $e$ is distributed as  ${\mathcal C}(\rho(e),\widehat \L(e_1),\widehat \L(e_2))$.
\end{enumerate}
\end{enumerate}

\end{theorem}

\begin{figure}[h!]
\begin{tikzpicture}[line cap=round,line join=round,>=triangle 45,x=0.5cm,y=.5cm]
\clip(-10.,-7) rectangle (24,7);
\draw [line width=3.6pt] (-2.,6.)-- (-2.,4.);
\draw [line width=3.6pt] (-2.,6.)-- (0.,6.);
\draw [line width=3.6pt] (-2.,4.)-- (0.,4.);
\draw [line width=3.6pt] (0.,6.)-- (0.,4.);
\draw [line width=2.8pt,dash pattern=on 7pt off 7pt] (0.,4.)-- (2.,4.);
\draw [line width=2.8pt,dash pattern=on 7pt off 7pt] (2.,4.)-- (2.,6.);
\draw [line width=2.8pt,dash pattern=on 7pt off 7pt] (2.,6.)-- (4.,6.);
\draw [line width=2.8pt,dash pattern=on 7pt off 7pt] (4.,6.)-- (4.,4.);
\draw [line width=3.6pt] (2.,4.)-- (4.,4.);
\draw [line width=2.8pt,dash pattern=on 7pt off 7pt] (4.,6.)-- (6.,6.);
\draw [line width=2.8pt,dash pattern=on 7pt off 7pt] (4.,4.)-- (6.,4.);
\draw [line width=3.6pt] (8.,4.)-- (8.,6.);
\draw [line width=2.8pt,dash pattern=on 7pt off 7pt] (8.,2.)-- (6.,2.);
\draw [line width=3.6pt] (4.,4.)-- (4.,2.);
\draw [line width=3.6pt] (4.,2.)-- (2.,2.);
\draw [line width=3.6pt] (2.,4.)-- (2.,2.);
\draw [line width=3.6pt] (2.,2.)-- (0.,2.);
\draw [line width=2.8pt,dash pattern=on 7pt off 7pt] (0.,4.)-- (0.,2.);
\draw [line width=2.8pt,dash pattern=on 7pt off 7pt] (0.,2.)-- (-2.,2.);
\draw [line width=2.8pt,dash pattern=on 7pt off 7pt] (-2.,4.)-- (-2.,2.);
\draw [line width=2.8pt,dash pattern=on 7pt off 7pt] (-2.,2.)-- (-2.,0.);
\draw [line width=2.8pt,dash pattern=on 7pt off 7pt] (-2.,0.)-- (0.,0.);
\draw [line width=3.6pt] (0.,0.)-- (0.,2.);
\draw [line width=3.6pt] (0.,0.)-- (2.,0.);
\draw [line width=3.6pt] (2.,0.)-- (2.,2.);
\draw [line width=2.8pt,dash pattern=on 7pt off 7pt] (2.,0.)-- (4.,0.);
\draw [line width=2.8pt,dash pattern=on 7pt off 7pt] (6.,0.)-- (6.,2.);
\draw [line width=3.6pt] (6.,0.)-- (8.,0.);
\draw [line width=2.8pt,dash pattern=on 7pt off 7pt] (8.,0.)-- (8.,2.);
\draw [line width=2.8pt,dash pattern=on 7pt off 7pt] (8.,0.)-- (8.,-2.);
\draw [line width=3.6pt] (8.,-2.)-- (6.,-2.);
\draw [line width=3.6pt] (6.,-2.)-- (6.,0.);
\draw [line width=2.8pt,dash pattern=on 7pt off 7pt] (6.,-2.)-- (4.,-2.);
\draw [line width=2.8pt,dash pattern=on 7pt off 7pt] (2.,-2.)-- (2.,0.);
\draw [line width=3.6pt] (8.,6.)-- (10.,6.);
\draw [line width=2.8pt,dash pattern=on 7pt off 7pt] (10.,6.)-- (10.,4.);
\draw [line width=3.6pt] (8.,4.)-- (10.,4.);
\draw [line width=2.8pt,dash pattern=on 7pt off 7pt] (8.,2.)-- (10.,2.);
\draw [line width=2.8pt,dash pattern=on 7pt off 7pt] (10.,4.)-- (10.,2.);
\draw [line width=3.6pt] (8.,0.)-- (10.,0.);
\draw [line width=2.8pt,dash pattern=on 7pt off 7pt] (10.,2.)-- (10.,0.);
\draw [line width=2.8pt,dash pattern=on 7pt off 7pt] (10.,0.)-- (10.,-2.);
\draw [line width=2.8pt,dash pattern=on 7pt off 7pt] (10.,-2.)-- (10.,-4.);
\draw [line width=3.6pt] (10.,-4.)-- (8.,-4.);
\draw [line width=3.6pt] (8.,-2.)-- (8.,-4.);
\draw [line width=2.8pt,dash pattern=on 7pt off 7pt] (4.,-2.)-- (4.,-4.);
\draw [line width=2.8pt,dash pattern=on 7pt off 7pt] (2.,-2.)-- (2.,-4.);
\draw [line width=2.8pt,dash pattern=on 7pt off 7pt] (0.,-2.)-- (0.,-4.);
\draw [line width=3.6pt] (0.,-4.)-- (-2.,-4.);
\draw [line width=2.8pt,dash pattern=on 7pt off 7pt] (-2.,-2.)-- (-2.,-4.);
\draw [line width=3.6pt] (10.,6.)-- (12.,6.);
\draw [line width=3.6pt] (12.,6.)-- (12.,4.);
\draw [line width=3.6pt] (10.,4.)-- (12.,4.);
\draw [line width=2.8pt,dash pattern=on 7pt off 7pt] (10.,2.)-- (12.,2.);
\draw [line width=2.8pt,dash pattern=on 7pt off 7pt] (12.,2.)-- (12.,0.);
\draw [line width=3.6pt] (10.,0.)-- (12.,0.);
\draw [line width=3.6pt] (12.,0.)-- (12.,-2.);
\draw [line width=2.8pt,dash pattern=on 7pt off 7pt] (10.,-2.)-- (12.,-2.);
\draw [line width=3.6pt] (12.,-2.)-- (12.,-4.);
\draw [line width=3.6pt] (10.,-4.)-- (12.,-4.);
\draw [line width=2.8pt,dash pattern=on 7pt off 7pt] (12.,-4.)-- (12.,-6.);
\draw [line width=2.8pt,dash pattern=on 7pt off 7pt] (10.,-6.)-- (10.,-4.);
\draw [line width=2.8pt,dash pattern=on 7pt off 7pt] (10.,-6.)-- (8.,-6.);
\draw [line width=2.8pt,dash pattern=on 7pt off 7pt] (8.,-4.)-- (8.,-6.);
\draw [line width=2.8pt,dash pattern=on 7pt off 7pt] (8.,-6.)-- (6.,-6.);
\draw [line width=2.8pt,dash pattern=on 7pt off 7pt] (6.,-6.)-- (4.,-6.);
\draw [line width=2.8pt,dash pattern=on 7pt off 7pt] (4.,-4.)-- (4.,-6.);
\draw [line width=2.8pt,dash pattern=on 7pt off 7pt] (4.,-6.)-- (2.,-6.);
\draw [line width=2.8pt,dash pattern=on 7pt off 7pt] (2.,-6.)-- (0.,-6.);
\draw [line width=3.6pt] (0.,-4.)-- (0.,-6.);
\draw [line width=3.6pt] (0.,-6.)-- (-2.,-6.);
\draw [line width=3.6pt] (-2.,-4.)-- (-2.,-6.);
\begin{scriptsize}
\draw [fill=black] (-2.,6.) circle (3.5pt);
\draw [fill=black] (-2.,4.) circle (3.5pt);
\draw [fill=black] (0.,6.) circle (3.5pt);
\draw [fill=black] (0.,4.) circle (3.5pt);
\draw [fill=black] (2.,4.) circle (3.5pt);
\draw [fill=black] (2.,6.) circle (3.5pt);
\draw [fill=black] (4.,6.) circle (3.5pt);
\draw [fill=black] (4.,4.) circle (3.5pt);
\draw [fill=black] (6.,6.) circle (3.5pt);
\draw [fill=black] (6.,4.) circle (3.5pt);
\draw [fill=black] (8.,6.) circle (3.5pt);
\draw [fill=black] (8.,4.) circle (3.5pt);
\draw [fill=black] (8.,2.) circle (3.5pt);
\draw [fill=black] (6.,2.) circle (3.5pt);
\draw [fill=black] (4.,2.) circle (3.5pt);
\draw [fill=black] (2.,2.) circle (3.5pt);
\draw [fill=black] (0.,2.) circle (3.5pt);
\draw [fill=black] (-2.,2.) circle (3.5pt);
\draw [fill=black] (-2.,0.) circle (3.5pt);
\draw [fill=black] (0.,0.) circle (3.5pt);
\draw [fill=black] (2.,0.) circle (3.5pt);
\draw [fill=black] (4.,0.) circle (3.5pt);
\draw [fill=black] (6.,0.) circle (3.5pt);
\draw [fill=black] (8.,0.) circle (3.5pt);
\draw [fill=black] (8.,-2.) circle (3.5pt);
\draw [fill=black] (6.,-2.) circle (3.5pt);
\draw [fill=black] (4.,-2.) circle (3.5pt);
\draw [fill=black] (2.,-2.) circle (3.5pt);
\draw [fill=black] (0.,-2.) circle (3.5pt);
\draw [fill=black] (-2.,-2.) circle (3.5pt);
\draw [fill=black] (10.,6.) circle (3.5pt);
\draw [fill=black] (10.,4.) circle (3.5pt);
\draw [fill=black] (10.,2.) circle (3.5pt);
\draw [fill=black] (10.,0.) circle (3.5pt);
\draw [fill=black] (10.,-2.) circle (3.5pt);
\draw [fill=black] (10.,-4.) circle (3.5pt);
\draw [fill=black] (8.,-4.) circle (3.5pt);
\draw [fill=black] (6.,-4.) circle (3.5pt);
\draw [fill=black] (4.,-4.) circle (3.5pt);
\draw [fill=black] (2.,-4.) circle (3.5pt);
\draw [fill=black] (0.,-4.) circle (3.5pt);
\draw [fill=black] (-2.,-4.) circle (3.5pt);
\draw [fill=black] (12.,6.) circle (3.5pt);
\draw [fill=black] (12.,4.) circle (3.5pt);
\draw [fill=black] (12.,2.) circle (3.5pt);
\draw [fill=black] (12.,0.) circle (3.5pt);
\draw [fill=black] (12.,-2.) circle (3.5pt);
\draw [fill=black] (12.,-4.) circle (3.5pt);
\draw [fill=black] (12.,-6.) circle (3.5pt);
\draw [fill=black] (10.,-6.) circle (3.5pt);
\draw [fill=black] (8.,-6.) circle (3.5pt);
\draw [fill=black] (6.,-6.) circle (3.5pt);
\draw [fill=black] (4.,-6.) circle (3.5pt);
\draw [fill=black] (2.,-6.) circle (3.5pt);
\draw [fill=black] (0.,-6.) circle (3.5pt);
\draw [fill=black] (-2.,-6.) circle (3.5pt);
\end{scriptsize}
\end{tikzpicture}
\caption{Illustration of assertion (1) in Theorem \ref{t:cluster} for the square lattice. Edges on which the local time hits $0$ are erased. Edges in solid lines are crossed an odd number of times, edges in dashed lines are crossed an even number of times (possibly zero). The graph in solid lines is chosen uniformly among all subgraphs of the clusters with even degrees.  }
\label{f:cluster}
\end{figure}
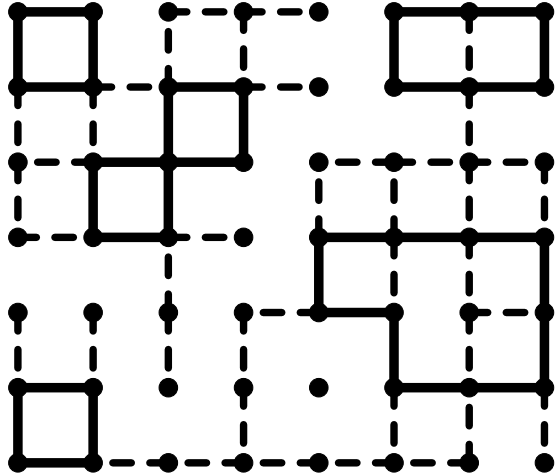

See Figure \ref{f:cluster} for an illustration of the theorem. In (2), we do not condition on the clusters. Conditioning on the clusters would yield for case (a) a conditioning of the loop soup ${\mathcal B}(\rho(e),\widehat \L(e_1),\widehat \L(e_2))$ on having positive local time field if the edge is in a cluster, or on having local time $0$ at some point of the edge if it is not contained in a cluster. We will prove the theorem via Theorem \ref{t:stargraph}. When taking for $W$ the set of all vertices of $\G$, the star graph $\G^*$ is obtained by extending each edge at its endpoints. Theorem  \ref{t:stargraph} says that the trace of the conditioned loop soup on the edges of $\G$ is distributed as  the trace of the conditioned loop soup on the old edges of $\G^*$. Under the measure $\n$ associated to $\G^*$ (call it $\n^*$), the extended edges of $\G^*$ can be crossed at most once by a loop. It has the following interpretation for the loop soup on $\G$: edges of $\G$ with an odd number of crossings are in correspondence with extended edges in $\G^*$ which are crossed once by the loop soup, while edges  of $\G$ with an even number of crossings are in correspondence with extended edges which are not crossed by the loop soup on $\G^*$. We will show that under $\n^*$, conditionally on the clusters, the configuration of edges crossed by loops is indeed uniform (Proposition \ref{p:uniform}), which will give assertion (1). Assertion (2) is almost immediate: under $\n^*$, the loops on the extended edges are loop soups on an interval of the real line.

\bigskip

Theorem \ref{t:cluster} can actually be  proved via the link with the random current model described by Werner \cite{wernermarkov} (see also Proposition 3.2 in Le Jan \cite{lejan17} and Proposition 6.7 in Lawler \cite{lawler18}). Werner shows  that conditionally on the local times at vertices of the graph, the number of unoriented jumps on the edges is distributed as a random current model. This model assigns an integer $n_e$ on each edge $e$ proportionally to 
$$
\prod_e {(\beta_e)^{n_e} \over n_e !}
$$

\noindent for all configurations $(n_e)_e$ such that at any vertex, the sum of the $n_e$'s over all incoming edges $e$ is even. When one chooses  the weight $\beta_e$ to be ${1\over \rho(e)}\sqrt{\widehat \L(e_1)\widehat \L(e_2)}$, one recovers the law of the number of crossings of the edges conditioned on the local times at vertices.  Notice that the random current model can be described as first choosing the parity of  $n_e$, before choosing $n_e$ independently for each edge. More precisely, one first picks  $\alpha_e \in \{0,1\}$  on the edges such that the sum of the $\alpha_e$'s over all incoming edges of a vertex is always even, proportionally to $\prod_e \cosh(\beta_e)^{1-\alpha_e}\sinh(\beta_e)^{\alpha_e}$.  Given $(\alpha_e)_e$,  one chooses $n_e$ independently on each edge $e$ with probability ${1\over \cosh(\beta_e)} {\beta^{n_e}\over n_e!}$ on the set of even integers if $\alpha_e=0$ and ${1\over \sinh(\beta_e)} {\beta^{n_e}\over n_e!}$ on the set of odd integers if $\alpha_e=1$.    Lupu and Werner in \cite{lw16} exhibit a link between the FK-Ising model and the current model. In the setting of the loop soup, a cluster of the FK-Ising model is a cluster of the loop soup  on the metric graph. By looking at their proof of this coupling,  one can check that, conditionally on  the clusters, the configuration $(\alpha_e)_e$ is chosen uniformly among all admissible configurations, which is basically assertion (1).  
One can then observe (see for example the decomposition of squared Bessel bridges given by Pitman and Yor in \cite{pitmanyor}) that the numbers of crossings $n_e$ match those of the loop soups ${\mathcal B}(\rho(e),\widehat \L(e_1),\widehat \L(e_2))$ and ${\mathcal C}(\rho(e),\widehat \L(e_1),\widehat \L(e_2))$, which implies assertion (2). Likewise, Theorems \ref{t:markov} and \ref{t:markov2} should have a proof using a similar direction. Consider the random current representation, look at the odd/even crossings configuration  inside the explored part of a cluster  and use a combinatorial argument as in Section \ref{s:n}.

\bigskip

\paragraph{Conditioning a loop-soup on its occupation field}

Theorem \ref{t:cluster} implies that the trace of the loop soup on the edges conditioned on the occupation field  is given by the conditioned versions of ${\mathcal B}(\rho(e),\widehat \L(e_1),\widehat \L(e_2))$ and ${\mathcal C}(\rho(e),\widehat \L(e_1),\widehat \L(e_2))$. The conditional law of a one-dimensional loop soup given its occupation field follows from \cite{eyzburglar}. From the link discovered by Lupu \cite{lupu18} with mu-processes, it amounts to describing the law of a mu-process conditioned on its local time. The paper \cite{eyzburglar} shows that this problem can be seen as a reformulation of the Perkins' disintegration theorem \cite{perkins}, \cite{ethmarch}. In particular, it recovers the link between the conditional version of the one-dimensional loop soup and the Bass-Burdzy flow introduced in \cite{bassburdzy}, which has been described by Warren and Yor  \cite{warren-yor} and Warren \cite{warren} (it is stated there for the Brownian motion, but can be interpreted in terms of a loop soup with intensity $1$) and Lupu, Sabot and Tarr\`es \cite{lst20} (in the case of the loop soup with intensity $\frac12$), extending it to all intensities of the loop soup. The similar problem for the Brownian motion  has also been studied in  \cite{aldouscondbr}, \cite{berestycki} and \cite{GKW}. We refer to Section \ref{s:loopsoupcond} for a presentation of the conditioned version of the one-dimensional loop soup and its link with Fleming--Viot processes but mostly refer to \cite{eyzburglar}.  In Section \ref{s:reconstruct}, we will complete this description by reconstructing the loops of the loop soup (and not only its trace on the edges), following the gluing operation of Werner \cite{wernermarkov}.

It is also possible to take another point of view, and dynamically reconstruct the loops. Let us present one way which  follows from a series of works of Sabot, Tarr\`es \cite{st16} and Lupu, Sabot, Tarr\`es \cite{lst19}, \cite{lst20}.  In \cite{st16}, considering the discrete loop soup, Sabot and Tarr\`es show that the discrete loops at a vertex $v_0$ can be traced via a process that can be roughly described as follows. At time $0$, the total local time at each vertex is known. The process starts at $v_0$, and eats its local time. Jump rates are given in terms of the remaining local time available at that vertex. When it jumps to a neighboring vertex, it starts eating  the local time at its new position then jumps according to the remaining local time. The process continues similarly until it ends at the starting vertex $v_0$, when all loops at $v_0$ will be reconstructed and the local time at $v_0$ will be exhausted. One can then choose another vertex and likewise build the remaining loops at that vertex. Step by step all loops will be discovered. In \cite{lst19}, Lupu, Sabot and Tarr\`es extended this work to invert Lupu's isomorphism between the signed Gaussian free field and the discrete loop soup. Incorporating the sign prevents for example jumps between vertices of opposite signs. They describe the process which traces back the loops  in terms of  Poisson point processes governing the jumps, whose clocks are related to the available local time at the current vertex. In \cite{lst20}, Lupu, Sabot and Tarr\`es study the closely related problem of inverting the Ray-Knight identity on the line.  It can be interpreted as considering a Brownian loop soup on the real line, and tracing back the loops at $0$ when conditioning on the occupation field. This reconstruction involves a self-interacting diffusion, related to the Bass--Burdzy flow, which traces loops up to a certain time when the diffusion will have exhausted the local time field at some position (which, unlike the discrete case, happens before having completed the discoveries of all loops touching $0$). Reconstructing the whole loop soup would involve a concatenation of ``excursions'' of this self-interacting diffusion. In the case of a general metric graph, one should further add a gluing procedure at the vertices.

\bigskip


{\bf Organization of the paper}. In Section \ref{s:model}, we recall some facts about loop soups. In Section \ref{s:loopsoupcond}, we discuss the conditioning of one-dimensional loop soups given their occupation fields. Section \ref{s:n} presents the construction of the measure $\n$. Given a graph and a set of vertices, called star vertices, we naturally associate a probability measure under which the loop soup has local time $0$ at the star vertices but has positive local time field on the adjacent edges. This section also contains the analogs of Theorems \ref{t:markov} and \ref{t:markov2} under $\n$, which are Theorems \ref{t:markovn} and \ref{t:markovn2}. In Section \ref{s:detachment}, we prove a key lemma to relate the loop soup on the graph $\G$ and on a graph where we detach edges adjacent to a vertex. By induction, it will lead to the proof of Theorem \ref{t:stargraph} in Section \ref{s:stargraph}. Theorem \ref{t:markov}, Theorem \ref{t:markov2} and  Theorem \ref{t:cluster} will be consequences and proved in Section \ref{s:markov} and Section \ref{s:clusterproof}. Section \ref{s:lejan} contains the relation with Le Jan's isomorphism theorem. Section \ref{s:reconstruct} explains how to reconstruct the loops of the loop soup.

\bigskip

{\bf Acknowledgements}. We thank Titus Lupu for explaining to us how the results of the paper  \cite{lst20} can be used to construct the loop soup conditioned on the occupation field, and  Yueyun Hu and Zhan Shi for stimulating discussions on the link with \cite{eyzburglar}.


\section{Generalities on loop soups on metric graphs}
\label{s:model}

We briefly present the model but refer  to Lupu \cite{lupu16} and the references therein for the precise definitions. We refer to Lawler \cite{lawler18}, Le Jan \cite{lejan} and Powell and Werner \cite{pw20} for generalities on discrete loop soups. 

\bigskip

We consider a connected discrete graph $(V,E)$ where the sets of vertices $V$ and edges $E$ are supposed finite. We allow self-loops and multiple edges.  An edge $e\in E$ is seen as a one-dimensional segment assigned with some length $\rho(e) \in (0,\infty)$. It defines a metric graph which we  denote by $\G$.  Each edge has two orientations, an edge together with an orientation being called a directed edge. When a directed edge is oriented from $v$ to $w$, one calls $v$ the root and $w$ the terminal point of the edge. We call $d_v$  the number of directed edges rooted at the vertex $v$ (a self-loop counts for two).

\bigskip

If we fix a set of vertices $V_0$, one can define a loop soup on the metric graph associated to the loop measure of the Brownian motion on $\G$ killed when hitting $V_0$. We will only consider the loop soup at intensity $1/2$. We call it for short loop soup with sink $V_0$. We  consider loops as unrooted and unoriented, unless stated otherwise. We attach to the graph a distinguished vertex ${\mathfrak v}$ and we denote by $\L$ the loop soup with sink $\{\mathfrak v\}$. By convention, we do not add $\mathfrak v$ to the set $V$, nor the edges adjacent to $\mathfrak v$ to the set $E$. We can describe the law of $\L$ as follows.

\bigskip

Consider a loop of $\L$ which crosses at least one edge. One can keep track of the vertices it visits, arbitrarily choosing a first vertex and an orientation. Forgetting about the starting point and the orientation,   it  defines a (unrooted, unoriented) discrete loop on the set $V$. We call the collection of all discrete loops thus obtained  the discrete loop soup and we denote it by $\overline \L$. A realization of $\overline \L$ is called a discrete loop configuration.  

\bigskip

A rooted oriented discrete loop can be seen  as a sequence $\ell=(e_0,\ldots,e_{n-1})$ of adjacent directed edges (the edges crossed by the loop) such that the terminal point of $e_i$ is the root of $e_{i+1}$, and the root of $e_0$ is the terminal point of $e_{n-1}$. Call $v_i$ the root of $e_i$, and for any vertex $v$, set 
$$
a_v := \sum_{e: v \rightarrow } {1\over \rho(e)}
$$

\noindent  where the sum runs over all directed edges $e$ rooted at $v$. While the rooted loop measure of $\ell$ is given by ${1\over n}\prod_{i=0}^{n-1} { 1/\rho(e_i) \over a_{v_i}}$, we express, see \cite{wernermarkov}, the loop measure of the corresponding unrooted  loop (still denoted by $\ell$)  as 
$$
\mu(\ell):= {1\over J(\ell)} \prod_{i=0}^{n-1} { 1/\rho(e_i) \over a_{v_i}}
$$

\noindent where $J(\ell)$ is the maximal number of times $\ell$ can be written as the concatenation of the same loop. Following \cite{wernermarkov}, we let $\nu$ be the projection of ${1\over 2}\mu$ (since the loop soup has intensity $\frac12$) on the space of unoriented loops. The discrete loop soup  $\overline \L$ is a Poisson point process of unrooted, unoriented loops with intensity $\nu$. Notice that if, for example, the loop $\ell$ crosses some edge only once, then $J(\ell)=1$ and the two orientations of $\ell$ give rise to two different loops, hence   (still denoting by $\ell$ the corresponding unrooted unoriented loop under $\nu$),
\begin{equation}\label{eq:nu}
\nu(\ell)=\mu(\ell) = \prod_{i=0}^{n-1} { 1/\rho(e_i) \over a_{v_i}}.
\end{equation}

\bigskip

Let $\overline L$ be a discrete loop configuration. For an edge $e\in E$, we define $n_{\overline L}(e)$ as the number of crossings (in any direction) of the edge $e$.  For a vertex $v\in W$, we let 
$$
k_{v}=k_v(\overline{ L}):= {1\over 2}\sum_{e:v \rightarrow} n_{\overline L}(e)
$$

\noindent where we recall that the sum runs over directed edges rooted at $v$. In words, $k_v$ is the number of departures from $v$, or equivalently arrivals to $v$. Conditioned on the event $\{\overline \L=\overline L\}$, the local times $\widehat \L(v)$ at vertices $v\in V$ are independent with $\widehat \L(v)$ being gamma($k_v+{1\over 2}$, ${a_v\over 2}$) distributed (the density at $x\in \mathbb{R}$ of the gamma($p,\alpha$) distribution is $ { x^{p-1}\alpha^p \over \Gamma(p)}{\mathrm e}^{-\alpha x}{\bf 1}_{{\mathbb R}^+}(x)$). Each visit of $v$ by a loop is accompanied by a gamma($1, {a_v\over 2}$) (i.e. exponential) random variable. The accumulated local time of the loops visiting $v$ which do not cross any adjacent edge gives a gamma(${1\over 2}, {a_v\over 2}$) random variable (partition it with respect to a Poisson--Dirichlet($0,{1\over 2}$) distribution to recover each contribution, see \cite{lejan}).

\bigskip

Conditionally on $\overline \L$ and on $(\widehat \L(v),\, v\in V)$,  one can reconstruct the trajectories inside the edges. Each crossing is a Brownian crossing from one extremity to the other. On an edge $e$ with endpoints $v$ and $w$, noting $x_v$ and $x_w$ their respective local times, one  adds away from each extremity independent Brownian excursions which do not cross the edge, up to local time $x_v$ and $x_w$ respectively, and loops inside $e$ which do not hit $v$ nor $w$, distributed as a Brownian loop soup inside $e$ with sink $\{v,w\}$. 

\bigskip

The following lemma follows from arguments of Powell and Werner \cite{pw20} and Werner \cite{wernermarkov}.
\begin{lemma}
For $v\in V$ and $k\ge 0$,  
$$
P\big(k_v(\overline{\mathcal L}) = k\big)
=
\frac{p(v,v)^k}{k!}\frac{\Gamma(k+\frac12)}{\Gamma(\frac12)} \sqrt{1-p(v,v)}
$$ 

\noindent  where $p(v,v)$ is the probability that the random walk starting at $v$ returns to $v$ without hitting $\mathfrak v$ (when at position $w$, this random walk crosses an adjacent edge $e$ with probability $\frac{1/\rho(e)}{a_w}$).
\end{lemma}

\noindent {\it Proof}. For simplicity, we suppose that $v$ is not adjacent to a self-loop. We use the trick presented in \cite{pw20} and \cite{wernermarkov} and replace each  edge adjacent to $v$ by $K$ copies, where $K$ is a big number. The length of a duplicated edge is multiplied by $K$ so that  the  random walk keeps the same distribution. The discrete loop soup  keeps the same distribution when projected onto the original graph. For each  edge adjacent to $v$, index the duplicated edges by $1,2,\ldots,K$. With  probability tending to $1$ when $K\to \infty$, each index is used by the loup soup at most once, which allows us to use equation \eqref{eq:nu} for the loop measure. Call $F_K$ this event. We compute the probability that $k_v(\overline{\mathcal L}) = k$ for this loop soup, on the event $F_K$. If $\ell_1,\ldots,\ell_p$ are $p$ (unoriented, unrooted) loops visiting $v$ which use distinct indices each time they cross an edge adjacent to $v$, the probability that these are the loops of the loop soup visiting $v$ is proportional to the product of the transition probabilities (the constant being the probability that no discrete loop visits $v$).   We want to sum over all collections of such loops which verify that the number of visits of $v$ is $2k$.

One can decompose the loops visiting $v$ into $k$ excursions. Each excursion corresponds to a pair of indices\footnote{this is why we earlier discarded self-loops}. We index the excursion by the minimum of the two indices, orientate it so that it leaves through the edge with  minimal index, and rank the excursions in increasing order of their indices. Projecting onto the original graph, we get $k$ (not necessarily distinct) oriented excursions. We fix $\widetilde \ell_1,\ldots,\widetilde \ell_k$  oriented excursions away from $v$ on the original graph, and sum  over collections of loops which give rise to these excursions by the previous procedure. Observe that they all carry the same weight which is also proportional to $K^{-2k}$ times the transition probabilities in the original graph (the term $K^{-2k}$ arises because we need to choose one of the $K$ duplicated edges each time we leave or return to $v$).    We count the number of  collections of loops associated to the excursions $\widetilde \ell_1,\ldots,\widetilde \ell_k$. To construct a suitable configuration,  we choose $k$ pairs of indices, $i_1<j_1,\ldots,i_k<j_k$, all indices being distinct, ordered so that $i_1<i_2<\ldots <i_k$. For each pair of indices $i_s<j_s$, $1\le s\le k$, draw the oriented excursion $\widetilde \ell_s$, leaving $v$ through index $i_s$ and coming back through index $j_s$. Erase the orientation of the excursion. Take a uniform pairing of $\{i_1,j_1,i_2,j_2,\ldots,i_k,j_k\}$. We get a collection of unoriented loops visiting $v$ associated to the excursions $\widetilde \ell_1,\ldots,\widetilde \ell_k$. The number of such collections is therefore 
$$
 \frac{ \frac{K(K-1)}{2}\frac{(K-2)(K-3)}{2}\ldots \frac{(K-2k+2)(K-2k+1)}{2}}{k!}.
$$ 

\noindent for the choice of the indices multiplied by
$$
(2k-1)\times (2k-3)\times \ldots \times 3\times 1 = 2^{k} {\Gamma(k+{1\over 2}) \over \Gamma({1\over 2})}
$$
 
 \noindent for the choice of the pairing. It remains to sum over all oriented excursions $\widetilde \ell_1,\ldots, \widetilde \ell_k$ and make $K\to \infty$. $\Box$

\bigskip

We have $p(v,v)=1- {H(v,\mathfrak v) \over a_v/2}$ where $H(v,\mathfrak v)$ is the mass under the Brownian excursion measure at $v$  of excursions which hit $\mathfrak v$. A computation shows that the accumulated local time at $v$ is gamma($\frac12,H(v,\mathfrak v)$) distributed. This fact is also a consequence of the Feynman--Kac formula, see \cite{lejan}. We can also check that conditionally on the local time at $v$ being $x$, excursions of the (non-discrete) loop soup away from $v$ are distributed as a Poisson point process of Brownian excursions up to local time $x$, conditioned on not hitting $\mathfrak v$.

\bigskip

We conclude this section by defining what we mean by trace of the loop soup on a subgraph. Let $\widetilde{\G}$ be a (not necessarily connected)  subgraph of $\G$. We naturally see $\widetilde{\G}$ as a metric graph. Continuous loops which intersect $\widetilde{\G}$  can either stay inside it, or hit also its complement. When a loop hits $\widetilde{\G}$ and its complement, it will make a countable number of excursions in $\widetilde{\G}$, away from the vertices which lie at the boundary of $\widetilde{\G}$ (or none if $\widetilde{\G}$ consists of a single vertex for example). The trace of the loop soup on $\widetilde{\G}$ is by definition the collection of loops which stay inside $\widetilde{\G}$, and all the excursions in $\widetilde{\G}$ made by loops. We only keep the collection of these trajectories and not the information on whether two excursions belong to the same loop or not. We recall from the introduction that the trace of the loop soup outside $\widetilde{\G}$ is made of loops which do not hit $\widetilde{\G}$ and excursions away from vertices on the boundary of $\widetilde{\G}$ which do not hit $\widetilde{\G}$ in between.  If $E'$ is a set of edges, by trace of the loop soup on/outside $E'$, we mean the trace of the loop soup on/outside $\widetilde{\G}$ where $\widetilde{\G}$ is the subgraph associated to the edges of $E'$ and their endpoints. The trajectories of the trace are considered unoriented.

\section{Description of one-dimensional loop soups conditioned on their occupation fields}
\label{s:loopsoupcond}

As shown in Lupu \cite{lupu18}, one way to construct the Brownian loop soup on $\mathbb R$ is to root each loop at its minimum, and see the loops as the excursions of a process above its current infimum. This process is actually a mu-process, also called perturbed reflected Brownian motion, see \cite{yor92} and the references therein. It is defined  for $t\ge 0$ as $X_t:=|B_t| - \mu {\mathfrak L}_t$ and for $t\le 0$ as $X_t:=|\widetilde B_{-t}| + \mu \widetilde {\mathfrak L}_{-t}$, where  $B$, $\widetilde B$ are independent standard Brownian motions, ${\mathfrak L}$, $\widetilde {\mathfrak L}$ their respective local time processes at position $0$ and $\mu$ is the inverse of the intensity of the loop soup, i.e. $\mu=2$ in our setting. 

\bigskip

One-dimensional loop soups enjoy (at any intensity) a spatial Markov property. If $[a_1,b_1] \subset [a_0,b_0]$, the law of the trace of the loop soup on $[a_1,b_1]$ given the trace on $[a_0,b_0] \backslash [a_1,b_1]$ only depends on the local times at $a_1$ and $b_1$ (see \cite{perman-werner}, \cite{werner95}, \cite{yor92}, \cite{eyzburglar} for different proofs of this independence property in the setting of mu-processes). 

\bigskip

\begin{figure}[h!]
\begin{tikzpicture}[line cap=round,line join=round,>=triangle 45,x=0.3cm,y=0.3cm]

\clip(-27.32924134140175,-7.712935517349647) rectangle (20.395670605682085,11.593746900579648);

\draw[color=black] (-16.32924134140175,0.) -- (15,0.);
\foreach \x in {-16.,-14.,-12.,-10.,-8.,-6.,-4.,-2.,2.,4.,6.,8.,10.,12.,14.,16.,18.}
\draw[shift={(\x,0)},color=black] (0pt,-2pt);
\draw[color=black] (0.,-6) -- (0.,12.393746900579648);

\draw[dashed] (-16.,4.9)--(15,4.9);
\node at (-17.,4.5) {$\rho$};
\node at (-17.,-0.5) {$0$};
\draw [dotted] (0.,0.)-- (-0.06659852520710975,0.16951021971730434)-- (-0.44043911647646866,0.15959451741861955)-- (-0.6166147630670642,0.4521844660194193)-- (-1.5268556037851408,0.4521844660194193)-- (-1.5897645972405297,0.8973626448783579)-- (-1.8498442892012326,0.9043689320388371)-- (-1.9502772770117538,1.3425774307763918)-- (-2.6279533949763625,1.356553398058255)-- (-2.804129041566958,1.8087378640776726)-- (-4.1988529104091725,1.835336950314109)-- (-4.404391164764867,2.234323243860654)-- (-5.211862878305096,2.2609223300970904)-- (-5.373357221013142,2.6333095374071993)-- (-8.28678105498576,2.605788121073672)-- (-8.38678105498576,3.0057881210736723)-- (-8.746781054985764,3.0257881210736723)-- (-8.846781054985769,3.605788121073672)-- (-10.206781054985772,3.605788121073672)-- (-10.366781054985772,4.005788121073673);
\draw [dotted] (-10.35989493282946,3.9885728156828923)-- (-10.687505823059531,3.955651056539124)-- (-10.853857429864824,4.235651056539124)-- (-11.41340374366445,4.215651056539125)-- (-11.519263857085999,4.575651056539124)-- (-12.41151338449621,4.555651056539125)-- (-12.472004877879954,4.915651056539124)-- (-12.744216598106796,4.935651056539125)-- (-12.850076711528345,5.275651056539125)-- (-13.969169339127593,5.275651056539125)-- (-14.09015232589508,5.635651056539125)-- (-15.526825293758975,5.655651056539124)-- (-15.693176900564268,5.955651056539125);
\draw  (-15.522982076505059,5.655597554909458)-- (-15.482625601951536,6.463577955454034)-- (-15.289965051492821,6.06459166190749)-- (-15.12482743681392,8.032924043403781)-- (-15.028497161584564,7.128555111364943)-- (-14.89088248268548,8.990491147915492)-- (-14.780790739566214,6.4901770416904725)-- (-14.629414592777223,7.421145059965741)-- (-14.560607253327682,6.117789834380363)-- (-14.422992574428566,6.94236150770989)-- (-14.093177367806444,5.635693168306785)-- (-14.051432941401075,6.463577955454034)-- (-13.940415068476739,5.275651056539125)-- (-13.734919179933183,7.15515419760138)-- (-13.624827436813916,5.958195316961745)-- (-13.445928354245108,8.032924043403781)-- (-13.225744868006576,5.559209023415199)-- (-13.060607253327674,6.330582524271853)-- (-12.802303307852291,5.122213121089567)-- (-12.757854959749693,6.011393489434616)-- (-12.404243321815532,4.555814016842433)-- (-12.303726519382717,5.585808109651636)-- (-12.12482743681391,4.947430039977163)-- (-11.945928354245101,6.037992575671054)-- (-11.849598079015745,4.947430039977163)-- (-11.711983400116662,5.665605368360944)-- (-11.519263857085999,4.575651056539124)-- (-11.340423767089138,6.383780696744727)-- (-11.230332023969869,5.293218161050835)-- (-11.092717345070788,6.011393489434616)-- (-10.853857429864824,4.235651056539124)-- (-10.624827436813906,4.920830953740727)-- (-10.268096308293375,3.7590762543426806)-- (-10.211983400116655,4.7612364363221085)-- (-10.033084317547832,4.043061107938327)-- (-9.96427697809829,5.585808109651636)-- (-9.799139363419407,5.0538263849229095)-- (-9.67528615241023,7.15515419760138)-- (-9.551432941401057,6.357181610508289)-- (-9.455102666171697,7.873329525985162)-- (-9.248680647823075,5.506010850942327)-- (-9.198221932226733,6.4255853797829845)-- (-9.069781565254264,4.468646487721309)-- (-8.94592835424509,5.293218161050835)-- (-8.849385310579928,3.5908766419189075)-- (-8.746781054985764,3.0257881210736723)-- (-8.52375857837244,3.697272986864653)-- (-8.28678105498576,2.605788121073672)-- (-8.28678105498576,2.605788121073672)-- (-8.230332023969858,3.7504711593375273)-- (-8.106478812960685,3.2982866933181096)-- (-7.8862953267221485,5.718803540833817)-- (-7.776203583602883,5.2666190748143995)-- (-7.638751228702499,7.075356938892071)-- (-7.462575582111903,6.197587093089671)-- (-7.253267803786367,11.517404340376945)-- (-7.005561381768016,8.219117647058845)-- (-6.881708170758843,9.176684751570548)-- (-6.761543321720154,6.995559680182758)-- (-6.49638706984141,8.751099371787568)-- (-6.179873308373519,5.000628212450036)-- (-6.056020097364344,5.904997144488871)-- (-5.904643950575354,4.069660194174763)-- (-5.780790739566179,5.346416333523709)-- (-5.377830430789857,2.659908623643638);
\draw  (0.,0.)-- (0.059580912712818145,-0.17689314656749086)-- (0.1205414173041166,-0.34659908623643676)-- (0.20590595728840583,0.4163590380494338)-- (0.28216184800737354,0.)-- (0.37536349221944526,1.6119263641665778)-- (0.45324463091470335,-0.34659908623643676)-- (0.5194023969108288,0.2905098458265765)-- (0.5288589976443822,-0.6865990862364366)-- (0.6210769178694524,0.)-- (0.6888599318418682,-0.4016607113991384)-- (0.765115822560836,1.0141427011080058)-- (0.8161935912171617,-0.7065990862364366)-- (0.8766850846009047,-1.0665990862364367)-- (0.9854106179711872,0.35343444193800516)-- (1.044720755197051,-0.37019841334342407)-- (1.1125037691694668,3.027729776673722)-- (1.2395969203677464,-0.49604760556628136)-- (1.3327985645798182,0.5107459322165767)-- (1.4816000184383353,-1.0665990862364367)-- (1.5530933599901693,-0.55897220167771)-- (1.5874601318598855,-1.3465990862364368)-- (1.7140780181746569,0.6680574224951483)-- (1.8411711693729365,-0.8735951822348531)-- (1.9258999368384562,0.32197214388229084)-- (2.0529930880367355,-0.8106705861234245)-- (2.1631404857419114,2.0838608350022927)-- (2.273287883447087,-0.8735951822348531)-- (2.3834352811522623,0.)-- (2.4192181658863525,-1.3465990862364368)-- (2.468164048617782,-0.6533590958448529)-- (2.5250782793079027,-1.6865990862364368)-- (2.6460944602953735,-0.7162836919562815)-- (2.7392961045074453,-1.251142758903425)-- (2.824024871972965,-0.33873611528770975)-- (2.948518732994104,-1.6865990862364368);
\draw  (3.0146813038825497,-2.0265990862364367)-- (3.171412818581596,-0.43312300945485094)-- (3.315451723272979,-1.345529653070566)-- (3.425599120978155,-0.30727381723199365)-- (3.603529532655746,-1.7545395277948521)-- (3.722149807107474,-0.6533590958448511)-- (3.8492429583057537,-1.6286903355719948)-- (3.9254988490247213,-0.7162836919562797)-- (4.01468130388255,-2.0265990862364367)-- (4.128102853977094,-2.306599086236437)-- (4.255941042140249,0.5736705283280071)-- (4.425398577071288,-0.6848213939005654)-- (4.603328988748879,5.037920083348011)-- (4.7050035097075025,2.650182200005152)-- (4.82362378415923,5.935848273365157)-- (5.110027072583374,2.7445690941722947)-- (5.228647347035101,5.0286221698902296)-- (5.431996388952349,-0.7162836919562797)-- (5.618399677376492,1.7063132583337224)-- (5.855640226279947,-1.723077229739138)-- (5.965787623985123,0.7624443166622931)-- (6.224322135640615,-2.346599086236437)-- (6.313175570593754,-1.1567558647362801)-- (6.406377214805825,-2.352323190853424)-- (6.482633105524794,-1.4084542491819947)-- (6.647010623248842,-2.7865990862364365)-- (6.802027037497622,-1.180306355979135)-- (6.923222696345495,-2.352323190853424)-- (7.0093218283660764,-1.536478184697562)-- (7.117475632297443,-2.397226770767094)-- (7.189578168251688,-0.67572959862803)-- (7.282923269856118,-2.7465990862364373)-- (7.423122424392062,-1.6916149316834235)-- (7.570590600574772,-3.1703429403019783);
\draw  (7.570590600574772,-3.1703429403019783)-- (7.8,-1.9)-- (8.0,-2.8)-- (8.153949586639712,-2.3675457850405586)-- (8.350070357685157,-3.610815113081979)-- (8.45,-1.)-- (8.55,-2.)-- (8.65,1)--(8.75,-0.8)-- (8.85,2)--(8.95,1.5)--(9.05,4.)-- (9.15,3.5)--(9.25,3.6)--(9.35,2.3)-- (9.45,5)--(9.55,4.5)--(9.685340573622164,7.605265419075407)-- 
(9.8,6.5)--(9.9,6.7)--(10.099141731982235,3.5389703745400323)-- (10.362469741847735,5.141743603772954)-- 
(10.45,4)--(10.6,2.7)-- (10.7,3.2)--
(10.8,2)--
(10.96436233582602,2.7079027741970356)-- (11.152453771444235,-0.7944535415341638)--
(11.3,1)-- (11.4,0.2)-- (11.5,-1)--(11.6,-0.5)--(11.7,-2.5)--(11.75,-2)--(11.88601037035527,-3.6141471855550447)-- (12.1,-4.)-- (12.2,-2)--(12.3,-2.5)--(12.4,-1.8036070562363737)-- (12.5,-2)--
(12.6,-1.5)--(12.788849261322698,-3.317337328289689)-- 
(12.9,-2)--(13.07098641475002,-3.008183813587479)-- 
(13.2,-2.5)--(13.385088248209,-3.9883626897505504)-- (13.701093555145642,-4.554684054753408);
\draw  (-5.211862878305096,2.2609223300970904)-- (-5.121236338564625,3.1535789688965736)-- (-4.9709204288879105,2.7131067961165733)-- (-4.820604519211196,4.349146295013718)-- (-4.655257018566811,2.6187199019494303)-- (-4.459846335987082,3.468201949453717)-- (-4.404391164764867,2.234323243860654)-- (-4.1988529104091725,1.835336950314109)-- (-4.099088152762968,2.555795305838002)-- (-4.,2.)-- (-3.7834247424418677,4.715467660016575)-- (-3.693235196635839,3.5625888436208597)-- (-3.542919286959125,4.443533189180861)-- (-3.4226665592177534,2.555795305838002)-- (-3.302413831476382,3.4996642475094313)-- (-3.1370663308319964,2.1153231330580016)-- (-2.9717188301876107,3.0277297766737163)-- (-2.804129041566958,1.8087378640776726)-- (-2.6279533949763625,1.356553398058255)-- (-2.445613146319111,4.752543063905147)-- (-2.2502024637393823,1.674850960278001)-- (-2.114918145030339,3.468201949453717)-- (-1.9502772770117538,1.3425774307763918)-- (-1.8498442892012326,0.9043689320388371)-- (-1.7069312181661775,1.7284302452213545)-- (-1.5268556037851408,0.4521844660194193)-- (-1.363338596646768,0.9826804030522861)-- (-1.2581174598730678,2.964805180562288)-- (-1.1078015501963536,0.8253689127737145)-- (-0.9424540495519679,2.3040969213922873)-- (-0.8222013218105965,0.668057422495143)-- (-0.7320117760045679,1.8321624505565726)-- (-0.6166147630670642,0.4521844660194193)-- (-0.44043911647646866,0.15959451741861955)-- (-0.27389331607556355,0.8083196876987495)-- (0.,0.);
\draw [dotted] (0.,0.)-- (0.1205414173041166,-0.34659908623643676)-- (0.45324463091470335,-0.34659908623643676)-- (0.5288589976443822,-0.6865990862364366)-- (0.8161935912171617,-0.7065990862364366)-- (0.8766850846009047,-1.0665990862364367)-- (1.4816000184383353,-1.0665990862364367)-- (1.5874601318598855,-1.3465990862364368)-- (2.4192181658863525,-1.3465990862364368)-- (2.5250782793079027,-1.6865990862364368)-- (2.948518732994104,-1.6865990862364368)-- (3.0146813038825497,-2.0265990862364367)-- (4.01468130388255,-2.0265990862364367)-- (4.128102853977094,-2.306599086236437)-- (6.124322135640615,-2.346599086236437)-- (6.547010623248842,-2.7865990862364365)-- (7.182923269856118,-2.7465990862364373)-- (7.470590600574772,-3.1703429403019783)-- (8.189141841553921,-3.1689323996570193)-- (8.250070357685157,-3.610815113081979)-- (11.78601037035527,-3.6141471855550447)-- (12.,-4.)-- (13.285088248209,-3.9883626897505504)-- (13.601093555145642,-4.554684054753408);
\end{tikzpicture}

\begin{tikzpicture}[line cap=round,line join=round,>=triangle 45,x=0.5cm,y=0.45cm]

\clip(-15.5,-1.) rectangle (15.395670605682085,6.7);

\draw[color=black] (-11.,0.) -- (12.,0.);
\foreach \x in {-16.,-14.,-12.,-10.,-8.,-6.,-4.,-2.,2.,4.,6.,8.,10.,12.,14.,16.,18.}
\draw[shift={(\x,0)},color=black] (0pt,-2pt);

\draw (-11.,5.) -- (12.,5.);
\node at (-11.7,5.) {$\rho$};
\node at (-11.7,0.) {$0$};
\draw (11.7,0.2) -- (11.7,-0.2);
\node at (11.8,-0.6) {$\tau_{\ell_1}$};

\draw [dotted] (0.,0.)-- (-0.06659852520710975,0.16951021971730434)-- (-0.44043911647646866,0.15959451741861955)-- (-0.6166147630670642,0.4521844660194193)-- (-1.5268556037851408,0.4521844660194193)-- (-1.5897645972405297,0.8973626448783579)-- (-1.8498442892012326,0.9043689320388371)-- (-1.9502772770117538,1.3425774307763918)-- (-2.6279533949763625,1.356553398058255)-- (-2.804129041566958,1.8087378640776726)-- (-4.1988529104091725,1.835336950314109)-- (-4.404391164764867,2.234323243860654)-- (-5.211862878305096,2.2609223300970904)-- (-5.373357221013142,2.6333095374071993)-- (-8.28678105498576,2.605788121073672)-- (-8.38678105498576,3.025)-- (-8.746781054985764,3.0257881210736723)-- (-8.846781054985769,3.705788121073672)-- (-10.206781054985772,3.705788121073672)-- (-10.366781054985772,4.005788121073673);
\draw [dotted] (-10.35989493282946,3.9885728156828923)-- (-10.357505823059531,3.955651056539124)-- 
(-10.45,4.235651056539124)--
(-10.853857429864824,4.235651056539124);

\draw (-10.952717345070788,5.)-- (-10.853857429864824,4.235651056539124)-- (-10.624827436813906,4.920830953740727)-- (-10.268096308293375,3.7590762543426806)-- (-10.211983400116655,4.7612364363221085)-- (-10.033084317547832,4.043061107938327)-- (-9.96427697809829,5.);
\draw  (-9.1,5.)-- (-9.069781565254264,4.468646487721309)-- (-8.94592835424509,5.)-- (-8.849385310579928,3.5908766419189075)-- (-8.746781054985764,3.0257881210736723)-- (-8.52375857837244,3.697272986864653)-- (-8.28678105498576,2.605788121073672);
\draw [line width=1.5pt] (-8.28678105498576,2.605788121073672)-- (-8.230332023969858,3.7504711593375273)-- (-8.106478812960685,3.2982866933181096)--
(-7.9862953267221485,5.);

\draw [line width=1.5pt] (-7.7,5.) -- (-7.65,4.5)--(-7.6,5.);
 
\draw [line width=1.5pt] (-6.16020097364344,5.)-- (-6.054643950575354,4.069660194174763)-- (-6.00790739566179,5.); 

\draw [line width=1.5pt](-5.75,5.)--(-5.377830430789857,2.659908623643638);
\draw (-5.377830430789857,2.659908623643638)--(-5.21,2.26);
\draw (0.05,0.)--(0.10590595728840583,0.4163590380494338)-- (0.28216184800737354,0.)-- (0.37536349221944526,1.6119263641665778)-- 
(0.4,0.);
\draw (0.7210769178694524,0.)-- ((0.765115822560836,1.0141427011080058)-- 
(0.8,0);
\draw (1.0054106179711872,0.)-- (1.1125037691694668,3.027729776673722)-- (1.1595969203677464,0.);

\draw (1.25,0.)-- (1.3327985645798182,0.5107459322165767)-- 
(1.4,0.) ;

\draw (1.65,0.)-- (1.7140780181746569,0.6680574224951483)-- 
(1.8,0.);

\draw  (1.9258999368384562,0.)--  (2.1631404857419114,2.0838608350022927)-- (2.203287883447087,0.);

\draw  (4.1,0.)-- (4.255941042140249,0.5736705283280071)-- 
(4.35,0.);

\draw (4.525398577071288,0.)-- (4.58,5.)
(4.7,5.)-- 
(4.8050035097075025,2.650182200005152)-- 
(4.85,5.)
(5.1,5.)--
(5.110027072583374,2.7445690941722947)-- (5.228647347035101,4.286221698902296)-- 
(5.381996388952349,0.);

\draw (5.551996388952349,0.)-- (5.618399677376492,1.7063132583337224)-- 
(5.728399677376492,0.) ;

\draw (5.9,0.)-- (5.965787623985123,0.7624443166622931)-- 
(6.1, 0.);

\draw   (8.4,0.)-- (8.55,1)--(8.6,0);

\draw [line width=1.5pt](8.8,0)-- (8.85,2)--(8.95,1.5)--(9.05,4.)-- (9.15,3.5)--(9.25,3.6)--(9.35,2.3)-- (9.45,5);
\draw (9.7,5)-- (9.75,4.5) -- (9.8,5.)
(10.15,5.)--(10.299141731982235,3.5389703745400323)-- (10.462469741847735,5.)-- 
(10.65,4)--(10.8,2.7)-- (10.9,3.2)--
(11.0,2)--
(11.16436233582602,2.7079027741970356)--  (11.3,0.)(11.4,0.)--
(11.5,1)-- (11.7,0); 
\draw  (-5.211862878305096,2.2609223300970904)-- (-5.121236338564625,3.1535789688965736)-- (-4.9709204288879105,2.7131067961165733)-- (-4.820604519211196,4.349146295013718)-- (-4.655257018566811,2.6187199019494303)-- (-4.459846335987082,3.468201949453717)-- (-4.404391164764867,2.234323243860654)-- (-4.1988529104091725,1.835336950314109)-- (-4.099088152762968,2.555795305838002)-- (-4.,2.)-- (-3.7834247424418677,4.915467660016575)-- (-3.693235196635839,3.5625888436208597)-- (-3.542919286959125,4.443533189180861)-- (-3.4226665592177534,2.555795305838002)-- (-3.302413831476382,3.4996642475094313)-- (-3.1370663308319964,2.1153231330580016)-- (-2.9717188301876107,3.0277297766737163)-- (-2.804129041566958,1.8087378640776726)-- (-2.6279533949763625,1.356553398058255)-- (-2.445613146319111,4.852543063905147)-- (-2.2502024637393823,1.674850960278001)-- (-2.114918145030339,3.468201949453717)-- (-1.9502772770117538,1.3425774307763918)-- (-1.8498442892012326,0.9043689320388371)-- (-1.7069312181661775,1.7284302452213545)-- (-1.5268556037851408,0.4521844660194193);

\draw [line width=1.5pt]  (-1.5268556037851408,0.4521844660194193) --(-1.363338596646768,0.9826804030522861)-- (-1.2581174598730678,2.964805180562288)-- (-1.1078015501963536,0.8253689127737145)-- (-0.9424540495519679,2.3040969213922873)-- (-0.8222013218105965,0.668057422495143)-- (-0.7320117760045679,1.8321624505565726)-- (-0.6166147630670642,0.4521844660194193); 
\draw (-0.6166147630670642,0.4521844660194193)--(-0.44043911647646866,0.15959451741861955)-- (-0.27389331607556355,0.8083196876987495)-- (0.,0.);
\end{tikzpicture}

\caption{Top: the mu-process $X$. Its infimum process is represented in dotted lines. Excursions above it are the loops of the loop soup. Bottom: the process looked on the space interval $[0,\rho]$ (when one does not close the time gaps), stopped at some local time $\ell_1$ at position $0$. The associated loop soup consists of crossings, excursions away from $0$ and $\rho$, and loops inside $(0,\rho)$. In the negative time-axis, an excursion above infimum which hits $\rho$ will be associated to excursions away from $\rho$ (the concatenation of the paths back and forth from the infimum to $\rho$ yields a Brownian excursion from $\rho$). Such an excursion is represented in bold on the left of the picture. An excursion above infimum which does not hit $\rho$ is associated to a loop inside $(0,\rho)$ (in bold in the middle). In the positive time-axis, excursions above infimum which hit $0$ will yield excursions away from $0$. Excursions above infimum which hit $\rho$ will further yield excursions away from $\rho$ and crossings of the interval (a crossing is traced in bold on the right of the picture). In this example, there are 4 crossings. If one conditions the local time at $\rho$ on being $\ell_2$, one gets ${\mathcal B}(\rho,\ell_1,\ell_2)$. }
\label{f:muprocess}
\end{figure}
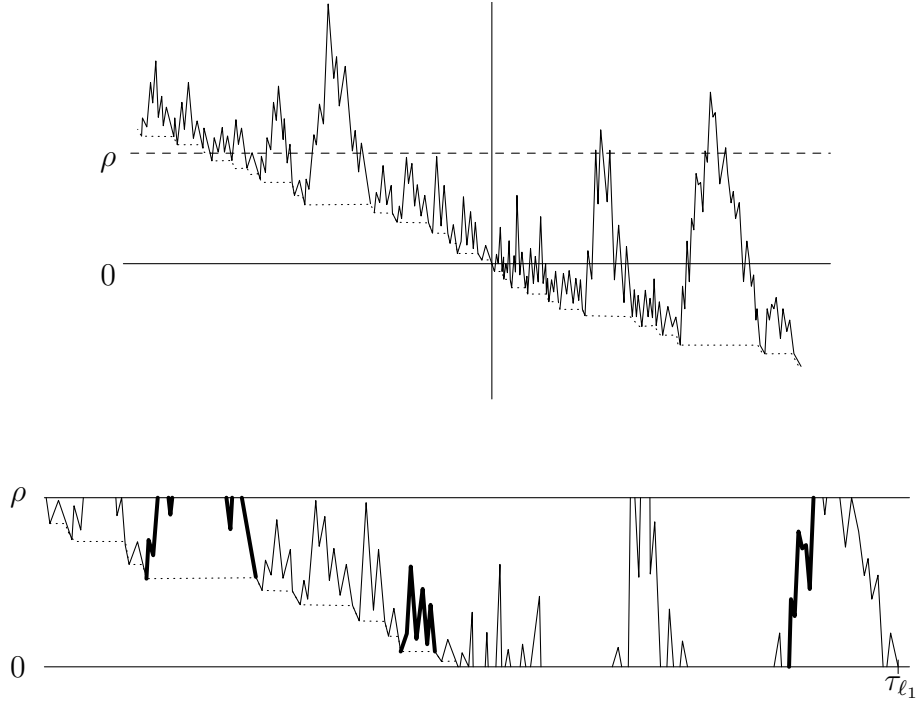

Fix $\rho>0$ and consider the interval $[0,\rho]$ in $\mathbb R$.  For $\ell_1,\ell_2 \ge 0$, the spatial Markov property implies that we can make sense of the loop soup on $[0,\rho]$ conditioned on having local time $\ell_1$ and $\ell_2$ respectively at $0$ and $\rho$, which we already denoted by ${\mathcal B}(\rho,\ell_1,\ell_2)$. One can take any $\rho_1<0$, $\rho_2>\rho$ and define it by looking at the loop soup inside $[\rho_1,\rho_2]$ with sink $\{\rho_1,\rho_2\}$, then condition it on having local times respectively $\ell_1$ at position $0$ and $\ell_2$ at position $\rho$. The loop soup ${\mathcal B}(\rho,\ell_1,\ell_2)$ is made of excursions away from $0$ which do not hit $\rho$, excursions away from $\rho$ which do not hit $0$, loops entirely contained in  $(0,\rho)$ and crossings of $[0,\rho]$. See Figure \ref{f:muprocess}. Conditioned on the number, say $2n$, of crossings (it is necessarily even), the loop soup is composed of $2n$ Brownian crossings, a loop soup of intensity ${1\over 2}$ in $[0,\rho]$ with sink $\{0,\rho\}$, and Brownian excursions away from $0$, resp. $\rho$, up to local time $\ell_1$, resp. $\ell_2$, conditioning on not hitting the other boundary. See Pitman and Yor \cite{pitmanyor} for this statement in case of squared Bessel bridge decompositions.

The local time field of the conditioned loop soup ${\mathcal B}(\rho,\ell_1,\ell_2)$ is a squared Bessel bridge of dimension $1$ from $\ell_1$ to $\ell_2$ of duration $\rho$.  If $\ell_1\ell_2>0$ the bridge hits $0$  with probability strictly between $0$ and $1$, and we denote by ${\mathcal B}^+(\rho,\ell_1,\ell_2)$ the loop soup ${\mathcal B}(\rho,\ell_1,\ell_2)$ conditioned on staying positive. When $\ell_1\ell_2=0$,  it almost surely hits $0$ on $(0,\rho)$, but taking limits $\ell_1\searrow$ or/and $\ell_2\searrow 0$, one can extend the definition of ${\mathcal B}^+(\rho,\ell_1,\ell_2)$ to this case. The loop soup ${\mathcal B}^+(\rho,\ell_1,\ell_2)$ when $\ell_1\ell_2=0$ has no crossing (having a crossing in order to have positive local time is too costly). By duality of Bessel processes, see chap. XI in \cite{revuz-yor}, the local time field of ${\mathcal B}^+(\rho,\ell_1,\ell_2)$ is a ${\rm BESQ}^3$ bridge from $\ell_1$ to $\ell_2$ of duration $\rho$ (in our case, the duality just says that a Brownian bridge conditioned on not touching $0$ is a Bessel(3) bridge).

 \bigskip

The problem of conditioning a one-dimensional Brownian loop soup on its occupation field is addressed in \cite{eyzburglar}. It is actually stated there in terms of a mu-process.  Let us briefly present the results. From now on, we consider the mu-process $X$  up to the time when the local time at position $0$ is equal, say, to $\ell_1$. Let $\widehat \L(u)$ denote the accumulated local time of $X$ at position $u$. The local time flow of  $X$  is the collection of processes $(Y_{u,v}(x),\, 0\le u\le v,\, 0\le x\le \widehat \L(u))$ giving the value at position $v$ of  the local time of $X$ at the time when the local time at position $u$ is equal to $x$. The idea of looking at the local time flow originates in a paper of T\'oth and Werner \cite{toth-werner} who construct the true self-repelling motion via a version of the Brownian web. We refer to \cite{eyzburglar} for a more detailed description of this flow in our case.  The  local time flow of $X$ is  actually a critical Feller continuous-state branching process (CSBP) flow with immigration, the immigration coming from loops encountered when exploring the positive half-line in the upwards direction. See Section 4 of \cite{pitmanyor}, also \cite{bertoin-legall00}, for the embedding of the ${\rm BESQ}^0$ flow in the  Brownian motion (i.e. when $\mu=1$).  Dawson and Li \cite{dawson-li12} construct the CSBP with immigration as the solution of the SDE
$$
Y_{u,v}(x) = x+ 2 \int_u^v  {\mathcal W}([0,Y_{u,r}(x)],\d r ) + {2\over \mu} (v-u)
$$

\noindent where ${\mathcal W}$ is a space-time white noise (and we recall that $\mu=2$ for us). A proof that the local time flow of $X$ satisfies this equation can be found in \cite{eyzburglar} by means of Tanaka's formula. A version of the Perkins' disintegration theorem stated there says that the image of this flow under a transformation which is measurable with respect to the accumulated occupation field turns it into a Fleming--Viot flow with immigration, independent of the accumulated occupation field. More precisely, if one defines
$$
\widetilde Y_{\widetilde u,\widetilde v}(\widetilde x)  := {Y_{u,v}(x) \over \widehat \L(v)}
$$

\noindent where $\widetilde x = {x\over \widetilde \L(u)}$, $\widetilde u=\int_0^u{\d r\over \widehat \L(r)} $ and $\widetilde v=\int_0^v{\d r\over \widehat \L(r)} $, then the flow $(\widetilde Y_{\widetilde u,\widetilde v}(\widetilde x),\, 0\le \widetilde u\le \widetilde v,\, \widetilde x\in [0,1])$ is a solution of the SDE
$$
\widetilde Y_{\widetilde u,\widetilde v}(\widetilde x) = \widetilde x + 2 \int_{r=\widetilde u}^{\widetilde v}  \int_{y=0}^1 \left({\bf 1}_{\{y \le \widetilde Y_{\widetilde u,\widetilde r}(\widetilde x)\}} - \widetilde Y_{\widetilde u,r}(\widetilde x)\right)\widetilde{\mathcal W}(\d y,\d r ) + {2\over \mu} \int_{\widetilde u}^{\widetilde v} (1-\widetilde Y_{\widetilde u,r}(\widetilde x))\d r
$$

\noindent where $\widetilde {\mathcal W}$ is a space-time white noise. See \cite{dawson-li12} where it is indeed defined  as a particular case of a Fleming--Viot with immigration.  Observe that the transformation dilates the space locally around $u$ by a factor ${1\over \widehat \L(u)}$. This flow is called a Jacobi($\frac2{\mu},0$) flow in \cite{eyzburglar}.

Consider now the Brownian loop soup with law ${\mathcal B}(\rho,\ell_1,\ell_2)$, with $\ell_1\ell_2>0$.  This loop soup can be traced via the mu-process $X$ conditioned on having local time $\ell_1$ and $\ell_2$ respectively at positions $0$ and $\rho$, restricted to times when it lies in $[0,\rho]$, see Figure \ref{f:muprocess}. Let $X^{(\rho)}$ denote this process\footnote{More precisely, take a mu-process $X$, and condition it on having local time $\ell_2$ at position $\rho$ at the time when the local time at position $0$ hits $\ell_1$. Restrict the process to the times when it is in $[0,\rho]$ and close the time gaps.}. From the independence of $\widetilde Y$ and $\widehat \L$ stated in the former paragraph, we deduce that the local time flow of the process $X^{(\rho)}$ will be sent by the same transformation to the same Fleming--Viot flow but restricted to the corresponding space window: in other words, its image is distributed as
$$
(\widetilde Y_{\widetilde u,\widetilde v}(\widetilde x),\, 0\le \widetilde u\le \widetilde v\le \widetilde \rho,\, \widetilde x\in [0,1])
$$

\noindent where $\widetilde \rho:= \int_0^{\rho} {\d r\over \widehat \L(u)}$ (and $\widehat \L$ denotes now the total local time of $X^{(\rho)}$).

 When one interprets  this Fleming--Viot  flow as the local time flow of some process $\widetilde X^{(\rho)}$, and draws from this process loops, excursions and crossings, one will be able to reconstruct the image of the  loop soup ${\mathcal B}(\rho,\ell_1,\ell_2)$ through a space-time scaling,  with the image being independent of $\widehat \L$. The process $\widetilde X^{(\rho)}$ is related to $X^{(\rho)}$ by the identity:
 \begin{equation}\label{def:Xtilde}
 \widetilde X^{(\rho)}_{\widetilde t} := \int_0^{X^{(\rho)}_t} {\d r \over \widehat \L(r)}
 \end{equation}
 
 \noindent where $\widetilde t = \int_0^t {\d s \over \widehat \L(X_s^{(\rho)})^2}$. It is a generalization of the definition of the so-called burglar process by Warren and Yor \cite{warren-yor}. We can consider $\widetilde X^{(\rho)}$ as the {\it contour function of a Fleming--Viot process with immigration}, or equivalently the {\it contour function of the {\rm Jacobi($\frac2{\mu},0$)} flow}, on the space interval $[0,\widetilde \rho]$. By inverting the transformation, it gives the law of the Brownian loop soup ${\mathcal B}(\rho,\ell_1,\ell_2)$ conditioned on its occupation field $\widehat \L$.

There is an issue with the previous definitions. We implicitly assumed that the local time field $\widehat \L$ does not vanish. In the case it hits $0$,  the previous operation will only reconstruct the Brownian loop soup inside a cluster. One should then iterate the operation in each cluster to   recover the whole loop soup. 


\bigskip

We can extend the previous considerations to the loop soup denoted by ${\mathcal C}(\rho,\ell_1,\ell_2)$ in the introduction. It can be constructed by taking a loop soup in $[\rho_1,\rho_2]$ with $\rho_1<0$ and $\rho_2>\rho$ with sink $\{\rho_1,\rho_2\}$, then add a Brownian crossing from $\rho_1$ to $\rho_2$ (or equivalently  a Bessel(3) process starting at $0$  stopped upon hitting $\rho_2-\rho_1$ translated to the starting position $\rho_1$), and  condition the obtained loop soup on having local time $\ell_1$, resp. $\ell_2$, at position $0$, resp. $\rho$. Again it does not depend on the choice of $\rho_1$ and $\rho_2$. The local time field of ${\mathcal C}(\rho,\ell_1,\ell_2)$ is a squared Bessel bridge of dimension $3$. Conditionally on the number, say, $2n+1$ of crossings, the loop soup ${\mathcal C}(\rho,\ell_1,\ell_2)$ is made of $2n+1$ Brownian crossings of the interval, a loop soup of intensity ${1\over 2}$ in $[0,\rho]$ with sink $\{0,\rho\}$, and Brownian excursions away from $0$ and $\rho$ up to local time $\ell_1$ and $\ell_2$ conditioned on not hitting the other boundary. In the case $\ell_1\ell_2=0$, ${\mathcal C}(\rho,\ell_1,\ell_2)$ comprises only one crossing almost surely (which is the extra crossing). Let us describe the law of the loop soup   conditioned on its accumulated occupation field.

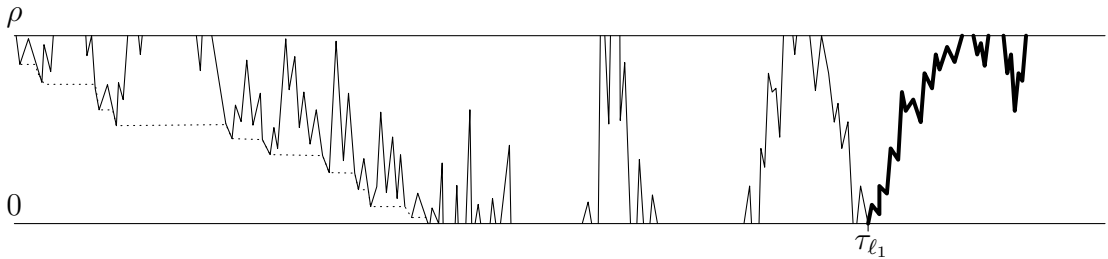
\begin{figure}[h!]
\begin{tikzpicture}[line cap=round,line join=round,>=triangle 45,x=0.5cm,y=0.5cm]

\clip(-12.5,-1.) rectangle (18,5.9);

\draw[color=black] (-11.,0.) -- (18.,0.);
\foreach \x in {-16.,-14.,-12.,-10.,-8.,-6.,-4.,-2.,2.,4.,6.,8.,10.,12.,14.,16.,18.}
\draw[shift={(\x,0)},color=black] (0pt,-2pt);

\draw (-11.,5.) -- (18.,5.);
\node at (-11.,5.5) {$\rho$};
\node at (-11.,0.5) {$0$};
\draw (11.7,0.2) -- (11.7,-0.2);
\node at (11.8,-0.6) {$\tau_{\ell_1}$};

\draw [dotted] (0.,0.)-- (-0.06659852520710975,0.16951021971730434)-- (-0.44043911647646866,0.15959451741861955)-- (-0.6166147630670642,0.4521844660194193)-- (-1.5268556037851408,0.4521844660194193)-- (-1.5897645972405297,0.8973626448783579)-- (-1.8498442892012326,0.9043689320388371)-- (-1.9502772770117538,1.3425774307763918)-- (-2.6279533949763625,1.356553398058255)-- (-2.804129041566958,1.8087378640776726)-- (-4.1988529104091725,1.835336950314109)-- (-4.404391164764867,2.234323243860654)-- (-5.211862878305096,2.2609223300970904)-- (-5.373357221013142,2.6333095374071993)-- (-8.28678105498576,2.605788121073672)-- (-8.38678105498576,3.025)-- (-8.746781054985764,3.0257881210736723)-- (-8.846781054985769,3.705788121073672)-- (-10.206781054985772,3.705788121073672)-- (-10.366781054985772,4.005788121073673);
\draw [dotted] (-10.35989493282946,3.9885728156828923)-- (-10.357505823059531,3.955651056539124)-- 
(-10.45,4.235651056539124)--
(-10.853857429864824,4.235651056539124);

\draw (-10.952717345070788,5.)-- (-10.853857429864824,4.235651056539124)-- (-10.624827436813906,4.920830953740727)-- (-10.268096308293375,3.7590762543426806)-- (-10.211983400116655,4.7612364363221085)-- (-10.033084317547832,4.043061107938327)-- (-9.96427697809829,5.);
\draw  (-9.1,5.)-- (-9.069781565254264,4.468646487721309)-- (-8.94592835424509,5.)-- (-8.849385310579928,3.5908766419189075)-- (-8.746781054985764,3.0257881210736723)-- (-8.52375857837244,3.697272986864653)-- (-8.28678105498576,2.605788121073672);
\draw (-8.28678105498576,2.605788121073672)-- (-8.230332023969858,3.7504711593375273)-- (-8.106478812960685,3.2982866933181096)--
(-7.9862953267221485,5.);

\draw  (-7.7,5.) -- (-7.65,4.5)--(-7.6,5.);
 
\draw  (-6.16020097364344,5.)-- (-6.054643950575354,4.069660194174763)-- (-6.00790739566179,5.); 

\draw (-5.75,5.)--(-5.377830430789857,2.659908623643638);
\draw (-5.377830430789857,2.659908623643638)--(-5.21,2.26);
\draw (0.05,0.)--(0.10590595728840583,0.4163590380494338)-- (0.28216184800737354,0.)-- (0.37536349221944526,1.6119263641665778)-- 
(0.4,0.);
\draw (0.7210769178694524,0.)-- ((0.765115822560836,1.0141427011080058)-- 
(0.8,0);
\draw (1.0054106179711872,0.)-- (1.1125037691694668,3.027729776673722)-- (1.1595969203677464,0.);

\draw (1.25,0.)-- (1.3327985645798182,0.5107459322165767)-- 
(1.4,0.) ;

\draw (1.65,0.)-- (1.7140780181746569,0.6680574224951483)-- 
(1.8,0.);

\draw  (1.9258999368384562,0.)--  (2.1631404857419114,2.0838608350022927)-- (2.203287883447087,0.);

\draw  (4.1,0.)-- (4.255941042140249,0.5736705283280071)-- 
(4.35,0.);

\draw (4.525398577071288,0.)-- (4.58,5.)
(4.7,5.)-- 
(4.8050035097075025,2.650182200005152)-- 
(4.85,5.)
(5.1,5.)--
(5.110027072583374,2.7445690941722947)-- (5.228647347035101,4.286221698902296)-- 
(5.381996388952349,0.);

\draw (5.551996388952349,0.)-- (5.618399677376492,1.7063132583337224)-- 
(5.728399677376492,0.) ;

\draw (5.9,0.)-- (5.965787623985123,0.7624443166622931)-- 
(6.1, 0.);

\draw   (8.4,0.)-- (8.55,1)--(8.6,0);

\draw (8.8,0)-- (8.85,2)--(8.95,1.5)--(9.05,4.)-- (9.15,3.5)--(9.25,3.6)--(9.35,2.3)-- (9.45,5);
\draw (9.7,5)-- (9.75,4.5) -- (9.8,5.)
(10.15,5.)--(10.299141731982235,3.5389703745400323)-- (10.462469741847735,5.)-- 
(10.65,4)--(10.8,2.7)-- (10.9,3.2)--
(11.0,2)--
(11.16436233582602,2.7079027741970356)--  (11.3,0.)(11.4,0.)--
(11.5,1)-- (11.7,0); 
\draw  (-5.211862878305096,2.2609223300970904)-- (-5.121236338564625,3.1535789688965736)-- (-4.9709204288879105,2.7131067961165733)-- (-4.820604519211196,4.349146295013718)-- (-4.655257018566811,2.6187199019494303)-- (-4.459846335987082,3.468201949453717)-- (-4.404391164764867,2.234323243860654)-- (-4.1988529104091725,1.835336950314109)-- (-4.099088152762968,2.555795305838002)-- (-4.,2.)-- (-3.7834247424418677,4.915467660016575)-- (-3.693235196635839,3.5625888436208597)-- (-3.542919286959125,4.443533189180861)-- (-3.4226665592177534,2.555795305838002)-- (-3.302413831476382,3.4996642475094313)-- (-3.1370663308319964,2.1153231330580016)-- (-2.9717188301876107,3.0277297766737163)-- (-2.804129041566958,1.8087378640776726)-- (-2.6279533949763625,1.356553398058255)-- (-2.445613146319111,4.852543063905147)-- (-2.2502024637393823,1.674850960278001)-- (-2.114918145030339,3.468201949453717)-- (-1.9502772770117538,1.3425774307763918)-- (-1.8498442892012326,0.9043689320388371)-- (-1.7069312181661775,1.7284302452213545)-- (-1.5268556037851408,0.4521844660194193);

\draw   (-1.5268556037851408,0.4521844660194193) --(-1.363338596646768,0.9826804030522861)-- (-1.2581174598730678,2.964805180562288)-- (-1.1078015501963536,0.8253689127737145)-- (-0.9424540495519679,2.3040969213922873)-- (-0.8222013218105965,0.668057422495143)-- (-0.7320117760045679,1.8321624505565726)-- (-0.6166147630670642,0.4521844660194193); 
\draw (-0.44043911647646866,0.15959451741861955)-- (-0.27389331607556355,0.8083196876987495)-- (0.,0.);

\draw [line width=1.5pt] (11.7,0) -- (11.8,0.5)--(12.,0.25)--(12,1)--(12.2,0.8)--(12.3,2)--(12.5,1.7)--(12.6,3.5)--(12.7,3.)--(12.9,3.3)--(13.1,2.7)--(13.2,4.)--(13.4,3.6)--(13.5,4.5)--(13.6,4.1)--(13.8,4.7)--(14.,4.3)--(14.2,5.) (14.5,5.) -- (14.6,4.5)-- (14.7,4.8)--(14.8,4.2)--(14.9,5.) (15.3,5.)--(15.4,4.)--(15.5,4.5)--(15.6,3.)--(15.7,4.)--(15.8,3.8)--(15.9,5.);

\end{tikzpicture}

\caption{Construction of ${\mathcal C}(\rho,\ell_1,\ell_2)$. After the mu-process reaches local time $\ell_1$ at position $0$, we concatenate it with a Bessel(3) process. We look at this process on the space interval $[0,\rho]$. The Bessel(3) process will give an extra crossing and excursions away from $\rho$ (in bold in the picture). Conditioning on the local time at position $\rho$ being $\ell_2$, we get the loop soup ${\mathcal C}(\rho,\ell_1,\ell_2)$. }
\label{f:bessel}
\end{figure}

Let now $X$ denote the mu-process up to local time $\ell_1$ at position $0$, followed by a Bessel(3) process starting at $0$. Again $X^{(\rho)}$ is defined as the process $X$ conditioned on having local time $\ell_2$ at position $\rho$, looked only in the space interval $[0,\rho]$. The loop soup ${\mathcal C}(\rho,\ell_1,\ell_2)$ can be defined via $X^{(\rho)}$, see Figure \ref{f:bessel}. 

Suppose $\ell_1\ell_2>0$. We still write $\widehat \L$ for the local time field of $X^{(\rho)}$ and define $\widetilde X^{(\rho)} $ by \eqref{def:Xtilde}. Then the process $\widetilde X^{(\rho)} $ is independent of $\widehat \L$. The loops, excursions and crossings traced by $\widetilde X^{(\rho)}$ can be interpreted as the loop soup conditioned on having local time $1$ at every position. Inverting the transformation \eqref{def:Xtilde} gives the law of the conditioned version of the loop soup ${\mathcal C}(\rho,\ell_1,\ell_2)$ given $\widehat \L$. The local time flow $\widetilde Y$ of  $\widetilde X^{(\rho)}$ is now solution of 
$$
\widetilde Y_{\widetilde u,\widetilde v}(\widetilde x) = \widetilde x + 2 \int_{r=\widetilde u}^{\widetilde v}  \int_{y=0}^1 \left({\bf 1}_{\{y \le \widetilde Y_{\widetilde u, r}(\widetilde x)\}} - \widetilde Y_{\widetilde u,r}(\widetilde x)\right)\widetilde {\mathcal W}(\d y,\d r ) +  \int_{\widetilde u}^{\widetilde v}{2\over \mu} (1-\widetilde Y_{\widetilde u,r}(\widetilde x)) - 2 \widetilde Y_{\widetilde u,r}(\widetilde x)\d r
$$

\noindent which we can  interpret as a Fleming--Viot process with two immigrants. The flow $\widetilde Y$ is called in \cite{eyzburglar} a Jacobi($\frac2{\mu},2$) flow, and $\widetilde X^{(\rho)}$ can be considered as its contour function on the space interval $[0,\widetilde \rho]$.


\section{Construction and properties of the measure $\n$}
\label{s:n}
In the setting of Section \ref{s:model}, consider a metric graph $\G=(V,E)$ to which we attached the distinct vertex $\mathfrak v$. Recall that by convention, $\mathfrak v\notin V$ and edges adjacent to $\mathfrak v$ are not in $E$.  Let $W\subset V$ be a set of vertices.  We call the vertices in $W$ star vertices. Let $P$ be the probability measure under which we defined the loop soup $\mathcal L$. Our goal is to define another probability measure $\n$ under which the local time at star vertices is $0$ but edges adjacent to a star vertex have positive local time fields. Therefore, a  graph and a set of star vertices will be naturally associated to a measure $\n$.

We denote by $\partial W$ the set of vertices $v\in V$ which are not in $W$ but have a neighbor in $W$. We denote by $E(W)$ the set of edges $e\in E$ which are adjacent to a vertex in $W$ and write $\mathcal E_W$ for the event that the local time field on all edges in $E(W)$ does not hit $0$. (Later, the graph $\G$ will represent our star graph, the set $W$ the set of star vertices, the set $\partial W$ the set of replicas and $E(W)$ the set of star edges.)

For a discrete loop configuration $\overline L$, recall that $n_{\overline L}(e)$ is the number of crossings of the edge $e$ and $k_v=k_v(\overline L)$ is the number of departures from $v$. We let $E^0(W):=\{e\in E(W)\,:\, n_{\overline L}(e)=0\}$ be the set of edges without any crossing, and, for $v\in W\cup\partial W$, we denote by $d_v^{W,0}$  the number of directed edges rooted at $v$ in $E(W)$ which are not crossed by any loop in $\overline L$. Recall that $d_v$ is the number of directed edges rooted at $v$ (we do not count edges adjacent to ${\mathfrak v}$).

\begin{proposition}\label{p:limit}
  Let $\overline L$ be a discrete configuration  of loops. 
\begin{enumerate}[(1)]
\item \label{i:limit1} The limit
\begin{equation}\label{eq:limit}
\lim_{(x_v)_{v \in W}\to 0}  \left(\prod_{v \in W} {1\over x_v^{(d_v+1)/2}}\right) 
P \left({\mathcal E}_W, \widehat \L(v) \le x_v \, \forall \, v\in W, \overline \L =\overline L \right)
\end{equation}

\noindent exists. 
\item \label{i:limit2}If $n_{\overline L}(e)\ge 2$ for some $e\in E(W)$, then the limit is $0$. 
\item \label{i:limit3} The law of $(\widehat \L(v),\, v\in V\backslash W)$ conditionally on $\{{\mathcal E}_W, \widehat \L(v) \le x_v \, \forall \, v\in W, \overline \L =\overline L \}$, converges as $(x_v)_{v \in W}\to 0$ to that of independent random variables, with  $\widehat \L(v)$ being gamma$\left(k_v+{d_v^{W,0}+1\over 2}, {a_v\over 2}\right)$ distributed. 
\end{enumerate}
\end{proposition}

\bigskip

\begin{remark} \label{r:limit} We will show in the course of the proof, see \eqref{eq:limit3}, that for $\overline L$ such that $n_{\overline L}(e)\le 1$ for all $e\in E(W)$, and $F$ any bounded measurable function,
\begin{eqnarray*}
 && \lim_{(x_v)_{v \in W}\to 0}  \left(\prod_{v \in W} {1\over x_v^{(d_v+1)/ 2}}\right) 
E \left[F(\widehat \L(v),\,v \in V\backslash W),\, {\mathcal E}_W, \widehat \L(v) \le x_v \, \forall \, v\in W ,\, \overline \L =\overline L\right] \\
 &=&
E\left[{\bf 1}_{\{\overline \L=\overline L\}}  F(\widehat \L(v),\,v \in V\backslash W) \prod_{v \in \partial W}  \left( {\widehat \L}(v) \right)^{d_v^{W,0}/2} \right]
 \prod_{v \in W}  {\left({a_v\over 2}\right)^{k_v+{1\over 2}} \over ({d_v+1 \over 2})\Gamma(k_v+{1\over 2})}\prod_{e\in E^0(W)}{1\over \rho(e)}.
\end{eqnarray*}
\end{remark}

\bigskip

{\it Proof of the proposition}.  Recall that  conditionally on $\{\overline \L = \overline L\}$, the random variables $(\widehat \L(v),\, v\in V)$  are independent, and $\widehat \L(v)$ is gamma$(k_v+{1\over 2}, {a_v\over 2})$ distributed.  Let $F$ be a bounded measurable function.

Observe that the local time field does not hit $0$ on an edge which is crossed by a loop. Consider now $e\in E^0(W)$ with endpoints $e_1\in W$ and $e_2$. We may have $e_1=e_2$ if it is a self-loop and we may have $e_2\in W$ if the edge joins two vertices in $W$. Conditionally on $\widehat \L(e_1)=a$, $\widehat \L(e_2)=b$ and $e\in E^0(W)$, the trace of the loop soup on the edge is composed of Brownian excursions from $e_1$ up to local time $a$, from $e_2$ up to local time $b$, all excursions being conditioned on not crossing the edge, and a Brownian loop soup inside the edge. The conditional probability that the local time field does not hit $0$ on the edge can be written 
as $q\left( {\sqrt{ab}\over \rho(e)}\right)$ for  $q(x)=1-e^{-x}$, see Lupu \cite{lupu18}, Corollary 3.6. We deduce that
\begin{align}
& E\left[F(\widehat \L(v),\, v \in V\backslash W),\, {\mathcal E}_W,\, \widehat \L(v)\le x_v \, \forall\, v\in W \, \big| \, \overline \L = \overline L\right]\label{eq:bound1} \\
&=
E\left[ F(Y_v,\,v \in V\backslash W) \,\prod_{e \in E^0(W)} q\left( {\sqrt{Y_{e_1}Y_{e_2}}\over \rho(e)}\right),\, Y_v \le x_v\, \forall \,  v\in W\right] \nonumber
\end{align}

\noindent  where $(Y_v,\,v\in V)$ are independent and gamma$(k_v+{1\over 2},{a_v\over 2})$ distributed. Since $q(x) \le x$ for $x\ge 0$, we have on the event $\{Y_v \le x_v\, \forall\, v \in W\}$,
\begin{equation}\label{eq:bound2}
\prod_{v\in W} {1\over x_v^{{d_v^{W,0} / 2}}} \prod_{e \in E^0(W)} q\left( {\sqrt{Y_{e_1}Y_{e_2}}\over \rho(e)}\right)
\le
 \prod_{e \in E^0(W),e_2\notin W}  {\sqrt{Y_{e_2}} \over \rho(e)}
  \prod_{e \in E^0(W),e_2 \in W}  {1  \over \rho(e)}.
\end{equation}

\noindent  We use the following fact: if $g$ is a gamma($p,\alpha$) random variable, then

(i) ${g \over x}$ conditioned on $g \le x$ converges in distribution as $x\searrow 0$ to a beta$(p,1)$ random variable, (the beta$(p,q)$ distribution has density $\frac{\Gamma(p+q)}{\Gamma(p)\Gamma(q)}x^{p-1}(1-x)^{q-1} {\bf 1}_{(0,1)}(x)$ for $p,q>0$)

(ii) $\lim_{x\searrow 0}{1\over x^p}P(g \le x) = {\alpha^p \over \Gamma(p+1)}$.

\bigskip

\noindent Dominated convergence and (i) imply that 
\begin{eqnarray*}
&& \lim_{(x_v)_{v\in W} \to 0} \left( \prod_{v\in W} {1\over x_v^{{d_v^{W,0} / 2}}}\right)E\left[F(Y_v,\, v \in V\backslash W) \prod_{e \in E^0(W)} q\left( {\sqrt{Y_{e_1}Y_{e_2}}\over \rho(e)}\right)\, \bigg| \, Y_v \le x_v,\, \forall  v\in W\right] \\
&=&
E\left[ F(Y_v,\, v \in V\backslash W)  \prod_{e \in E^0(W)}  {\sqrt{Z_{e_1}  Z_{e_2}}\over \rho(e)}\right]
\end{eqnarray*}

\noindent  where $(Z_v, \, v\in V)$ are independent,  $Z_v$ being beta($k_v+{1\over 2},1$) distributed if $v\in W$ and $Z_v := Y_v$ if $v  \in V\backslash W$. From (ii), we have that 
\begin{equation}\label{eq:prooflimit}
 \lim_{(x_v)_{v\in W} \to 0} \left( \prod_{v\in W} {1\over x_v^{{k_v+ {1\over 2}}}}\right)
P\left( Y_v \le x_v\, \forall  v \, \in W
\right) 
=
 \prod_{v \in W}  {\left({a_v\over 2}\right)^{k_v+{1\over 2}} \over\Gamma(k_v+{3\over 2})}.
\end{equation}

\noindent Hence,
\begin{eqnarray*}
&& \lim_{(x_v)_{v \in W}\to 0}  \left(\prod_{v \in W} {1\over x_v^{{d^{W,0}_v\over 2}+ k_v+{1\over 2}}}\right) 
E \left[ F(\widehat \L(v),\, v \in V\backslash W)  ,\, {\mathcal E}_W, Y_v \le x_v \, \forall \, v\in W \, \big|\, \overline \L =\overline L\right] \\
&=&
E\left[ F(Y_v,\, v \in V\backslash W)  \prod_{e \in E^0(W)}  {\sqrt{Z_{e_1} Z_{e_2}}\over \rho(e)}\right]
 \prod_{v \in W}  {\left({a_v\over 2}\right)^{k_v+{1\over 2}} \over \Gamma(k_v+{3\over 2})}.
\end{eqnarray*}

\noindent We compute that 
$$
\prod_{e \in E^0(W)}  {\sqrt{Z_{e_1} Z_{e_2}}\over \rho(e)} 
=
\prod_{v \in \partial W}  \left( Y_v \right)^{d_v^{W,0}/2}
\prod_{v \in  W}  \left(  Z_v \right)^{d_v^{W,0}/2}\prod_{e\in E^0(W)}{1\over \rho(e)}.
$$

\noindent  Since, for $v\in W$,
$$
E\left[  \left( Z_v \right)^{d_v^{W,0}/2}	\right]
=
  {k_v+{1\over 2}\over k_v +{d_v^{W,0}+1 \over 2} },
$$

\noindent we get 
\begin{align*}
 \lim_{(x_v)_{v \in W}\to 0}  & \left(\prod_{v \in W} {1\over x_v^{{d^{W,0}_v \over 2}+ k_v+{1\over 2}}}\right) 
E \left[F(\widehat \L(v),\,v \in V\backslash W),\, {\mathcal E}_W, \widehat \L(v) \le x_v \, \forall \, v\in W \, \big|\, \overline \L =\overline L\right]  \\
=
&E\left[ F(Y_v,\,v \in V\backslash W) \, \prod_{v \in \partial W}  \left( Y_v \right)^{{d_v^{W,0}/ 2}} \right]
 \prod_{v \in W}  {\left({a_v \over 2}\right)^{k_v+{1\over 2}}(k_v+{1\over 2}) \over (k_v +{d_v^{W,0}+1 \over 2}) \Gamma(k_v+{3\over 2})} \prod_{e\in E^0(W)}{1\over \rho(e)}.  
\end{align*}

\noindent  Notice that for all $v\in W$, $k_v+{d_v^{W,0}+1\over 2} \ge {d_v+1\over 2}$ and we have equality if and only if every edge adjacent to $v$ is crossed by $\overline L$  at most once. It readily proves assertions \eqref{i:limit1} and \eqref{i:limit2}. Using $\Gamma(k_v+{3\over 2})=(k_v+{1\over 2})\Gamma(k_v+{1\over 2})$, we get when $n_{\overline L}(e)\le 1$ for all $e\in E(W)$, 
\begin{eqnarray}
\label{eq:limit3} && \lim_{(x_v)_{v \in W}\to 0}  \left(\prod_{v \in W} {1\over x_v^{{d_v+1\over 2}}}\right) 
E \left[F(\widehat \L(v),\,v \in V\backslash W),\, {\mathcal E}_W, \widehat \L(v) \le x_v \, \forall \, v\in W \, \big|\, \overline \L =\overline L\right] \\
\nonumber &=&
E\left[  F(Y_v,\,v\in V\backslash W) \prod_{v \in \partial W}  \left( Y_v \right)^{d_v^{W,0}/2} \right]
 \prod_{v \in W}  {\left({a_v \over 2}\right)^{k_v+{1\over 2}} \over ({d_v+1 \over 2})\Gamma(k_v+{1\over 2})}\prod_{e\in E^0(W)}{1\over \rho(e)}.
\end{eqnarray}

\noindent  From there we can deduce \eqref{i:limit3}. $\Box$

\bigskip
\bigskip

Notice that by \eqref{eq:bound1} and \eqref{eq:bound2}, for $c=\prod_{e\in E} \max(1,{1\over \rho(e)})$, for all $\overline L$,
\begin{eqnarray*}
&& 
\prod_{v\in W} {1\over x_v^{{d_v^{W,0} / 2}}}  P\left(  {\mathcal E}_W,\,\widehat \L(v)\le x_v, \, \forall\, v\in W  \, \big| \, \overline \L = \overline L\right)\\
&\le&
c\, E\left[  \prod_{e \in E^0(W),e_2\notin W}  \sqrt{Y_{e_2}}\right] P\left(Y_v \le x_v, \, \forall\, v\in W
   \right).
\end{eqnarray*}

\noindent By \eqref{eq:prooflimit} (the right-hand side of  \eqref{eq:prooflimit} is actually an upper bound), we get  for some constant $c'>0$,
\begin{eqnarray*}
&& \prod_{v\in W} {1\over x_v^{{ (d_v+1)/ 2}}}  P\left(  {\mathcal E}_W,\,\widehat \L(v)\le x_v, \, \forall\, v\in W   \, \big| \, \overline \L = \overline L\right) \\
&\le&
c' \, E\left[   \prod_{e \in E^0(W),e_2\notin W}  \sqrt{Y_{e_2}}\right]\\
&=&
c' \, E\left[   \prod_{e \in E^0(W),e_2\notin W}  \sqrt{\widehat {\mathcal L}(e_2)}\, \bigg| \, \overline \L = \overline L \right].
\end{eqnarray*}

\noindent We have $E\left[\prod_{e \in E^0(W),e_2\notin W}  \sqrt{\widehat {\mathcal L}(e_2)}\right]<\infty$ since for $v\in V$, the variable $\widehat {\mathcal L}(v)$ is gamma distributed hence has all moments. By Proposition \ref{p:limit}, we deduce that the limit $$
 \left(\prod_{v \in W} {1\over x_v^{(d_v+1)/2}}\right) 
P \left({\mathcal E}_W, \widehat \L(v) \le x_v \, \forall \, v\in W\right)
$$ 

\noindent exists. By  Remark \ref{r:limit}, we can see that the distribution of $(\overline {\mathcal L}, (\widehat {\mathcal L}(v),\, v\in V\backslash W))$ conditionally on the event $\{{\mathcal E}_W, \widehat \L(v) \le x_v \, \forall \, v\in W\}$ converges. And by the description of  the loop soup conditionally on $\overline{\mathcal L}$ and the local times at vertices  given in Section \ref{s:model},  the loop soup  converges in law. The setting in which the convergence in distribution takes place is the following. We consider the space of loops on the metric graph $\G$. Renormalizing the time parametrization of a Brownian loop by its duration, we can define a rooted oriented loop (see for example \cite{lawler}) as a couple $(\gamma,r)$ where $\gamma$ is a continuous map from the unit circle to $\G$ and $r\in (0,\infty)$ is its duration. A unrooted unoriented loop is a couple $(\bar{\gamma},r)$ where $\bar{\gamma}$ is an equivalence class of loops $\gamma$ with respect to time-shift and time-reversal. We equip this space with a product metric (the distance between $\bar{\gamma}$ and $\bar{\gamma'}$ is the infimum of the uniform distance between $\gamma$ and $\gamma'$ over all $\gamma \in \bar{\gamma}$, $\gamma'\in \bar{\gamma'}$). A loop soup on $\G$ can be seen as a subset of the space of unrooted unoriented loops.  We equip the space of loop soups with the distance $d:=\sum_{n\ge 1} 2^{-n}\min(d_n,1)$ where $d_n$ is the Hausdorff distance when restricting to loops of duration $r\ge \frac1n$.  

\bigskip

\begin{definition}\label{def:n}
We denote by $\n$ the limiting probability measure  
$$
\lim_{(x_v)_{v \in W}\to 0}  
P \left(\cdot \,\bigg|\, {\mathcal E}_W, \widehat \L(v) \le x_v \, \forall \, v\in W \right).
$$
\end{definition}
\bigskip

We draw some consequences of  Proposition \ref{p:limit}, and usual properties of the  loop soup recalled in Section \ref{s:model}.   Let $\overline L$ be a discrete loop configuration such that $n_{\overline L}(e)\le 1$ for all $e\in E(W)$. Recall the notation ${\mathcal B}^+(\rho,\ell_1,\ell_2)$ and ${\mathcal C}(\rho,\ell_1,\ell_2)$ in Section \ref{s:loopsoupcond}. 

\bigskip

\begin{enumerate}[1)]
\item Under $\n(\cdot\, |\, \overline \L=\overline L)$, the random variables $(\widehat \L(v),\, v\in V\backslash W)$ are independent, the random variable $\widehat \L(v)$ being  gamma$\left(k_v+{d_v^{W,0}+1\over 2}, {a_v\over 2}\right)$ distributed.
\item Conditioning further on $(\widehat \L(v)=y_v,\, v\in \partial W)$, the trace of the loop soup on different edges of $E(W)$ are independent and with distribution given as follows. ($e_1$ and $e_2$ are the endpoints of $e$ with the convention that $e_1\in W$)
\begin{enumerate}
\item  The trace of the loop soup on the edge $e\in E(W)$ is distributed as ${\mathcal B}^+(\rho(e),0,y_{e_2})$ if $n_{\overline L}(e)=0$.
\item The trace of the loop soup on the edge $e\in E(W)$ is distributed as ${\mathcal C}(\rho(e),0,y_{e_2})$ if $n_{\overline L}(e)=1$.
\end{enumerate}
\item Still under $\n(\cdot\, |\, \overline \L=\overline L,\, \widehat \L(v)=y_v,\, v\in \partial W)$, the trace of the loop soup outside $E(W)$  is distributed as under $P$ conditioned on the same event.
\end{enumerate}


\bigskip

\begin{theorem}\label{t:markovn}
Under $\n$, conditionally on $(\widehat \L(v),\, v \in \partial W)$, the trace of the loop soup outside $E(W)$  is independent  of the occupation  field in $E(W)$.
\end{theorem}

\noindent {\it Proof}. Conditionally on $(\widehat \L(v),\, v \in \partial W)$, the discrete loop soup $\overline \L$ is independent of the occupation field in $E(W)$.   
Indeed, conditionally on $\{\overline \L =\overline L\}$ and on the local times at vertices in $\partial W$, the  occupation fields on distinct edges of $E(W)$ are independent and distributed as ${\rm BESQ}^3$ bridges. Their (conditional) distribution does not depend on $\overline L$ hence the independence.
Then, conditionally on $(\widehat \L(v),\, v \in \partial W)$ and on $\overline \L$, the trace of the loop soup outside $E(W)$ is independent of the trace of the loop soup on $E(W)$ (it is true under $P$). $\Box$

\bigskip

 We give now the conditional law of the trace of the loop soup outside $E(W)$ when $(W,E(W))$ is connected meaning that any two vertices of $W$ can be joined by a path staying in $E(W)$.  The following result is the analog of the domain Markov property for the Gaussian free field. For any pair $v\neq v'\in \partial W$, we let $n_{\mathcal L}(\{v,v'\})$ stand for the number of excursions  outside $E(W)$ which  join $v$ and $v'$  and  $H_{\overline{W}^c}(v,v')$ be the mass under the excursion measure at $v$ of excursions which hit $v'$ before hitting any other point in $W\cup \partial W\cup\{\mathfrak v\}$, which is also $H_{\overline{W}^c}(v',v)$. In the statement of the theorem and in its proof, we write for short $v\neq v'\in \partial W$ whereas we rigorously index  with the set of pairs $\{\{v,v'\},\, v\neq v' \in \partial W\}$. The proof of the following theorem  follows the lines of the proof of the Markov property  by Werner \cite{wernermarkov}, see also Powell and Werner \cite{pw20}, together with an extra combinatorial argument.

\begin{theorem}\label{t:markovn2}
 Suppose that $(W,E(W))$ is connected. Under $\n$, conditionally on $(\widehat \L(v)=x_v,\, v \in \partial W)$, the random variables $(n_{\mathcal L}(\{v,v'\}),\, v\neq v'\in \partial W)$  are independent, Poisson distributed with respective mean $2H_{\overline{W}^c}(v,v')\sqrt{x_v x_{v'}}$. 
\end{theorem}
\noindent {\bf Remark}. To complete the picture of Theorem \ref{t:markovn2}: 
\begin{itemize}
\item independently at each $v\in \partial W$, one draws excursions away from  $E(W)$ conditioned on hitting back $\partial W$ at $v$ up to local time $x_v$.  It can be proved by the property that loops  which visit $v\in V$ under $P$  conditionally on $\widehat L(v)=x_v$ trace Brownian excursions away from $v$ up to local time $x_v$ (conditioned on not hitting $\mathfrak v$). By excursion theory, excursions outside $E(W)$ which leave and return to $v$ are independent of excursions which hit $E(W)$. By theory of Poisson point processes, they are also independent of the loops which do not visit $v$. We deduce  that  the collection of such excursions at $v$   are independent for $v\in \partial W$(conditionally on the local times at $v\in \partial W$), and independent of the occupation field in $E(W)$. From Definition \ref{def:n}, it is also true under $\n$. 
\item  Loops which stay outside $E(W)$ form a loop soup with sink $\partial W$.
It comes from the independence property of Poisson point processes.  
\item An excursion outside $E(W)$ from $v\in \partial W$ to $v'\in \partial W\backslash\{v\}$ has the law of an excursion from $v$ stopped at the hitting time of $v'$ under the excursion measure at $v$ conditioned on hitting back $E(W)$ at $v'$. Conditionally on $(n(\{v,v'\}),\, v\neq v'\in \partial W)$, all excursions joining $v\neq v' \in \partial W$ are independent.
\end{itemize}
\bigskip

\noindent {\it Proof of the theorem}. We take  $(n(e),\, e\in E(W))$ a collection of integers in $\{0,1\}$ and  $(n(\{v,v'\}),\, v\neq v' \in \partial W)$ a collection of integers in $\{0,1,2,\ldots\}$, which verify that for all $v\in W$, the sum of $n(e)$ over all directed edges rooted at $v$ is even (equal to, say, $2k_v$), and for all $v\in \partial W$, the sum of $n(e)$ and $n(\{v,v'\})$ over directed edges $e$ rooted at $v$ and vertices $v'$ in $\partial W\backslash\{v\}$ is even (equal to, say, $2k_v$).  We let $A$ be the event that for all $e \in E(W)$, the number of crossings of the edge $e$ is $n(e)$,  and for all pairs $v\neq v' \in \partial W$, $n_{\mathcal L}(\{v,v'\})=n(\{v,v'\})$.  We will first show that the conditional probability $\n(A\,|\, \widehat \L(v)=x_v,\, v \in \partial W)$ is
\begin{equation}\label{eq:na2}
\frac1{Z}\prod_{v\neq v' \in \partial W} \frac{(2H_{\overline{W}^c}(v,v')\sqrt{x_v x_{v'}})^{n(\{v,v'\})}}{n(\{v,v'\})!}
\end{equation}

\noindent where $Z$ is the renormalizing constant (we do not  need that $(W,E(W))$ is connected to prove this equation). Equation \eqref{eq:na2} is  of course reminiscent of  Proposition 7 of \cite{wernermarkov}.  We give the proof to be self-contained. As in Section \ref{s:model}, a key step  (duplicating the graph, see below) is borrowed from \cite{wernermarkov} and \cite{pw20}.  We  compute $\n(A,\, \widehat {\mathcal L}(v) \in \d x_v,\, \forall \, v\in \partial W)$.  We use Remark \ref{r:limit}. This probability is proportional to 
\begin{equation}\label{eq:nA2}
P(A,\,\widehat {\mathcal L}(v) \in \d x_v,\, \forall \, v\in \partial W ) \prod_{v\in W} \frac{a_v^{k_v}}{2^{k_v}\Gamma(k_v + \frac12)}\prod_{e\in E^0(W)} \frac{1}{\rho(e)} \prod_{v\in \partial W} x_v^{d_v^{W,0}/2}
\end{equation}
 
\noindent where $d_v^{W,0}$ is the number of edges $e$ in $E(W)$ adjacent to $v$ with $n(e)=0$. For $v \in \partial W$, we let $H_{\overline{W}^c}(v)$ be the mass under the excursion measure at $v$ of excursions which hit $W\cup \partial W\cup\{\mathfrak v\} \backslash\{v\}$. The probability that the Brownian motion starting at $v \in \partial W$ hits $W\cup \partial W\cup\{\mathfrak v\} \backslash \{v\}$ at $v'$ is
\begin{equation}\label{eq:prob1}
\frac{H_{\overline{W}^c}(v,v')}{H_{\overline{W}^c}(v)}.
\end{equation}

\noindent The probability that a Brownian motion starting at $v\in \partial W$ hits $W$ before $\partial W\cup\{\mathfrak v\} \backslash \{v\}$ by crossing an adjacent edge $e\in E(W)$ is
\begin{equation}\label{eq:prob2}
\frac1{H_{{\overline W}^c}(v)}\frac1{2\rho(e)} .
\end{equation}

\noindent We compute $P(A)$. We erase the graph outside $E(W)$ and replace it by  $K$ copies of it, glued along $\partial W$. In this manner,  vertices which are not in $W\cup\partial W$ and edges outside $E(W)$ have been duplicated, and now appear $K$ times. Edges in $E(W)$ are not duplicated nor vertices of $W\cup\partial W$. We give length $K\rho(e)$ to an edge $e$ which has been duplicated. Considering the rooted loop measure, we check that projecting onto the original graph, the discrete loop soup $\overline {\mathcal L}$ keeps the same distribution (holding times at vertices outside $W\cup\partial W$ differ but they will not appear in the computations). As $K\to \infty$, with probability tending to $1$, any excursion outside $E(W)$ joining distinct vertices of $\partial W$ will take place in a distinct duplicated graph, making computations easier. 

We let $N_{\overline W}$ be the number of ways to  trace the excursions outside $E(W)$ between distinct vertices of $\partial W$ (considered globally as a bridge between two vertices of $\partial W$, without taking care of its trajectory) and locally pair them with the crossings in $E(W)$ to get the configuration $(n(e),n(\{v,v'\}),e \in E(W),v\neq v'\in \partial W)$, when no two excursions belong to the same duplicated graph.  For every excursion outside $E(W)$, we need to discuss the duplicated graph they belong to. The number of choices is 
$$
\frac{K!}{(K-n(\cdot,\cdot))!\prod_{v\neq v' \in \partial W} n(\{v,v'\})!}
$$

\noindent where $n(\cdot,\cdot)$ denotes the sum of $n(\{v,v'\})$ over all $v\neq v'\in \partial W$. Then we locally glue the crossings and the excursions at each vertex $v$, and the number of ways to do that is
$$
\prod_{v\in W\cup\partial W}(2k_v-1)\times (2k_v-3)\times \ldots \times 3\times 1 = \prod_{v\in W\cup\partial W} 2^{k_v} {\Gamma(k_v+{1\over 2}) \over \Gamma({1\over 2})}.
$$

\noindent It gives
\begin{equation}\label{eq:nw2}
N_{\overline W} = \prod_{v\in W\cup\partial W} 2^{k_v} {\Gamma(k_v+{1\over 2}) \over \Gamma({1\over 2})} \frac{K!}{(K-n(\cdot,\cdot))!\prod_{v\neq v' \in \partial W} n(\{v,v'\})!}.
\end{equation}

\noindent Let $\widetilde A$ be the event $A$ intersected with the event that every excursion outside $E(W)$ occurs on a distinct graph. Let $n(v,\cdot)$ denote the sum of $n(\{v,v'\})$ over all $v' \in \partial W$ different from $v$. Recalling \eqref{eq:prob1} and \eqref{eq:prob2}, and using the expression of the loop measure in \eqref{eq:nu}, the probability $P(\widetilde A)$ is proportional to (the term $K^{-n(\cdot,\cdot)}$ appears due to the fact that knowing that there is an excursion between some $v\in \partial W$ and $v'\in \partial W$, there is a probability $1/K$ to choose a given duplicated graph)
$$
N_{\overline W} K^{-n(\cdot,\cdot)}\prod_{n(e)=1} \frac{1}{2\rho(e)}\prod_{v\neq v' \in \partial W} H_{{\overline W}^c}(v,v')^{n(\{v,v'\})} \prod_{v\in W} \frac{2^{k_v}}{a_v^{k_v}}\prod_{v\in \partial W} \frac{1}{H_{{\overline W}^c}(v)^{k_v}}.
$$

\noindent   The local time at a vertex $v\in \partial W$ being gamma$(k_v+\frac12, H_{{\overline{W}}^c}(v))$ distributed conditionally on $\widetilde A$, we get that $P(\widetilde A,\,\widehat {\mathcal L}(v) \in \d x_v,\, \forall \, v\in \partial W )$ is proportional to 
$$
N_{\overline{W}} K^{-n(\cdot,\cdot)} \prod_{n(e)=1} \frac{1}{2\rho(e)}\prod_{v\neq v' \in \partial W} H_{{\overline W}^c}(v,v')^{n(\{v,v'\})} \prod_{v\in W} \frac{2^{k_v}}{a_v^{k_v}} \prod_{v\in \partial W} \frac{x_v^{k_v}}{\Gamma(k_v+\frac12)} \d x_v.
$$

\noindent We deduce that $P(A, \widehat {\mathcal L}(v)\in \d x_v,\, \forall \, v\in \partial W)$ is proportional to the limit of the latter as $K\to\infty$. From \eqref{eq:nA2} and \eqref{eq:nw2}, we see that $\n( A,\, \widehat {\mathcal L}(v) \in \d x_v,\, \forall \, v\in \partial W)$  is proportional to
$$
\prod_{n(e)=1} \frac{1}{2}\prod_{v\neq v' \in \partial W} \frac{H_{{\overline W}^c}(v,v')^{n(\{v,v'\})}}{n(\{v,v'\})!} \prod_{v\in \partial W} x_v^{ k_v + \frac{d_v^{W,0}}{2} } \d x_v\prod_{v\in W\cup\partial W} 2^{k_v}.
$$

\noindent By definition, $\sum_{e\in E(W)} n(e) +\sum_{v\neq v' \in \partial W} n(\{v,v'\}) = \sum_{v\in W\cup \partial W} k_v$. We deduce that 
$$
 \prod_{n(e)=1}  \frac{1}{2} \prod_{v\in W\cup\partial W} 2^{k_v} = \prod_{e\in E(W)}  \frac{1}{2^{n(e)}} \prod_{v\in W\cup\partial W} 2^{k_v}
 = \prod_{v\neq v' \in \partial W} 2^{n(\{v,v'\})}.
$$

\noindent We also have for any $v\in \partial W$, 
$$
2k_v + d_v^{W,0} = \sum_{v'\in \partial W\backslash\{v\}} n(\{v,v'\}) + \sum_{e\in E(W),\, v\in e} n(e) + d_v^{W,0} = \sum_{v'\in \partial W\backslash\{v\}} n(\{v,v'\}) + d_v^W
$$

\noindent  where $d_v^W$ is the number of neighbors of $v$ which are in $W$. From these observations,  we can deduce equation \eqref{eq:na2}.

Fix a configuration $(n(\{v,v'\}),\, v\neq v'\in \partial W)$. It remains to sum over all admissible configurations $(n(e),\, e\in E(W))$. By admissible, we mean that for all $v\in W$, the sum of the $n(e)$'s over directed edges rooted at $v$ is even and for $v\in \partial W$, the sum of the $n(e)$'s and $n(\{v,v'\})$'s over directed edges $e\in E(W)$ rooted at $v$ and vertices $v'\in \partial W\backslash\{v\}$ is also even. We can restate the last property by saying that for each $v\in \partial W$, the sum of the $n(e)$'s over edges $e\in E(W)$ adjacent to $v$ has the parity of $n(v,\cdot)$. A priori, the number of admissible configurations should depend on the configuration $(n(\{v,v'\}),\, v\neq v' \in \partial W)$, or more exactly on the parity of the $n(v,\cdot)$'s. It is actually not the case when $(W, E(W))$ is connected .

Notice that $\sum_{v\in \partial W} n(v,\cdot)$ is even. Since $(W, E(W))$ is connected, the set of admissible configurations $(n(e),\, e\in E(W))$ is not empty.  We can show it by induction. Let $E^\partial(W)$ denote the set of edges joining a vertex in $W$ and a vertex in $\partial W$. Call self-avoiding path a path which does not go twice through the same vertex. When $n(v,\cdot)$ is even for all vertices $v\in \partial W$, we can take the configuration $n(e)=0$ for all $e\in E(W)$. Suppose that $2k+2$ vertices $v\in \partial W$ verify that $n(v,\cdot)$ is odd. For each such vertex $v$, choose an adjacent edge $e\in E^\partial(W)$ and set $n(e)=1$. We join by induction $2k$ of these edges by  self-avoiding paths in $E(W)$ which do not share any edges. For the last two of the $2k+2$ edges, we trace a self-avoiding path $P$ in $E(W)$. If $P$ does not share any edge with the other ones, we set $n(e)=1$ if $e$ is crossed by a path, and $0$ otherwise. In the other case, $P$ has a  common edge with another path $P'$. Call $e_1$ and $e_2$, resp. $e_1'$ and $e_2'$, the edges of $E^\partial(W)$ joined by $P$, resp. $P'$. We follow the path $P$ from $e_1$ to $e_2$ until it intersects $P'$ at a vertex $w$. Without loss of generality, we suppose that we chose $P'$ as the first path that $P$ intersects  among paths sharing a common edge with $P$. Let $w'$ be the last vertex on the path $P$ from $e_1$ to $e_2$ which belongs to $P'$. We have $w \neq w'$ since $P'$ is self-avoiding. When removing the vertex $w'$, the path $P'$ is disconnected into two components containing respectively $e_1'$ and $e_2'$, so that $w$ can be joined to  one of the two edges $e_1'$ or $e_2'$ (say it is $e_1'$) by a subpath of $P'$ which does not visit $w'$. We can form a self-avoiding path from $e_1$ to $e_1'$ by following $P$ up to the vertex $w$ then following $P'$ up to the edge $e_1'$. Moreover, if we follow $P'$ from $e_2'$ to $w'$, then  $P$ from $w'$ to $e_2$, one has another self-avoiding path which is disjoint from the first one. The path from $e_1$ to $e_1'$ has no common edges with any other self-avoiding path joining the $2k$ other edges.  We  remove it and still have a connected graph with $2k$ edges to join, so we can  use our induction hypothesis.

We now want to show the following claim: suppose $(W,E(W))$ forms a connected graph and $(n(v,\cdot),\, v\in \partial W)$ is a collection of integers in $\{0,1,2,\ldots\}$ whose sum is even. The number of admissible configurations $(n(e),\, e\in E(W))$, in the sense that at each $v\in W$, the sum of the $n(e)$'s over directed edges $e\in E(W)$ rooted at $v$ is even, and at each $v\in \partial W$, the sum has the parity of $n(v,\cdot)$ does not depend on $(n(v,\cdot),\, v\in \partial W)$. Let us prove it. Fix a configuration $(n(v,\cdot),\, v \in \partial W)$ with even sum. Fix an admissible configuration $(n'(e),\, e\in E(W))$ which exists by the paragraph above. We define a bijection between the set of admissible configurations associated to $(n(v,\cdot),\, v\in \partial W)$ and that associated to the configuration  $n(v,\cdot)=0$ for all $v\in \partial W$.  This bijection is defined as follows: for each admissible configuration $(n(e),\, e\in E(W))$, and each $e\in E(W)$, we let $N(e):= n(e)+n'(e)$ modulo $2$. At any $v\in W\cup\partial W$, the sum of the $N(e)$'s over  edges $e\in E(W)$ rooted at $v$ is even since $n$ and $n'$ have the same parity conditions. Its inverse map is $N\to N+n'$.  It proves the claim and the theorem. $\Box$


 \section{The detachment lemma}
 \label{s:detachment}
 
\begin{figure}[h!]
a)
\begin{minipage}{0.5\textwidth}
\begin{tikzpicture}[x=0.6cm,y=0.6cm]
\clip(-1.5,0) rectangle (5,5);
\node[draw,circle, fill=black,label={[label distance=0.1cm]below:$v_0$}] (1) at (2,2) {};
\node[draw,circle, fill=black]  at (4,2) {};
\node[draw,circle, fill=black]  at (0.5,3.5) {};
\node[draw,circle, fill=black]  at (0.5,0.5) {};

\draw (1)  -- (4,2);
\draw [line width=2pt] (1)  -- (3,2);

\draw (1) -- (0.5,3.5);

\draw (1) -- (0.5,0.5);

\end{tikzpicture}
\end{minipage}
\begin{minipage}{0.5\textwidth}
\begin{tikzpicture}[x=0.6cm,y=0.6cm]
\clip(0,0) rectangle (9,5);
\node[draw,circle, fill=black,label={[label distance=0.1cm]below:$v_0$}] (1) at (2,2) {};
\node[draw,circle, fill=black]  at (4,2) {};
\node[draw,circle, fill=black]  at (0.5,3.5) {};
\node[draw,circle, fill=black]  at (0.5,0.5) {};

\draw (1) -- (0.5,3.5);
\draw (1) -- (0.5,0.5);

\node[draw,star, star points=7,fill=black,label={below:$v_0^*$}]  at (4,2) {};
\node[draw,rectangle, fill=black,label={below:$\bar{v}_0$}]  at (6,2) {};

\draw[dotted] (1) -- (4,2);
\draw[dotted] (4,2) -- (6,2);

\draw(6,2)--(8,2);
\draw [line width=2pt] (6,2)  -- (7,2);
\node[draw,circle,fill=black] at (8,2) {};
\end{tikzpicture}
\end{minipage}

b)
\begin{minipage}{0.5\textwidth}
\begin{tikzpicture}[x=0.6cm,y=0.6cm]
\clip(-1.2,-1) rectangle (6,6);
\node[draw,circle, fill=black,label={[label distance=0.1cm]below:$v_0$}] (1) at (3,1) {};

\draw (3,3) circle (2);
\draw[line width=2pt] (3,1) arc[start angle=-90, end angle=0, radius=2];
\end{tikzpicture}
\end{minipage}
\begin{minipage}{0.5\textwidth}
\begin{tikzpicture}[x=0.6cm,y=0.6cm]
\clip(0,-1) rectangle (9,5);
\node[draw,circle, fill=black,label={[label distance=0.1cm]below:$v_0$}] (1) at (3,1) {};
\node[draw,star, star points=7,fill=black,label={[label distance=0.0cm]below:$v_0^*$}] (2) at (5,1) {};
\node[draw,rectangle, fill=black,label={[label distance=0.1cm]below:$\bar{v}_0$}] (3) at (7,1) {};

\draw[dotted] (1)--(2)--(3);

\draw (1)   to [out=90,in=180] (5.,4) to  [out=0,in=90] (3);

\draw[line width=2pt] (3)   to [out=90,in=0] (5.,4);
\end{tikzpicture}
\end{minipage}
\caption{We present two examples of detachment. The beginning of the directed edge $e_0$  is drawn in bold. The star edges are represented in dotted lines. Example b) is the case of a self-loop.}
\label{f:detachment}
\end{figure}
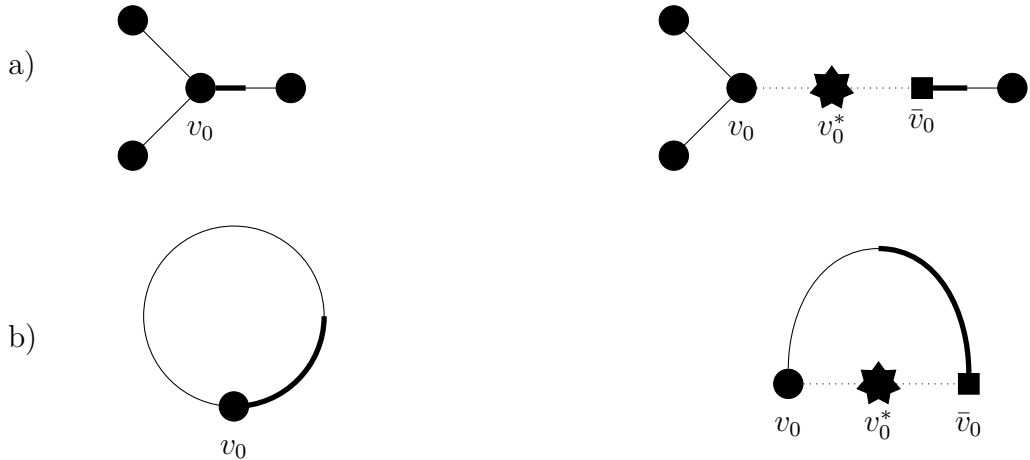

Let $v_0$ be a vertex in $V$, and $e_0\in E$ be a directed edge rooted at $v_0$. We could also take for $e_0$ a directed edge  from $v_0$ to $\mathfrak v$ (such an edge has been added to the graph when attaching $\mathfrak v$ so is not contained in $E$ in the definition of our model in Section \ref{s:model}).  We construct a new graph $\G^0$ as follows. We  detach the edge at $v_0$, thereby duplicating $v_0$. We write $\bar{v}_0$ for the vertex we created. Then we  introduce a star vertex $v_0^*$ and draw two edges between $v_0^*$ and $v_0$, and between $v_0^*$ and $\bar{v}_0$. We refer to these edges as {\it star edges}.  We naturally identify the {\it old edges} of $\G^0$ with the edges of $\G$. See Figure \ref{f:detachment}.

We can define the loop soup on $\G^0$ with sink $\mathfrak v$ as usual, under some probability measure $P^0$. We denote it by $\L^0$, and the associated local time at $v$ by $\widehat \L^0(v)$. Let $\n^0$ be the measure $\n$ constructed in Section \ref{s:n} by taking on the graph $\G^0$, $W:=\{v_0^*\}$ for the set of star vertices. Hence, $E(W)$ is only composed of the two star edges we introduced. 

In the following lemma, the trace of the loop soup outside $\{v_0\}$ comprises the loops which do not hit $v_0$ and the excursions away from $v_0$. It is important that we lose the orientation of the loops and the information on how these excursions are connected inside the loop soup.

\begin{lemma}\label{l:detachment}
Let $x>0$. The trace of the loop soup $\L$ on $\G$  outside $\{v_0\}$ under $P(\cdot\, | \, \widehat \L(v_0)=x)$  is equal in law to the trace of the  loop soup $\L^0$  outside the star edges under $\n^0( \cdot\, |\, \widehat \L^0(v_0)=\widehat \L^0(\bar{v}_0)=x)$.  
\end{lemma}

\noindent {\it Proof}. The trace of $\mathcal L$ under $P(\cdot\, | \, \widehat \L(v_0)=x)$ is that of a loop soup with sink $\{v,\mathfrak v\}$ plus  a Poisson point process of excursions away from $v_0$ up to local time $x$, conditioned on not hitting $\mathfrak v$. Lemma \ref{l:detachment} is a consequence of Theorem \ref{t:markovn2} with $W:=\{v_0^*\}$ there, and the remark following it. $\Box$

\bigskip

\begin{figure}[h!]
\begin{minipage}{0.4\textwidth}
\begin{tikzpicture}[x=0.8cm,y=0.8cm]
\clip(-1.7,0) rectangle (5,5);
\node[draw,circle, fill=black,label={[label distance=0.1cm]below:$v_0$}] (1) at (2,2) {};
\node[draw,circle, fill=black]  at (4,2) {};
\node[draw,circle, fill=black]  at (0.5,3.5) {};
\node[draw,circle, fill=black]  at (0.5,0.5) {};

\draw (1)  -- (4,2);
\draw [line width=2pt] (1)  -- (3,2);
\draw [line width=2pt] (1)  -- (1.25,2.75);
\draw [line width=2pt] (1)  -- (1.25,1.25);
\draw (1) -- (0.5,3.5);

\draw (1) -- (0.5,0.5);

\end{tikzpicture}
\end{minipage}
\begin{minipage}{0.5\textwidth}
\begin{tikzpicture}[x=0.7cm,y=0.7cm]
\clip(-5.5,-1.5) rectangle (9,6.);
\node[draw,star, star points=7, fill=black,label={[label distance=0.1cm]below:$v_0$}] (1) at (2,2) {};

\node[draw,rectangle, fill=black] (2)  at (0.5,3.5) {};
\node[draw,rectangle, fill=black]  at (0.5,0.5) {};

\node[draw,circle, fill=black]  at (-1.,5.) {};

\draw[dotted] (1) -- (0.5,3.5);
\draw[dotted] (1) -- (0.5,0.5);
\draw (2) -- (-1.,5.);

\node[draw,rectangle, fill=black]  at (4,2) {};

\draw[dotted] (1) -- (4,2);

\draw(4,2)--(6,2);
\node[draw,circle,fill=black] at (6,2) {};

\node[draw,circle,fill=black] at (-1,-1) {};

\draw (0.5,0.5) -- (-1,-1);

\end{tikzpicture}
\end{minipage}
\caption{We detach the three edges rooted at $v_0$. The created replicas are represented by squares and the star edges by dotted lines.}
\label{f:detach2}
\end{figure}
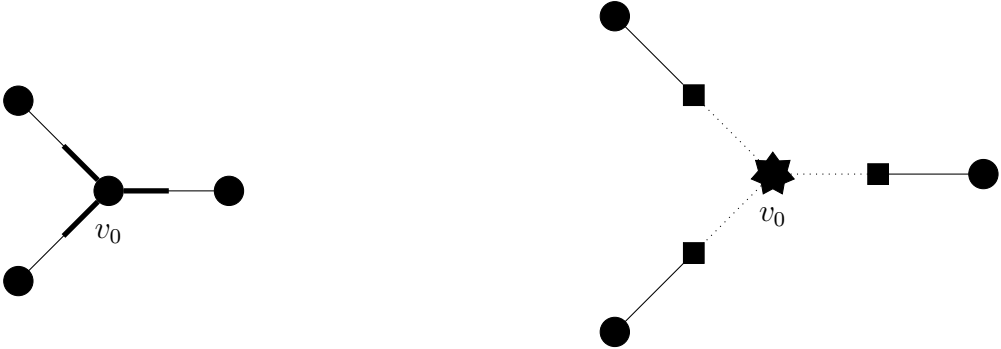

We will actually use the following version of the detachment lemma. Let $v_0$ be a vertex in $V$, and detach all directed edges rooted at $v_0$ (even those joining $\mathfrak v$ if they exist), creating as many replicas of $v_0$. Join the replicas to the original vertex $v_0$ (now being a star vertex) by edges (the star edges), see Figure \ref{f:detach2}. Let  $\G_{\{v_0\}}^*$ be the graph obtained, and $\n^*_{\{v_0\}}$ be the measure $\n$ of Section \ref{s:n}  on the graph $\G^*_{\{v_0\}}$ with $W=\{v_0\}$.

\begin{lemma}\label{l:detach2}
Let $x>0$. The trace of the loop soup $\L$ on $\G$ outside $\{v_0\}$  under $P(\cdot\, | \, \widehat \L(v_0)=x)$ has the same law as the trace of the loop soup on $\G^*_{\{v_0\}}$ outside the star edges under $\n^*_{\{v_0\}}$, conditioned on the event that all replica of $v_0$ have local time $x$. 
\end{lemma}

\noindent {\it Proof}. It is  again a consequence of Theorem \ref{t:markovn2} with now $W:=\{v_0\}$ there, and the remark following it. $\Box$

\section{Proof of Theorem \ref{t:stargraph}}

\label{s:stargraph}

We keep the setting of Section \ref{s:detachment}. Let $x>0$ and  $A$ be some event in the $\sigma$-field of the trace of $\L$ outside $\{v_0\}$ with $P(A\,|\, \widehat {\mathcal L}(v_0)=x)>0$ and $P(A)>0$. Lemma \ref{l:detach2} can be extended when replacing $P$ by $P(\cdot \,|\, A)$, and $\n^*_{\{v_0\}}$ by $\n^*_{\{v_0\}}(\cdot | A)$ (where $A$ is under $\n^*_{\{v_0\}}$ naturally seen as the corresponding event on the old edges of $\G^*_{\{v_0\}}$). It suffices to write $P(\cdot \,|\, A)$ as ${P(\cdot,\, A)\over P(A)}$ and use the lemma. Then, by continuity arguments, the lemma holds when replacing $P$ and $\n^*_{\{v_0\}}$ by $P(\cdot \,|\, \widehat \L(v)=x_v,\, v\in V')$ and  $\n^*_{\{v_0\}}(\cdot \,|\, \widehat \L(v)=x_v,\, v\in V')$ where $V'$ is a set of vertices disjoint from $\{v_0,\, \mathfrak v\}$. Let us draw some consequences of this remark. 

Let $W \subset V$ be a set of vertices (which, as $V$, does not contain $\mathfrak v$ ). Pick $v_0 \in W$ and condition the local time at $v_0$ on being $x$. One can detach all directed edges rooted at $v_0$ and work under $\n^*_{\{v_0\}}$  conditioned on the event that all replicas of $v_0$ have local time $x$. One then chooses another vertex $v_1 \in W$ (if it exists). We can check for example from Remark \ref{r:limit} that the law of the trace of the loop soup on $\G^*_{\{v_0\}}$ outside the star edges under $\n^*_{\{v_0\}}$ is absolutely continuous with respect to the law of the trace of the usual loop soup on the same graph. By the former paragraph, we deduce that Lemma \ref{l:detach2} also holds (for $v_1$ instead of $v_0$) when replacing $P$ by the probability  
$\n^*_{\{v_0\}}$ conditioned on replicas of $v_0$ having local time $x$, restricted to the $\sigma$-field of the trace of the loop soup outside the star edges. We can detach all directed edges rooted at $v_1$ and continue until we detached all directed edges rooted at  vertices of $W$. In the end, we transformed the graph $\G$ into the star graph of Theorem \ref{t:stargraph}. The proof of the theorem is complete. We restate it in the following form.

\begin{theorem}\label{t:detachment}
Let $W\subset V$ be a set of vertices of $\G$ ($W$ does not contain $\mathfrak v$). Define a graph $\G^*_W$ in the following way: 
\begin{itemize}
\item replace each vertex $v\in W$ by a star vertex that we still call $v$; 
\item add to each directed edge rooted at a star vertex $v$ a small line segment at the root, and call the endpoint of this line segment a replica of $v$. The added line segments are called star edges.
\end{itemize}
Let $\n^*_W$ be the measure $\n$ associated to the graph $\G^*_W$ and the set of star vertices (which we identify with $W$). Then, the trace of the loop soup $\L$ on $\G$ outside the set $W$ and conditioned on $\{\widehat \L(v)=x_v,\, v\in W\}$ has the same law as the trace of the loop soup on $\G^*_W$ outside the star edges under $\n^*_W$, conditioned on the event that all replica of $v\in W$ have local time $x_v$. 
\end{theorem}

 \section{Proof of Theorems \ref{t:markov} and \ref{t:markov2}}
 
 \label{s:markov}
 
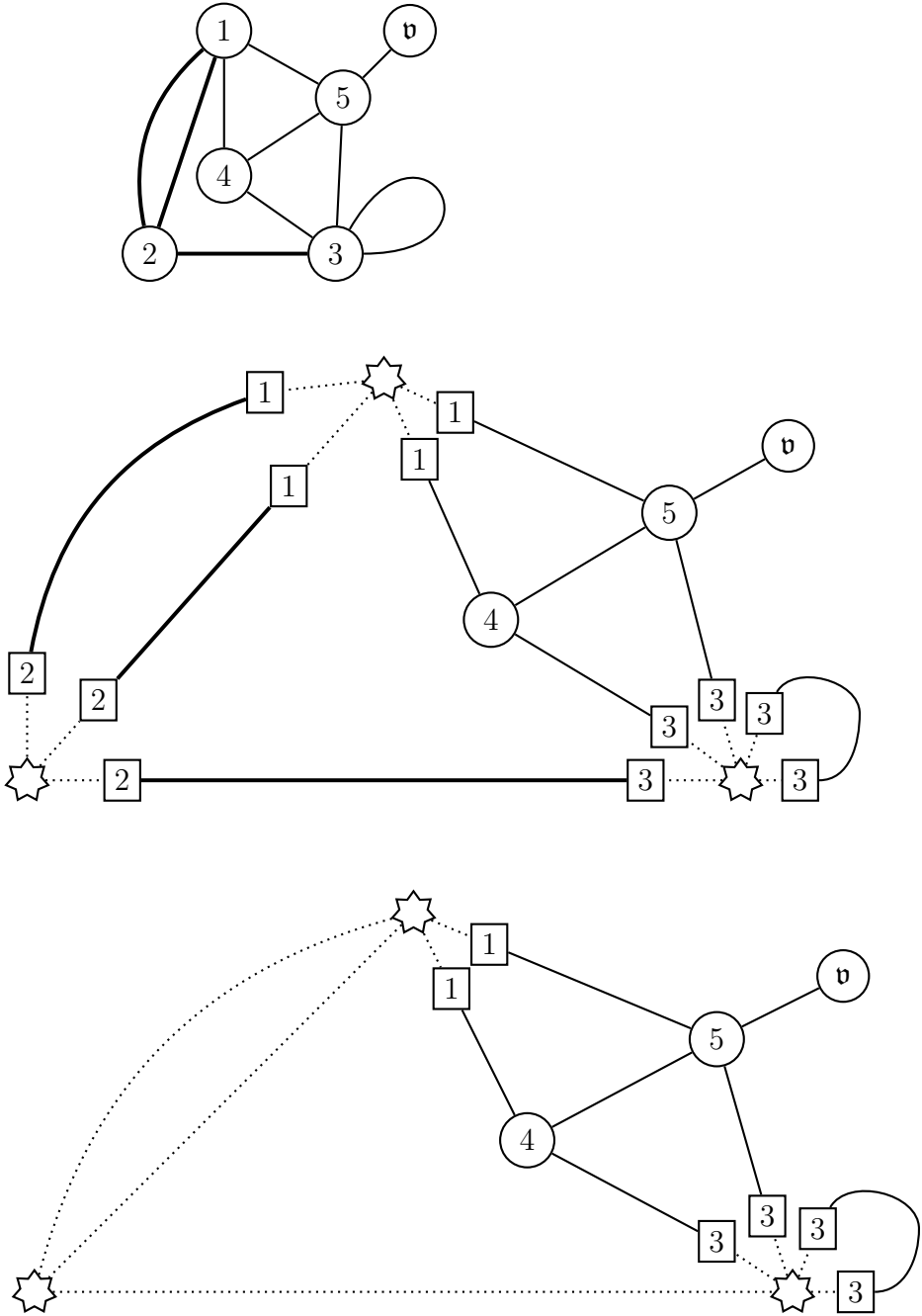
\begin{figure}
\begin{tikzpicture}[ thick,main node/.style={circle,draw},x=0.5cm,y=1.5cm]
\tikzset{every loop/.style={min distance=2cm,in=0,out=60}}
\clip(-8.5,-2.5) rectangle (6.5,0.5);
  \node[main node] (1) at (0,0){1};
  \node[main node] (2) at (-2,-2) {2};
  \node[main node] (3) at (3,-2) {3};
  \node[main node] (4) at (0.,-1.3) {4};
  \node[main node] (5) at (3.2,-0.6) {5};
  \node[main node] (0) at (5,0) {$\mathfrak v$};
  \path[every node/.style={font=\sffamily\small}]
    (1) edge[line width=1.5pt]  (2)
        edge [bend right,line width=1.5pt]  (2)
        edge (5)
    (2) edge[line width=1.5 pt]  (3)
     (1)   edge  (4)
     (4) edge (3)
     (4) edge (5)
     (5) edge (3)
     (3) edge[loop] node {} (3)
     (0) edge (5)
     ;
\end{tikzpicture}

\begin{tikzpicture}[ thick,main node/.style={rectangle,draw},x=1.6cm,y=1.8cm]
\tikzset{every loop/.style={min distance=2cm,in=0,out=60}}
\clip(-4,-3.5) rectangle (5,0.5);
  \node[draw,star,star points=7] (1) at (0,0) {};
  
      \node[main node] (11) at (-1,-0.1) {1};
      \node[main node] (12) at (-0.8,-0.8) {1};
      \node[main node] (13) at (0.3,-0.6) {1};
	  \node[main node] (14) at (0.6,-0.25) {1};
  \node[draw,star,star points=7] (2) at (-3,-3) {};
      \node[main node] (21) at (-3,-2.2) {2};
      \node[main node] (22) at (-2.4,-2.4) {2};
      \node[main node] (23) at (-2.2,-3) {2};

  \node[draw,star,star points=7] (3) at (3,-3) {};	
       \node[main node] (31) at (2.2,-3) {3};
      \node[main node] (32) at (2.4,-2.6) {3};
      \node[main node] (33) at (3.2,-2.5) {3};
      \node[main node] (34) at (3.5,-3) {3};
      \node[main node] (35) at (2.8,-2.4) {3};
      
  \node[circle,draw] (4) at (0.9,-1.8) {4};
  \node[circle,draw] (5) at (2.4,-1.0) {5};
  \node[circle,draw] (0) at (3.4,-0.5) {$\mathfrak v$};
   \path[every node/.style={font=\sffamily\small}]
      (1) edge[dotted] node [left] {} (11)
          edge[dotted] node [left] {} (12)
          edge[dotted] node [left] {} (13)
          edge[dotted] node [left] {} (14)
      (2) edge[dotted] node [left] {} (21)
          edge[dotted] node [left] {} (22)
          edge[dotted] node [left] {} (23)
      (3) edge[dotted] node [left] {} (31)
          edge[dotted] node [left] {} (32)
          edge[dotted] node [left] {} (33)
     	  edge[dotted] node [left] {} (34)
    	  edge[dotted] node [left] {} (35)
    (11) edge[bend right,line width=1.5pt]  (21)
    (12)   edge[line width=1.5pt]  (22)
    (13)   edge   (4)
    (4) edge (32)
    (4) edge (5)
    (5) edge (35) 
    (14) edge (5)
    (23) edge[line width=1.5pt]  (31)
    (0) edge (5)
    ;   
    
    \draw (33)   to [out=60,in=90] (4.,-2.5) to  [out=270,in=0] (34);
\end{tikzpicture}

\begin{tikzpicture}[ thick,main node/.style={rectangle,draw},x=1.7cm,y=1.7cm]
\tikzset{every loop/.style={min distance=2cm,in=0,out=60}}
\clip(-4,-4) rectangle (5,0.5);
  \node[draw,star,star points=7] (1) at (0,0) {};
  
      \node[main node] (13) at (0.3,-0.6) {1};
	  \node[main node] (14) at (0.6,-0.25) {1};
  \node[draw,star,star points=7] (2) at (-3,-3) {};

  \node[draw,star,star points=7] (3) at (3,-3) {};	
      \node[main node] (32) at (2.4,-2.6) {3};
      \node[main node] (33) at (3.2,-2.5) {3};
      \node[main node] (34) at (3.5,-3) {3};
      \node[main node] (35) at (2.8,-2.4) {3};
      
  \node[circle,draw] (4) at (0.9,-1.8) {4};
  \node[circle,draw] (5) at (2.4,-1.0) {5};
  \node[circle,draw] (0) at (3.4,-0.5) {$\mathfrak v$};
   \path[every node/.style={font=\sffamily\small}]
      (1) edge[dotted,bend right] node [left] {} (2)
          edge[dotted] node [left] {} (2)
          edge[dotted] node [left] {} (13)
          edge[dotted] node [left] {} (14)
      (2) 
          edge[dotted] node [left] {} (3)
      (3) 
          edge[dotted] node [left] {} (32)
          edge[dotted] node [left] {} (33)
     	  edge[dotted] node [left] {} (34)
    	  edge[dotted] node [left] {} (35)
    (13)   edge   (4)
    (4) edge (32)
    (4) edge (5)
    (5) edge (35) 
    (14) edge (5)
    (0) edge (5)
    ;   
    
    \draw (33)   to [out=60,in=90] (4.,-2.5) to  [out=270,in=0] (34);
\end{tikzpicture}

\caption{Proof of Theorem \ref{t:markov}. Edges of $\widetilde \G$ are in bold lines (top). $W=\{1,2,3\}$. The boundary of $\widetilde \G$ consists of $\{1,3\}$. We first use Theorem \ref{t:detachment} (middle). Since the local times are positive on the bold lines, we can remove the corresponding replica and work on the graph $\G^+$ (bottom), under the associated measure $\n^+_W$. Then we use Theorem \ref{t:markovn}.}
\label{f:markov2}
\end{figure}

  Figure \ref{f:markov2} illustrates the proof. Let $W$ be the set of vertices of $\widetilde \G$. Let $\widetilde {\mathcal E}$ be the event that the local times on the edges of $\widetilde \G$ are always positive. We use Theorem \ref{t:detachment}. The trace of $\L$ outside $W$ under $P$ conditionally on $\{\widehat \L(v)=x_v,\,\forall \, v\in W\}$ is distributed as the trace of the  loop soup  on $\G^*_W$ outside the star edges, under the measure $\n^*_W$ conditioned on the event that the local times at replica of $v$ are all equal to $x_v$ for all $v\in W$. We identify the old edges in $\G^*_W$ with the edges of $\G$. The event $\widetilde {\mathcal E}$ for the loop soup on $\G$ corresponds for the loop soup on $\G^*_W$  to the event that the local times fields on the edges of $\widetilde \G$ in $\G^*_W$ do not hit $0$. Let $\G^+_W$ be the graph $\G^*_W$ where we discard from the set of vertices the replica which are endpoints in $\G^*_W$ of an edge belonging to $\widetilde \G$. We let $\n_W^+$ be the measure $\n$ constructed in Section \ref{s:n} associated to the graph $\G^+_W$ (and we take for $W$ there the set of star vertices). By Definition \ref{def:n}, we see that $\n_W^*(\cdot\, |\, \widetilde {\mathcal E})=\n_W^+(\cdot)$. Theorem \ref{t:markov} is then a consequence of the Markov property for $\n_W^+$ stated in Theorem \ref{t:markovn}. When $\widetilde \G$ is connected, Theorem \ref{t:markov2} is a consequence of Theorem \ref{t:markovn2}. $\Box$

\section{Proof of Theorem \ref{t:cluster}}

\label{s:clusterproof}

Consider the graph $\G=(V,E)$ with added vertex $\mathfrak v$. Let $\G^*$ be the star graph of the introduction. It is the graph $\G^*_W$ of Theorem \ref{t:detachment} obtained by taking for $W$ the set of all vertices in $V$. We write $\n^*$ for the measure $\n$ associated.

When erasing the replicas, we will denote by $\G^+=(V,E^+)$ the obtained graph (we do not include in $E^+$ edges adjacent to $\mathfrak v$). We let ${\overline \L^+}$ be the discrete loop soup on the graph $\G^+$. We write $e^+$ for an edge of $\G^+$ (i.e. an edge between star vertices, or an edge between a star vertex and $\mathfrak v$). If $\overline L^+$ is a discrete loop configuration on $\G^+$, we let $n_{\overline L^+}(e^+)$ be the number of crossings by $\overline L^+$ of $e^+$. 

Under $\n^*$, we call cluster a maximal set of adjacent edges in $E^+$ with nonzero local time in their interior. For $\widetilde E$ a set of edges in $E^+$, we mean by admissible set of edges in $\widetilde E$ a subset $E'\subset \widetilde E$ such that every vertex $v\in V$ has an even number of adjacent incoming edges in $E'$ (a self-loop is counted twice).

\begin{proposition}\label{p:uniform}
Let $\widetilde E\subset E^+$. Under $\n^*$, conditionally on the edges of the clusters being $\widetilde E$, the set $\left\{e^+\in E^+\,:\, n_{\overline \L^+}(e^+)=1\right\}$ is  uniform among all admissible sets of edges in $\widetilde E$ and is independent of the occupation field.  
\end{proposition}

\noindent {\it Proof}. We consider the graph $\G^*$ under the measure $\n^*$.  The conditional measure $\n^*$ given  that the edges of $\widetilde E$ have positive local time fields in their interior,  is the measure $\n$  associated to the graph $\G^*$ where we erase the replicas on any edge of $\widetilde E$, with star vertices being the star vertices of $\G^*$, identified with $V$. Let $\widetilde \n^*$ be this measure. In our setting, $\partial V$ is the set of replicas which are not on the edges of $\widetilde E$. The measure $\n^*$ conditioned on the edges of the clusters being $\widetilde E$ is the measure $\widetilde \n^*$ conditioned on the event that the occupation fields on edges between vertices of $\partial V$ hit $0$.

Let $S$ be an admissible set of edges in $\widetilde E$ and $\mathcal S:=\left\{e^+\in E^+\,:\, n_{\overline {\mathcal L}^+}(e^+)=1\right\}$.   By \eqref{eq:na2}, we have
$$
\widetilde \n^*(\mathcal S=S,\, n(\{v,v'\})=0,\, \forall v\neq v'\in \partial V\, |\, \widehat {\mathcal L}(v)=x_v,\, v\in \partial V) = \frac{1}{Z}. 
$$

\noindent It  implies that under $\widetilde \n^*$, conditionally on $\{n(\{v,v'\})=0,\, \forall v\neq v'\in \partial V\}$ and on $\{\widehat {\mathcal L}(v)=x_v,\, v\in \partial V\}$, the set $\mathcal S$ is uniform among admissible sets of edges in $\widetilde E$. Under this conditioning, the occupation field on an edge between vertices of $\partial V$ is independent of ${\mathcal S}$: it is the occupation field coming from excursions from the endpoints conditioned on not crossing the edge together with a loop soup inside the edge. We deduce that the statement still holds when further conditioning on the event that the occupation fields between vertices of $\partial V$ hit $0$. What we showed is that the  set $\mathcal S$ is uniform among admissible sets of edges in $\widetilde E$ under  $\n^*$, conditionally on the edges of the cluster being $\widetilde E$, and is independent of $\widehat {\mathcal L}(v)$ for $v$ being the replicas on the edges which are not in $\widetilde E$. To conclude that it is independent of the occupation field, observe that, conditionally on ${\mathcal S}$,  on $\widehat \L(v),\, v\in \partial V$, and on the event that edges of the clusters are the edges of $\widetilde E$,   the occupation field on an edge of $\widetilde E$ is a ${\rm BESQ}^3$ bridge whether the edge is crossed or not. On an edge which is not in $\widetilde E$: the occupation field between a star vertex and a replica is a ${\rm BESQ}^3$ bridge, between two replicas it is a ${\rm BESQ}^1$ bridge conditioned on hitting $0$, and  between a replica and $\mathfrak v$, it is a ${\rm BESQ}^1$ bridge. We see that the law of the occupation field does not depend  on the value of ${\mathcal S}$. $\Box$

\bigskip

\begin{figure}[h!]
\begin{tikzpicture}[ thick,main node/.style={circle,draw},x=0.8cm,y=1.2cm]
\tikzset{every loop/.style={min distance=2cm,in=0,out=60}}
\clip(-6.5,-2.5) rectangle (5.5,0.5);
  \node[main node] (1) at (0,0){1};
  \node[main node] (2) at (-2,-2) {2};
  \node[main node] (3) at (3,-2) {3};
  \node[main node] (4) at (0.,-1.3) {4};
  \node[main node] (5) at (3.2,-0.6) {5};
  \node[main node] (0) at (5,0) {$\mathfrak v$};
  \path[every node/.style={font=\sffamily\small}]
    (1) edge[line width=1.5pt]  (2)
        edge [bend right]  (2)
        edge (5)
    (2) edge[line width=1.5pt]  (3)
     (1)   edge  (4)
     (4) edge (3)
     (4) edge[line width=1.5pt] (5)
     (5) edge (3)
     (3) edge[loop] node {} (3)
     (0) edge (5)
     ;
\end{tikzpicture}

\begin{tikzpicture}[ thick,main node/.style={rectangle,draw},x=1.4cm,y=1.8cm]
\tikzset{every loop/.style={min distance=2cm,in=0,out=60}}
\clip(-4,-3.5) rectangle (5,0.2);
  \node[draw,star,star points=7] (1) at (0,0) {};
  
      \node[main node] (11) at (-1,-0.1) {1};
      \node[main node] (12) at (-0.8,-0.8) {1};
      \node[main node] (13) at (0.1,-0.6) {1};
	  \node[main node] (14) at (0.6,-0.25) {1};
  \node[draw,star,star points=7] (2) at (-3,-3) {};
      \node[main node] (21) at (-3,-2.2) {2};
      \node[main node] (22) at (-2.4,-2.4) {2};
      \node[main node] (23) at (-2.2,-3) {2};

  \node[draw,star,star points=7] (3) at (3,-3) {};	
       \node[main node] (31) at (2.2,-3) {3};
      \node[main node] (32) at (2.2,-2.6) {3};
      \node[main node] (33) at (3.2,-2.5) {3};
      \node[main node] (34) at (3.5,-3) {3};
      \node[main node] (35) at (2.8,-2.4) {3};
     
   \node[draw,star,star points=7] (4) at (0.2,-1.9) {};  
  	\node[main node] (41) at (0.2,-1.4) {4};
	\node[main node] (42) at (0.8,-1.65) {4};
	\node[main node] (43) at (0.8,-2.1) {4};	
   \node[draw,star,star points=7] (5) at (2.4,-1.0) {};  
  	\node[main node] (51) at (1.7,-0.7) {5};
	\node[main node] (52) at (1.9,-1.15) {5};
	\node[main node] (53) at (2.55,-1.5) {5};	\node[main node] (54) at (2.9,-0.7) {5};
  \node[circle,draw] (0) at (4.,0.) {$\mathfrak v$};
   \path[every node/.style={font=\sffamily\small}]
      (1) edge[dotted] node [left] {} (11)
          edge[dotted] node [left] {} (12)
          edge[dotted] node [left] {} (13)
          edge[dotted] node [left] {} (14)
      (2) edge[dotted] node [left] {} (21)
          edge[dotted] node [left] {} (22)
          edge[dotted] node [left] {} (23)
      (3) edge[dotted] node [left] {} (31)
          edge[dotted] node [left] {} (32)
          edge[dotted] node [left] {} (33)
     	  edge[dotted] node [left] {} (34)
    	  edge[dotted] node [left] {} (35)
	  (4) edge[dotted] node [left] {} (41)
          edge[dotted] node [left] {} (42)
          edge[dotted] node [left] {} (43)
      (5) edge[dotted] node [left] {} (51)
          edge[dotted] node [left] {} (52)
          edge[dotted] node [left] {} (53)  
		  edge[dotted] node [left] {} (54)
    (11) edge[bend right]   (21)         
    (12)   edge[line width=1.5pt]  (22)
    (13)   edge   (41)
    (43) edge (32)
    (42) edge[line width=1.5pt] (52)
    (53) edge (35) 
    (14) edge (51)
    (23) edge[line width=1.5pt]  (31)
    (0) edge (54)
    ;   
    
    \draw (33)   to [out=60,in=90] (4.,-2.5) to  [out=270,in=0] (34);
\end{tikzpicture}

\begin{tikzpicture}[ thick,main node/.style={rectangle,draw},x=1.6cm,y=1.6cm]
\tikzset{every loop/.style={min distance=2cm,in=0,out=60}}
\clip(-4,-3.5) rectangle (5,0.5);
  \node[draw,star,star points=7] (1) at (0,0) {};
  
      \node[main node] (11) at (-1,-0.1) {1};
      \node[main node] (13) at (0.1,-0.6) {1};
	  \node[main node] (14) at (0.6,-0.25) {1};
  \node[draw,star,star points=7] (2) at (-3,-3) {};
      \node[main node] (21) at (-3,-2.2) {2};

  \node[draw,star,star points=7] (3) at (3,-3) {};	
      \node[main node] (32) at (2.2,-2.6) {3};
      \node[main node] (33) at (3.2,-2.5) {3};
      \node[main node] (34) at (3.5,-3) {3};
      \node[main node] (35) at (2.8,-2.4) {3};
     
   \node[draw,star,star points=7] (4) at (0.2,-1.9) {};  
  	\node[main node] (41) at (0.2,-1.4) {4};
	\node[main node] (43) at (0.8,-2.1) {4};	
   \node[draw,star,star points=7] (5) at (2.4,-1.0) {};  
  	\node[main node] (51) at (1.7,-0.7) {5};
	\node[main node] (53) at (2.55,-1.5) {5};	\node[main node] (54) at (2.9,-0.7) {5};
  \node[circle,draw] (0) at (4.,0.) {$\mathfrak v$};
   \path[every node/.style={font=\sffamily\small}]
      (1) edge[dotted] node [left] {} (11)
          edge[dotted] node [left] {} (2)
          edge[dotted] node [left] {} (13)
          edge[dotted] node [left] {} (14)
      (2) edge[dotted] node [left] {} (21)
          edge[dotted] node [left] {} (3)
      (3) 
          edge[dotted] node [left] {} (32)
          edge[dotted] node [left] {} (33)
     	  edge[dotted] node [left] {} (34)
    	  edge[dotted] node [left] {} (35)
	  (4) edge[dotted] node [left] {} (41)
          edge[dotted] node [left] {} (5)
          edge[dotted] node [left] {} (43)
      (5) edge[dotted] node [left] {} (51)
          edge[dotted] node [left] {} (53)  
		  edge[dotted] node [left] {} (54)
    (11) edge[bend right]   (21)         
    (13)   edge   (41)
    (43) edge (32)
    (53) edge (35) 
    (14) edge (51)
    (0) edge (54)
    ;   
    
    \draw (33)   to [out=60,in=90] (4.,-2.5) to  [out=270,in=0] (34);
\end{tikzpicture}

\caption{ Proof of Theorem \ref{t:cluster}. Graph $\G$ (top). Edges in bold lines have positive occupation fields (there are 2 clusters). Graph $\G^*$ (middle). We can erase replicas at the ends of the bold edges (bottom). In the new star graph, choose uniformly an admissible set of crossed edges inside the clusters (in the example, it is necessarily empty).  Crossed edges  correspond to edges with an odd number of crossings in $\G$. Use the the definition of $\n$ in Section \ref{s:n} to give the law of the trace of the loop soup on the edges.}
\label{f:clusterproof}
\end{figure}
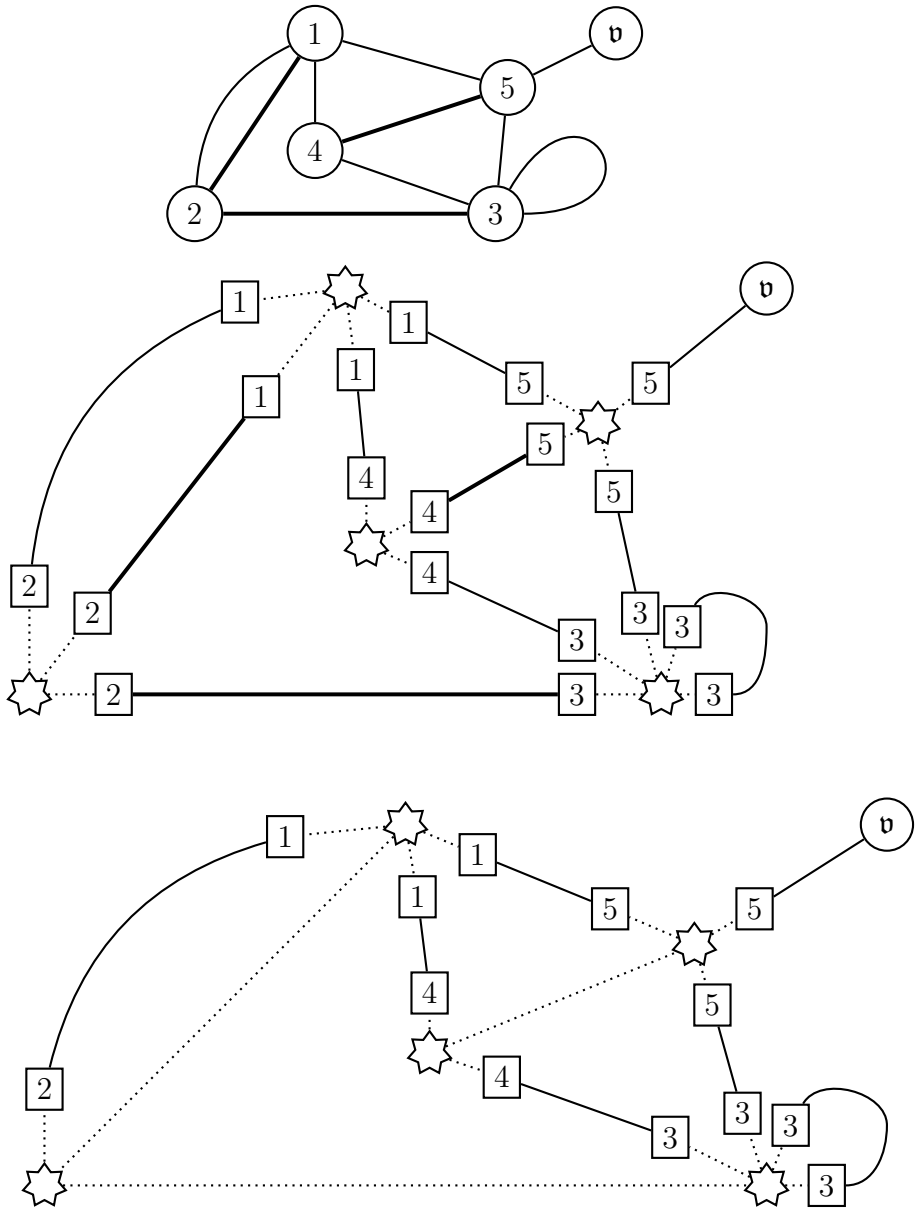

{\it Proof of Theorem \ref{t:cluster}}.
 By  Theorem \ref{t:detachment}, the trace of the loop soup $\L$ outside the set of vertices $V$ conditioned on having local time $x_v$ at vertex $v$ is equal in law to the trace of the  loop soup  on $\G^*$ outside the star edges under $\n^*$ conditioned on the event that all replicas of $v$ have local time $x_v$. Notice that $\alpha_e = n_{\overline \L^+}(e^+)$ where $e^+$ is the extended edge in the graph $\G^+$ which corresponds to $e$.  Assertion (1) is then a consequence of Proposition \ref{p:uniform}. Let us prove (2). Conditionally on $\overline \L^+$, the traces on the edges $e^+$ are independent. Recall the one-dimensional loop soups introduced in Section \ref{s:loopsoupcond}.  If $n_{\overline \L^+}(e^+)=1$, the loop soup on $e^+$ is distributed as ${\mathcal C}(\rho(e^+),0,0)$.  If  $n_{\overline \L^+}(e^+)=0$, the trace of the loop soup on the edge $e^+$ is that of ${\mathcal B}(\rho(e^+),0,0)$, conditioned on having positive local time along the two star edges lying at the extremities of $e$ (we can make sense of it by taking limits  $\ell_1,\ell_2 \searrow 0$ of ${\mathcal B}(\rho(e^+),\ell_1,\ell_2)$ conditioned on the same event).  Assertion (2) is then a consequence of the spatial Markov property of one-dimensional loop soups recalled in Section \ref{s:loopsoupcond}. See Figure \ref{f:clusterproof} for an illustration of the proof. $\Box$

\section{Relation with Le Jan's isomorphism theorem}

\label{s:lejan}
 
We use the following exploration process.   We start at time $T_0:=0$ at some vertex ${\mathfrak e}_0\neq \mathfrak v$. If $v_k = {\mathfrak e}_{T_k}$ is defined for some $k\ge 0$ and is adjacent to some unexplored edge $e$, we explore $e$ at speed $1$, and we stop at the first point with local time $0$ that we meet on our way, or at the endpoint of $e$ in the case local times are positive throughout the edge. We let $T_{k+1}$ be such that $T_{k+1}-T_k$ represents the distance covered from $v_k$ to the stopping point, and define ${\mathfrak e}_{t},\, t\in (T_k,T_{k+1}]$ as the point at distance $t-T_k$ of $v_k$. If we met a point with local time $0$ on our way, we artificially cut the edge at that point, and identify the root of the unexplored part of the edge with $\mathfrak v$\footnote{This operation will not change the conditional distribution of the trace of the loop soup in the unexplored region given the occupation field inside the explored one. It comes from the simple Markov property stated in Theorem \ref{t:markov2} and a continuity argument.}. In the case where all adjacent edges at $v_k$ have been explored, we choose  a vertex $v_\ell$ for some $\ell < k$ with  some unexplored adjacent edge, which we explore as before, thus defining some $T_{k+1}$ and ${\mathfrak e}_t,\, t\in (T_k,T_{k+1}]$. The exploration of the cluster terminates when such an $\ell$ does not exist. We let $\zeta$ be the time at which the exploration terminates.

We define for $t\in [0,\zeta]$, ${\mathscr X}_t:=\sqrt{\widehat \L({\mathfrak e}_t)}$, i.e. the square root of the local time at the location of the explorer at time $t$. Call $\G_t$ the unexplored part of the graph at time $t$, i.e.,   
$$
\G_t:=  \G\backslash \{{\mathfrak e}_s,\, s\in [0,t]\} . 
 $$

 \noindent For $v\neq v' \in \{{\mathfrak e}_s,\, s\in [0,t]\}$, let $H_t(v,v')$ be the mass under the excursion measure at $v$ of excursions inside $\G_t$ which hit back $\{ {\mathfrak e}_s,\, s\in [0,t]\}$ at $v'$.  In particular, $H_t(v,v')=0$ if $v$ or $v'$ is not on the boundary of $\G_t$. For $t\ge 0$, let $\mathcal F_t$ be the $\sigma$-field generated by the random set $\{\mathfrak e_s,\, s\le t\land \zeta\}$  and the occupation field inside $\{\mathfrak e_s,\, s\le t\land \zeta\}$. For convenience, we let $(\gamma_t,\, t\ge 0)$ be a Brownian motion independent of everything else, and we add in ${\mathcal F}_t$ the $\sigma$-field generated by $(\gamma_s,\, s\le t)$.

\begin{theorem}
There exists  a one-dimensional standard ${\mathcal F}_t$-Brownian motion $B$ such that on each $(T_k,T_{k+1})$ before time $\zeta$, the process ${\mathscr X}$ satisfies
\begin{equation}\label{eq:sde}
\d {\mathscr X}_t = \d B_t  +  2 \sum_{s <t} H_t({\mathfrak e}_t,{\mathfrak e}_s) ({\mathscr X}_s-{\mathscr X}_t)  \, \d t.
\end{equation}
\end{theorem}
\noindent {\it Proof}. For $t\in [0,\zeta]$, let ${\mathscr L}_t := \widehat \L({\mathfrak e}_t)$. Consider the trace of $\L$ on $\G_t$ at a time $t$. It is composed of excursions from the boundary to itself,  and of loops which stay inside $\G_t$.  Let $N_t(\mathfrak e_s)$ be the number of excursions inside $\G_t$ between $\mathfrak e_t$ 
and $\mathfrak e_s$, and $N_t :=\sum_{s<t} N_t(\mathfrak e_s)$. Let $H_t := \sum_{s<t} H_t(\mathfrak e_t,\mathfrak e_s)$.  Consider a time interval $(T_k,T_{k+1})$, related to some edge $e$. Suppose for simplicity that it is not a self-loop, which we can suppose by creating an artificial vertex in the middle of the edge. Call $e_1$ and $e_2$ the endpoints of $e$ and suppose that we are exploring the edge from $e_1$ to $e_2$.  One-dimensional Ray--Knight theorems show that, on $(T_k,T_{k+1})$, ${\mathscr L}$  is solution of 
\begin{equation}\label{eq:sdeproof}
\d {\mathscr L}_t = 2 \sqrt{{\mathscr L}_t} \d \widetilde B_t -4 {\mathscr L}_t H_t \d t + (2 N_t +1) \d t
\end{equation}

\noindent for some one-dimensional $\widetilde {\mathcal F}_t$-Brownian motion $\widetilde B$ where $\widetilde {\mathcal F}_t$ is generated by the random set $\{\mathfrak e_s,\, s<t\land \zeta\}$, the occupation field in $\{\mathfrak e_s,\, s<t\land \zeta\}$ and $(N_s,\, T_k< s\le t\land \zeta)$. Let us say a few words on why this SDE holds. Each excursion from ${\mathfrak e}_t$ to another point (there are a number $N_t$ of them) gives a drift $2\, \d t$ to the local time  (they look like Bessel(3) excursions on the vicinity of ${\mathfrak e}_t$). The Brownian loop soup in $\G_t$ gives a drift $\d t$ (it locally looks like a Brownian loop soup on ${\mathbb R}_+$). The term $2 \sqrt{{\mathscr L}_t} \d \widetilde B_t$ is the second Ray-Knight theorem for the local times of a one-dimensional Brownian motion from which we need to remove excursions from $\mathfrak e_t$ to itself which hit $\{\mathfrak e_s,\, s<t\}$, which makes appear the term $-4 {\mathscr L}_t H_t \d t$, the local time process of such an excursion being locally a ${\rm BESQ}^4$ process. To rigorously prove the SDE, one can first look at all excursions away from the boundary of $\G_{T_k}$ which intersect the interior of the edge $e$. An excursion joining points of $\partial \G_{T_k}$ distinct of $e_1$ which enters the edge $e$ without hitting $e_1$ will be associated to a pair of excursions to a point $u$ on the edge, the point $u$ being the tip of the excursion inside the edge $e$. One conditions on such points $u$ and on $\lim_{t\searrow T_k} N_t$ (it is the number of excursions in $\G_{T_k}$ from $e_1$ to another point of $\partial \G_{T_k}$, or to itself when the excursion crosses $e$ and returns to $e_1$ via another edge). One can then couple the trace of the loop soup on $e$ with one-dimensional Brownian paths: the traces of excursions in $\G_{T_k}$ from $e_1$ to itself which leave and return via $e$ are one-dimensional Brownian paths which, in the case they hit $e_2$, give an exponentially distributed local time with parameter $H_{T_k}(e_2)$ (the mass of excursions at $e_2$ hitting $\partial \G_{T_k}$) at $e_2$ before coming back to $e_1$.  The trace on $e$ of loops inside $\G_{T_k}$ visiting $e_2$ form a collection of Brownian excursions away from $e_2$ up to a local time which is gamma($\frac12, H_{T_k}(e_2)$) distributed. Loops inside $e$ form a loop soup.  The collection of these trajectories is distributed as the trace on $[0,\rho(e)]$ of a one-dimensional loop soup conditioned on having local time $\widehat {\mathcal L}(e_1)$ at $0$ and local time $0$ at $\frac{1}{2H_{T_k}(e_1)}$, where $H_{T_k}(e_1):=\lim_{t\searrow T_k} H_t$. Its local time process is a ${\rm BESQ}^1$ bridge from $\widehat {\mathcal L}(e_1)$ to $0$ of duration $\frac{1}{2H_{T_k}(e_1)}$. An excursion starting from a point $u$  inside the edge behaves as a Bessel(3) process stopped when hitting $e_2$, then spends an exponentially distributed local time there (the parameter depends on the location of $u$) and leaves. It is the trace of a (translated) one-dimensional Bessel(3) process stopped when hitting  $\frac{1}{2H_{T_k}(e_1)}$. Its local time process is a ${\rm BESQ}^2$ bridge from $0$ to $0$ with adequate duration. We then use the additivity of ${\rm BESQ}$ bridges (see Chapter XI of \cite{revuz-yor} for the SDE of ${\rm BESQ}$ bridges). 
 
Taking the square root in \eqref{eq:sdeproof} yields 
$$
\d {\mathscr X}_t = \d \widetilde B_t - 2{\mathscr X}_tH_t  \, \d t + {N_t \over {\mathscr X}_t}\, \d t.
$$

\noindent We use the following lemma:
\begin{lemma}
In some completed filtration $(\widetilde {\mathcal F}_t,\, t\ge 0)$, consider a continuous semi-martingale of the form $X_t = \widetilde B_t + \int_0^t \widetilde A_s \d s$ where $\widetilde B$ is a standard $\widetilde {\mathcal F}_t$-Brownian motion and $\widetilde A$ is a progressive stochastic  process such that $s\to E[|\widetilde A_s|]$ is finite and locally integrable.   Suppose that $X$ is adapted with respect to a smaller completed filtration $ ({\mathcal F}_t,\, t\ge 0)$.  Let $A_s:=E[\widetilde A_s \,|\, {\mathcal F}_s]$ and suppose that it defines a progressive stochastic process with respect to $ ({\mathcal F}_t,\, t\ge 0)$. 
Then $X_t= B_t +  \int_0^t A_s \d s$, where $B$ is a ${\mathcal F}_t$-Brownian motion.
\end{lemma}
 \noindent {\it Proof}. The integral $\int_0^t A_s \d s$ is almost surely well-defined since $s\mapsto A_s$ is progressive and $\int_0^t E[ |A_s|] \d s \le \int_0^t E[ |\widetilde A_s|]<\infty$. It is also measurable with respect to $\mathcal F_t$, again because  $s\mapsto A_s$ is progressive. Let $B_t:= X_t - \int_0^t A_s \d s= \widetilde B_t + \int_0^t (\widetilde A_s-A_s)\d s$. Then $B$ is continuous and adapted with respect to $ ({\mathcal F}_t,\, t\ge 0)$. Moreover, it is a ${\mathcal F}_t$-martingale. Indeed, for $s\le t$, 
$E[ B_t \, |\, {\mathcal F}_s]  = E[ \widetilde B_t \, |\, {\mathcal F}_s] + \int_0^t E[ \widetilde A_u - A_u \, |\, {\mathcal F}_s] \d u $. We have $E[ \widetilde B_t \, |\, \widetilde {\mathcal F}_s] = \widetilde B_s$, hence $E[ \widetilde B_t \, |\, {\mathcal F}_s] = E[ \widetilde B_s \, |\, {\mathcal F}_s]$ since ${\mathcal F}_s \subset \widetilde {\mathcal F}_s$.  When $u\ge s$, $E[ \widetilde A_u  \, |\, {\mathcal F}_s] = E[ A_u  \, |\, {\mathcal F}_s]$ since ${\mathcal F}_s \subset {\mathcal F}_u$. We get that $E[ B_t \, |\, {\mathcal F}_s]= E[ \widetilde B_s \, |\, {\mathcal F}_s] + \int_0^s E[ \widetilde A_u - A_u \, |\, {\mathcal F}_s] \d u = E[B_s \, |\, {\mathcal F}_s]=B_s $. Finally, the quadratic variation process of $B$ is  $t$. By L\'evy's characterization theorem, $ B$ is a ${\mathcal F}_t$-Brownian motion. $\Box$
 
 \bigskip

We use the lemma on $(T_k,T_{k+1})$. The process $\mathscr X$ is adapted to the filtration $({\mathcal F}_t,\, t\ge 0)$. By Theorems \ref{t:markov} and \ref{t:markov2}, $N_t(\mathfrak e_s)$ conditionally on ${\mathcal F}_t$ is  Poisson distributed with parameter $2H_t(\mathfrak e_t,\mathfrak e_s){\mathscr X}_t{\mathscr X}_s$.  We deduce that the conditional expectation of $N_t$ is $\sum_{s<t} 2H_t(\mathfrak e_t,\mathfrak e_s){\mathscr X}_t{\mathscr X}_s$ then the theorem (we define $B_t:=\gamma_t$ for $t>\zeta$ so that it is defined on $\r_+$). $\Box$

\bigskip
 
\begin{corollary}
The local times $(\widehat {\mathcal L}(v),\, v\in V)$ are distributed as $(\frac12 \phi_v^2,\, v\in V)$ where $(\phi_v,\,v\in V)$ is a Gaussian free field on $\G$ with boundary value $0$ at $\mathfrak v$.
\end{corollary} 
 
 \noindent {\it Proof}. We identify equation \eqref{eq:sde} with equation 5 in Lemma 2 of Lupu and Werner \cite{lw16}.  We note that in our case, we actually explore only one cluster at a time, so we need to iterate the exploration process each time we terminate the exploration of a cluster. $\Box$

\section{Reconstructing the loops}
\label{s:reconstruct}

Theorem \ref{t:cluster} together with Section \ref{s:loopsoupcond} gives the conditional law of the loop soup on the edges given the occupation field. This section completes the description. We will show how to glue excursions on adjacent edges to recreate the loops of the loop soup, following Werner \cite{wernermarkov}.  We condition on the occupation field $\widehat {\mathcal L}$ on the edges of the metric graph.

The clusters are measurable with respect to $\widehat {\mathcal L}$. We choose uniformly an admissible configuration $(\alpha_e,\, e\in E)$ of $0$ and $1$ as in Theorem \ref{t:cluster}. An edge with $\alpha_e=1$ will have an odd number  of crossings, an edge with $\alpha_e=0$ will have an even number of crossings.  

The traces of the loop soup on the edges are then independent. On an edge $e$ with $\alpha_e=0$: we know that the trace of the loop soup is that of a one-dimensional Brownian loop soup conditioned on the same occupation field. As explained in Section \ref{s:loopsoupcond}, we can reconstruct it from a Jacobi($1,0$) flow. Let
$$
T:=\int_0^{\rho(e)} \frac{\d s}{\widehat {\mathcal L}(e_s)}
$$ 

\noindent where $\widehat {\mathcal L}(e_s)$ denotes the local time at the point $e_s$, and $(e_s,\, s\in [0,\rho(e)])$ are the points on the edge $e$ which are naturally indexed  by the interval $[0,\rho(e)]$. We then consider the Jacobi($1,0$) flow from time $0$ to time $T$. We can define its contour function: it is the process whose local time flow is the Jacobi flow (it is the process $\widetilde X^{(\rho)}$ of Section \ref{s:loopsoupcond}). Inverting \eqref{def:Xtilde}, we can trace the crossings of the edge $e$, the loops which stay in $e$ and the excursions away from the endpoints of $e$ which do not cross $e$, using  the contour function. 

For example, the number of crossings corresponds in the Jacobi flow to twice the number of points $\widetilde x \in (0,1]$ at time $0$  such that $\widetilde Y_{0,T}(\widetilde x-) < \widetilde Y_{0,T}(\widetilde x)$ (in the notation of Section \ref{s:loopsoupcond}). It is also the number of particles in a Fleming-Viot process with immigration  which have a progeny at time $T$. From the well-known connection between Kingman coalescent and Fleming--Viot processes (\cite{eth},\cite{foucart11},\cite{fh}), one can describe it also as the number of blocks in a Kingman coalescent with killing. This killed Kingman coalescent is described as follows: the process has values in the set of subpartitions of $\{1,2,\ldots\}$.  The sets of the subpartition are called blocks.  At time $0$,  the partition is composed of the singletons $ \{1\}, \{2\}, \ldots$ Then, blocks merge at rate $4$ and disappear  at rate $1$. The number of blocks of this process at time $T$ is distributed as half the number of crossings of $e$.

Similarly, on an edge with an odd number of crossings, we can reconstruct the loop soup with a Jacobi($1,2$) flow. The number of crossings minus $1$, and divided by $2$, is given by a killed Kingman coalescent, where the merging rate is $4$ and the killing rate is now $3$.

Once we know the crossings, locally at every vertex, we glue pairs of crossings arriving at this vertex uniformly, see \cite{wernermarkov}. It will recreate the loops of the discrete loop configuration $\overline{\mathcal L}$. It remains to glue the excursions in the edges which are not crossings. Consider a vertex $v$. Say that the total number of crossings of the directed edges rooted at $v$ is $2k_v$. Excursions away from $v$ which do not cross any adjacent edge are naturally labelled by the local time process at $v$ (say divided by $\widehat {\mathcal L}(v)$ so that its range is $[0,1]$): for example one can take their label in the Jacobi flows (i.e. the associated local time in the process $\widetilde X^{(\rho)}$ in Section \ref{s:loopsoupcond}). We take a beta($\frac12,k_v$)-distributed random variable $\beta$ which represents the accumulated local time of all the loops visiting $v$ which did not cross any adjacent edge. We divide the interval $[0,\beta]$ with a Poisson--Dirichlet($0,\frac12$) partition. For each interval of the partition, we concatenate the excursions whose labels belong to the interval, in increasing order of their label. We will have then reconstructed all the loops which did not cross any edge. We then divide the interval $[\beta,1]$  into $k_v$ intervals throwing $k_v-1$ i.i.d. uniform random variables. Each pair of crossings that we earlier glued together corresponds to an interval. When a loop arrives at the vertex $v$ through some crossing, concatenate it with the excursions with labels belonging to the corresponding interval. In this way, we will have reconstructed all loops.


\end{document}